\documentclass[12pt]{article}
\usepackage{amsfonts}
\usepackage{amsmath}
\usepackage{epsfig}

\title{The Flapping Birds in the Pentagram Zoo}
\author{Richard Evan Schwartz \thanks{Supported by N.S.F. Grant DMS-21082802}}

\newtheorem{theorem}{Theorem}[section]

\newtheorem{lemma}[theorem]{Lemma}

\newtheorem{corollary}[theorem]{Corollary}
\newtheorem{conjecture}[theorem]{Conjecture}
\newtheorem{question}[theorem]{Question}

\def\startproof{{\bf {\medskip}{\noindent}Proof: }}

\def\endproof{$\spadesuit$  \newline}

\def\C{\mbox{\boldmath{$C$}}}%
\def\P{\mbox{\boldmath{$P$}}}%
\def\R{\mbox{\boldmath{$R$}}}%
\def\Z{\mbox{\boldmath{$Z$}}}%

\begin{document}
\maketitle

\begin{abstract}
  We study the $(k+1,k)$ diagonal map for
  $k=2,3,4,...$.  We call this map $\Delta_k$.
  The map $\Delta_1$ is the pentagram map and
  $\Delta_k$ is a generalization.
  $\Delta_k$ does not preserve
  convexity, but we prove that $\Delta_k$
  preserves a subset $B_k$ of certain star-shaped
  polygons which we call $k$-{\it birds\/}.  The action of
    $\Delta_k$ on $B_k$ seems similar to the action
  of $\Delta_1$ on the space of convex polygons.  We
  show that some classic geometric results
  about $\Delta_1$ generalize to this setting.
\end{abstract}

\section{Introduction}

\subsection{Context}

When you visit the pentagram zoo you should certainly make
the pentagram map itself your first stop.  This old
and venerated animal has been around since the place opened
up and it is very friendly towards children.  When defined
on convex pentagons, this map has a very long history. See
e.g. \cite{MOT}.  In modern times \cite{SCH1}, the pentagram
is defined and studied much more generally. The
easiest case to explain is the action on convex $n$-gons.
One starts with
a convex $n$-gon $P$, for $n \geq 5$, and then forms a new
convex $n$-gon
$P'$ by intersecting the consecutive diagonals, as shown Figure 1.1 below.

The magic starts when you iterate the map.  One of the first things
I proved in \cite{SCH1}
about the pentagram map is the successive iterates
shrink to a point.  Many years later, M. Glick  \cite{GLICK1}
proved that this limit point is an algebraic function of the
vertices, and indeed found a formula for it.
See also \cite{IZOS2} and \cite{AI}.

\begin{center}
\resizebox{!}{2in}{\includegraphics{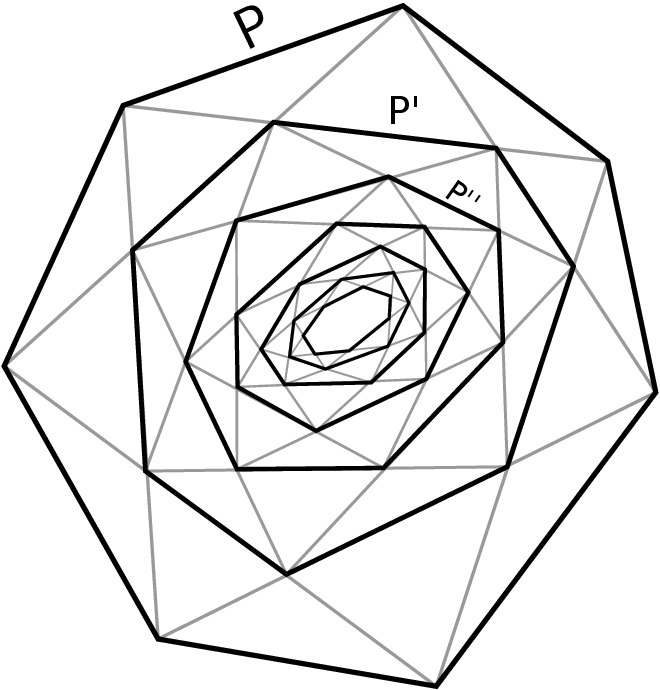}}
\newline
Figure 1.1: The pentagram map iterated on a convex $7$-gon $P$.
\end{center}

Forgetting about convexity, the pentagram map
is generically defined on polygons in the projective plane over any field
except for $\Z/2$.  In all cases, the pentagram map commutes with
projective transformations and thereby defines a birational map on the
space of $n$-gons modulo projective transformations.  The action
on this moduli space has a beautiful structure.
As shown in \cite{OST1} \cite{OST2}, and independently in
\cite{SOL}, the pentagram map
is a discrete completely integrable system when the ground
field is the reals.  (\cite{SOL} also treats the complex case.)
Recently, M. Weinreich
\cite{WEIN} generalized the integrability result, to a large extent,
to fields of positive characteristic.

The pentagram map has many generalizations.
See for example \cite{GSTV},   \cite{MF}, \cite{OV}, \cite{KS1}, \cite{KS2},  \cite{IK}.
The paper \cite{GSTV} has the first general complete integrability
result.
The authors prove the complete integrability of the
$(k,1)$ diagonal maps, i.e. the maps
obtained by intersecting successive $k$-diagonals.
Figure 1.3 below shows the $(3,1)$ diagonal map.
(Technically, \cite{GSTV} concentrates on what happens when these maps
act on so-called {\it corregated polygons\/} in higher dimensional
Euclidean spaces.)
The paper \cite{IK} proves an integrability result
for a very wide class of generalizations, including
the ones we study below. (Technically, for the maps we
consider here, the result in \cite{IK} does not 
establish the algebraic independence of invariants needed for
complete integrability.)
The pentagram map and its many generalizations are
related to a number of
topics: alternating sign matrices \cite{SCH2}, dimers \cite{GK}, cluster algebras \cite{GLICK2},
the KdV hierarchy \cite{MB1}, \cite{MB2},  spin networks \cite{GSTV}, Poisson Lie groups \cite{IZOS1},
Lax pairs \cite{SOL}, \cite{KS1}, \cite{KS2}, \cite{IK}, \cite{IZOS1},
and so forth. The zoo has many cages and sometimes you have to
get up on a tall ladder to see inside them.

\begin{center}
\resizebox{!}{2.2in}{\includegraphics{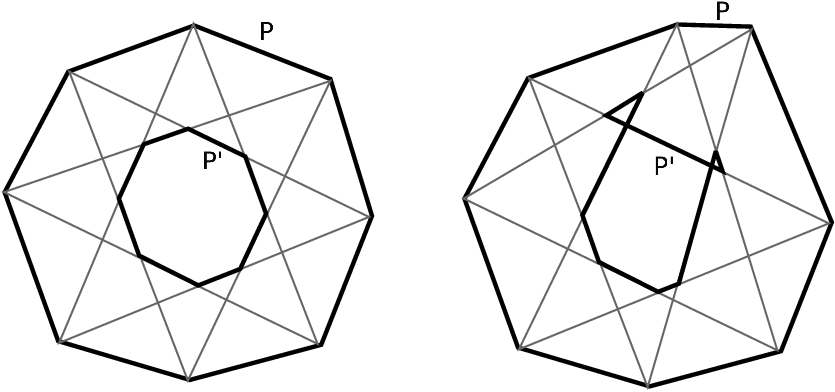}}
\newline
Figure 1.2: The $(3,1)$-diagonal map acting on $8$-gons.
\end{center}

The algebraic side of the pentagram zoo is extremely well
developed, but
the {\it geometric side\/} is hardly
developed at all. In spite of all the algebraic
results, we don't really know, geometrically speaking,
much about what the pentagram map and
its relatives really do to polygons.

Geometrically speaking, there seems to be a dichotomy
between convexity and non-convexity. The generic
pentagram orbit of a projective equivalance class of a
convex polygon lies on a smooth torus, and
you can make very nice animations.  What you will see,
if you tune the power of the map and pick suitable
representatives of the projective classes, is a convex
polygon sloshing around as if it were
moving through water waves.
If you try the pentagram map on a non-convex polygon,
you see a crazy erratic picture no matter how you
try to normalize the images.  The situation is even worse
for the other maps in the pentagram zoo, because these
generally do not preserve convexity.  Figure 1.2
shows how the $(3,1)$-diagonal map does not
necessarily preserve convexity, for instance.
See \cite{SCH3}, \cite{SCH4} for more details.

If you want to look at pentagram map generalizations, you have
to abandon convexity.  However, in this paper, I
will show that sometimes there are geometrically appealing replacements.
The context for these replacements is the
$(k+1,k)$-diagonal map, which I call $\Delta_k$,
acting on what I call $k$-{\it birds\/}.
$\Delta_k$ starts with the polygon $P$
and intersects the $(k+1)$-diagonals which differ
by $k$ clicks.  (We will give a more formal definition in the next section.)
$\Delta_k$ is well (but not perfectly) understood
algebraically \cite{IK}. Geometrically it is not well understood at all.

\subsection{The Maps and the Birds}

\noindent
{\bf Definition of a Polygon:\/}
For us, a {\it polygon\/} is a choice of both vertices and the edges
connecting them.
Each polygon $P$ we consider will all be {\it planar\/}, in the
sense that there is some projective transformation that
maps $P$, both vertices and edges, to the affine patch.
Our classical example is a regular $n$-gon, with the
obvious short edges chosen.
\newline
\newline
\noindent
{\bf The Maps:\/}
Given a polygon $P$, let
$P_a$ denote the $(a)$th vertex of $P$.
Let $P_{ab}$ be the line through $P_a$ and $P_b$.
The vertices of
$\Delta_k(P)$ are
\begin{equation}
  P_{j,j+k+1} \cap P_{j+1,j-k}.
\end{equation}
Here the indices are taken mod $n$.
Figure 1.3 shows this for $(k,n)=(2,7)$.
The polygons in Figure 1.3 are examples of
a concept we shall define shortly, that of a
$k$-bird.
\begin{center}
\resizebox{!}{1.8in}{\includegraphics{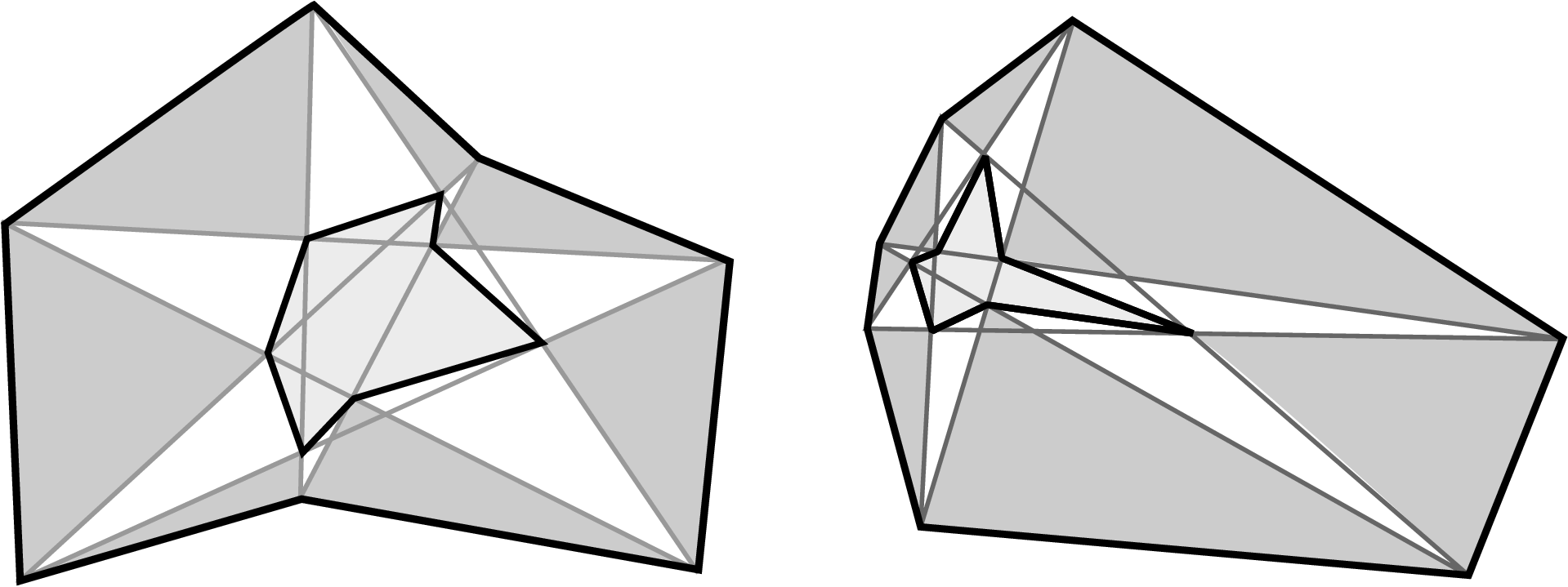}}
\newline
Figure 1.3: $\Delta_2$ acting on $2$-birds.
\end{center}
We should say a word about how the edges are defined.
In the case for the regular $n$-gon we make the
obvious choice, discussed above.  In general,
we define the class of polygons we consider in terms
of a homotopy from the regular $n$-gon.  So, in general,
we make the edge choices so that the edges vary continuously.
\newline
\newline
{\bf The Birds:\/}
Given an $n$-gon $P$, we let $P_{a,b}$ denote the line containing
the vertices $P_a$ and $P_b$.
We call $P$ $k$-{\it nice\/} if $n > 3k$, and $P$ is planar, and the $4$ lines
\begin{equation}
  P_{i,i-k-1}, \hskip 10 pt P_{i,i-k},  \hskip 10 pt P_{i,i+k}, \hskip 10 pt P_{i,i+k+1}
  \end{equation}
  are distinct for all $i$.
  It is not true that the generic $n$-gon is $k$-nice, because there
  are open sets of non-planar polygons.  (Consider a neighborhood of
  $P$,  where $P$
  the regular $100$-gon with the opposite choice of edges.)  However,
  the generic perturbation of a planar $n$-gon is also $k$-nice.

We call $P$ a $k$-{\it bird\/} if $P$ is the endpoint
of a path of $k$-nice $n$-gons that starts with the regular
$n$-gon. 
We let $B_{k,n}$ be the subspace of $n$-gons which
are $k$-birds.  Note that $B_{k,n}$ contains the set of
convex $n$-gons, and the containment is strict when $k>1$.
As Figure 1.3 illustrates, a $k$-bird need not be convex for
$k \geq 2$.  We will show that $k$-birds are always star-shaped,
and in particular embedded.  As we mentioned above, we use the homotopic definition
of a $k$-bird, to define the edges of $\Delta_k(P)$ when
$P$ is a $k$-bird.
\newline
\newline
{\bf Example:\/}
The homotopy part of our definition looks a bit strange, but it
is necessary.  To illustrate this, we consider the picture further
for the case $k=1$.
In this case, a $1$-bird must be convex, though
the $1$-niceness condition just means planar and locally
convex. Figure 1.4 shows how we might attempt a homoropy
from the regular octagon to a locally convex octagon
which essentially wraps twice around a quadrilateral.
The little grey arrows give hints about how the points are moved.
At some times,
the homotopy must break the $1$-niceness condition.
The two grey polygons indicate failures and the
highlighted vertices indicate the sites of the failures.
There might be other failures as well;  we are
taking some jumps in our depiction.

\begin{center}
\resizebox{!}{2.4in}{\includegraphics{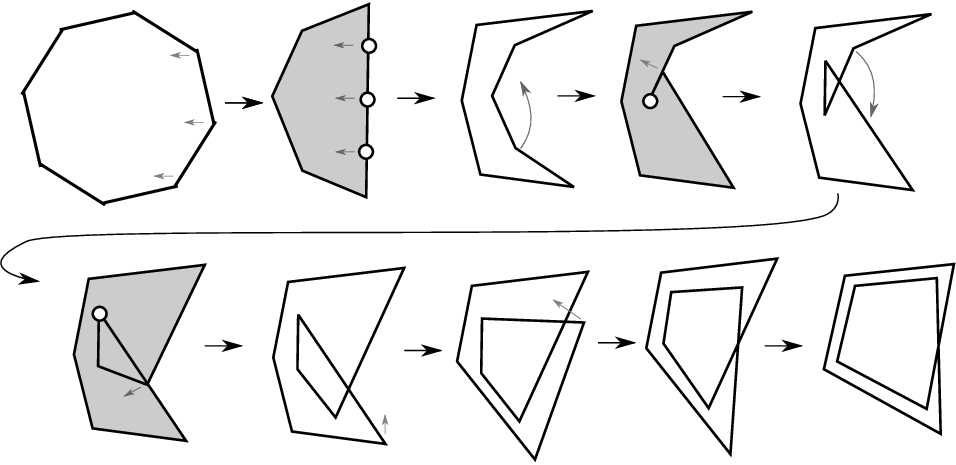}}
\newline
Figure 1.4: A homotopy that cannot stay $1$-nice.
\end{center}

One could make similar pictures when $k \geq 1$, but the pictures
might be harder to understand.

\subsection{The Main Result}

Given an embedded planar polygon
$P$, let $P^I$ denote the interior of
region bounded by $P$.
We say that $P$ is
{\it strictly star shaped\/} with respect to
$x \in P^I$ if
each ray emanating from $x$ intersects $P$ exactly once.
More simply, we say that $P$ is {\it strictly star shaped\/}
if it is strictly star shaped with respect to some point
$x \in P^I$.
Here is the main result.

\begin{theorem}
  \label{main}
   Let $k \geq 2$ and $n>3k$ and
   $P \in B_{k,n}$.  Then
   \begin{enumerate}
   \item $P$ is strictly star-shaped, and in particular embedded.
   \item $\Delta_k(P) \subset P^I$.
   \item $\Delta_k(B_{k,n})=B_{k,n}$.
   \end{enumerate}
\end{theorem}

\noindent
{\bf Remark:\/}
The statement that $n>3k$ is present just for emphasis.  $B_{n,k}$ is
by definition empty when $n \leq 3k$.
The restriction $n>3k$ is necessary.
    Figure 1.5 illustrates what would be a counter-example
    to Theorem \ref{main} for the pair $(k,n)=(3,9)$.
The issue is that a certain triple of $4$-diagonals
has a common intersection point.  This does not happen
for $n>3k$.  See Lemma \ref{struct}.

\begin{center}
\resizebox{!}{3.6in}{\includegraphics{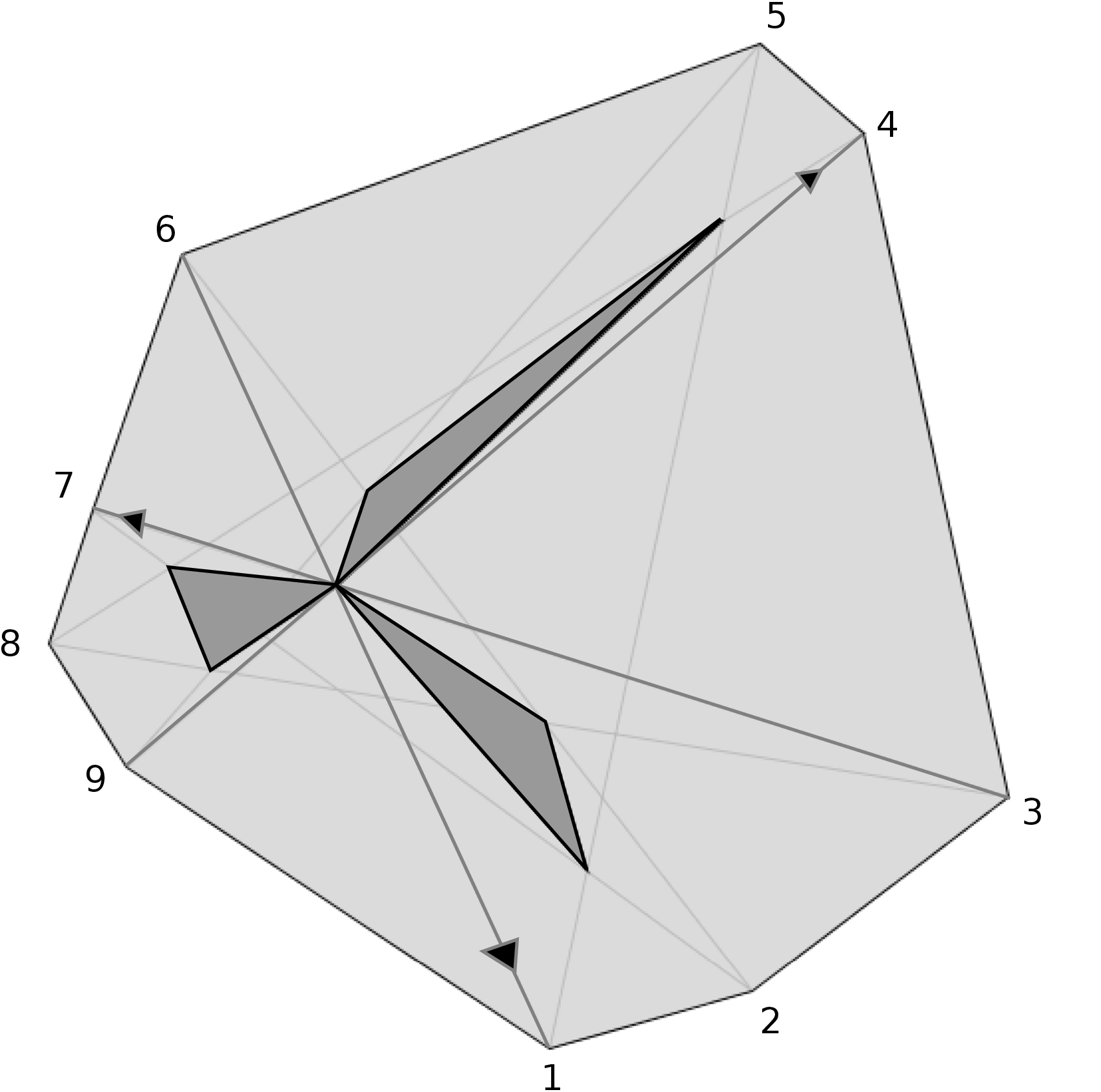}}
\newline
Figure 1.5: $\Delta_3$ acting on a certain convex $9$-gon.
\end{center}

\subsection{The Energy}
\label{eee}

We will deduce Statements 1 and 2 of Theorem \ref{main}
in a geometric way.  The key to proving Statement 3 is
a natural quantity associated to a $k$-bird.
We let $\sigma_{a,b}$ be the slope of the line $P_{a,b}$
and we define the {\it cross ratio\/}
\begin{equation}
  \chi(a,b,c,d)=\frac{(a-b)(c-d)}{(a-c)(b-d)}.
  \end{equation}
We define
{\small
  \begin{equation}
    \label{ENERGY}
  \chi_k(P)=\prod_{i=1}^n \chi(i,k,P), \hskip 20 pt
  \chi(i,k,P)=\chi(\sigma_{i,i-k},\sigma_{i,i-k-1},\sigma_{i,i+k+1},\sigma_{i,i+k})
\end{equation}
  \/}
Here we are taking the cross ratio the slopes the lines
involved in our definition of $k$-nice. 
When $k=1$ this is the familiar invariant $\chi_1=OE$ for
the pentagram map $\Delta_1$. See \cite{SCH1}, \cite{SCH2}, \cite{OST1},
\cite{OST2}. When $n=3k+1$, a
suitable star-relabeling of our polygons converts
$\Delta_k$ to $\Delta_1$ and $\chi_k$ to $1/\chi_1$.
So, in this case $\chi_k \circ \Delta_k=\chi_k$.
Figure 1.5 illustrates this for $(k,n)=(3,10)$.
Note that the polygons suggested by the dots in Figure 1.5 are \underline{not}
convex.  Were we to add in the edges we would get a highly non-convex pattern.

\begin{center}
\resizebox{!}{2in}{\includegraphics{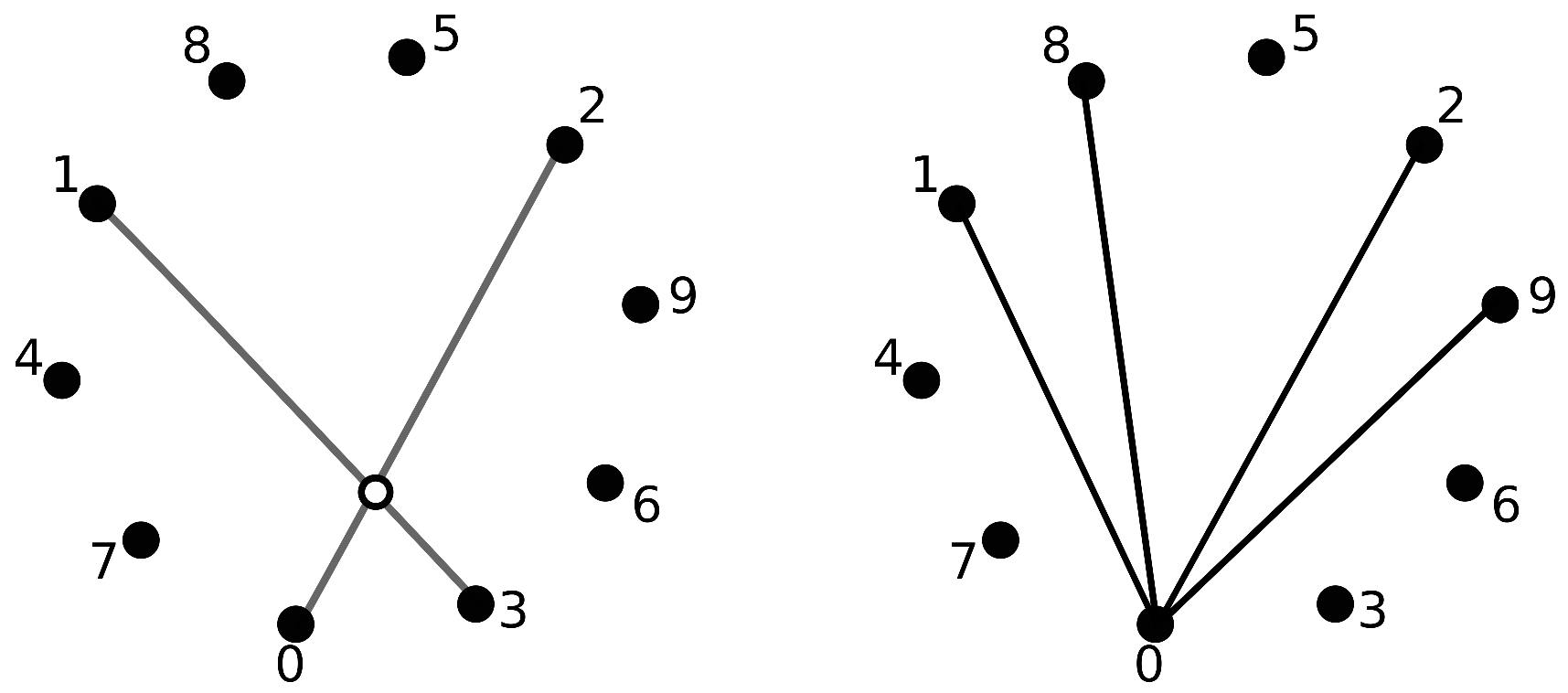}}
\newline
Figure 1.6: A star-relabeling converts $\Delta_1$ to $\Delta_3$ and $1/\chi_1$ to $\chi_3$.
\end{center}

In general, $\chi_k$ is not as clearly related to $\chi_1$.  Nonetheless,
we will prove
\begin{theorem}
  \label{energy}
  $\chi_k \circ \Delta_k=\chi_k$.
\end{theorem}
Theorem \ref{energy} is meant to hold for all
$n$-gons, as long as all quantities are defined.
There is no need to restrict to birds.

\subsection{The Collapse Point}

When it is understood that $P \in B_{k,n}$ it is convenient to write
\begin{equation}
  P^{\ell}=\Delta_k^{\ell}(P)
\end{equation}
We also let $\widehat P$ denote the closed planar region bounded by $P$.
Figure 1.7 below shows
$\widehat P=\widehat P^0,\widehat P^1,\widehat P^2,\widehat P^3,\widehat P^4$ for some
$P \in B_{4,13}$.

\begin{center}
\resizebox{!}{2.53in}{\includegraphics{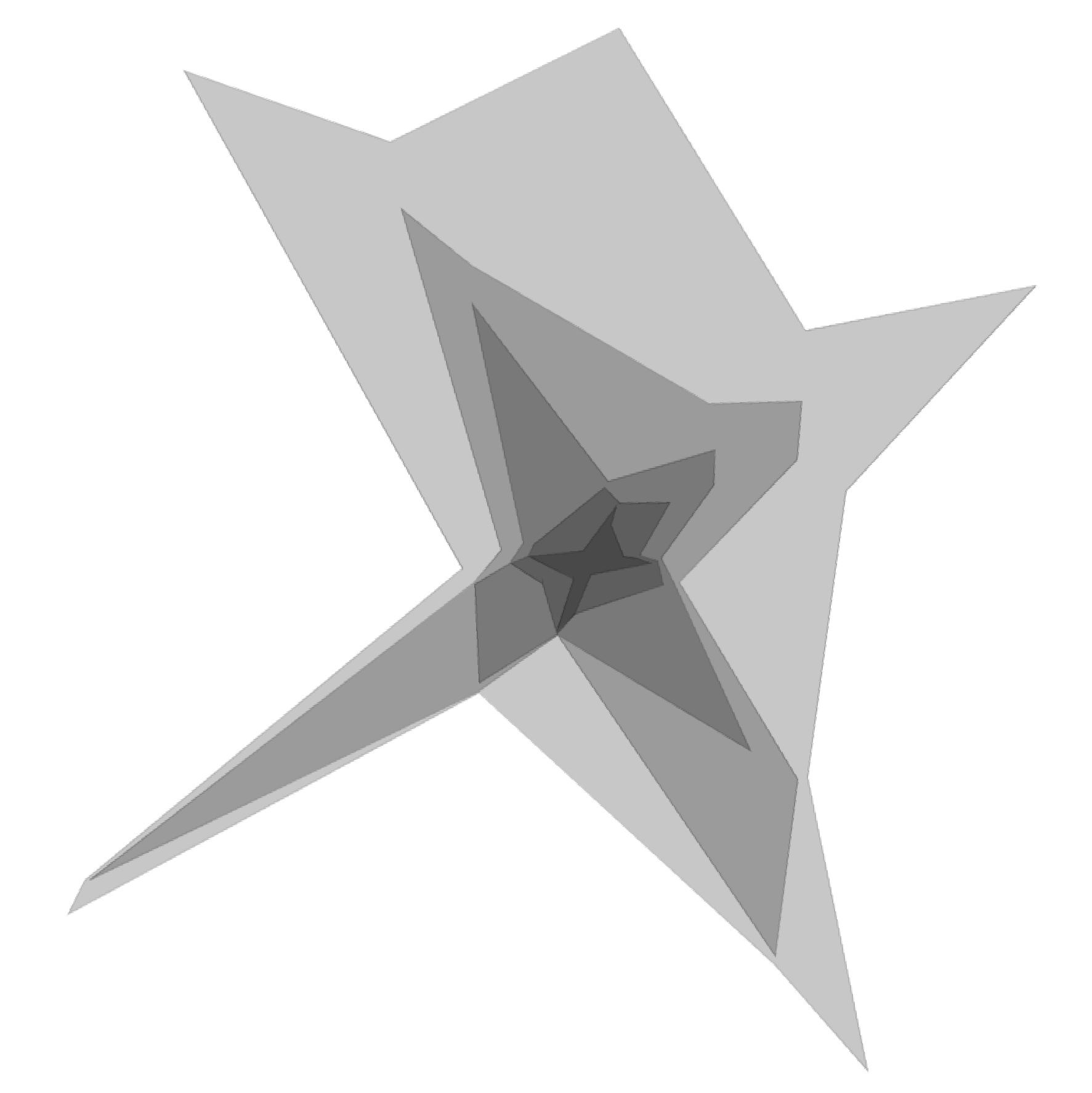}}
\newline
Figure 1.7: $\Delta_4$ and its iterates acting on a member of $B_{4,13}$.
\end{center}

Define
\begin{equation}
  \widehat P_{\infty}=\bigcap_{\ell \in \Z} \widehat P^{\ell}, \hskip 30 pt
  \widehat P_{-\infty}=\bigcup_{\ell \in \Z} \widehat P^{\ell}.
\end{equation}

\begin{theorem}
  \label{aux}
  If $P \in B_{k,n}$ then
  $\widehat P_{\infty}$ is a point and
  $\widehat P_{-\infty}$ is an affine plane.
\end{theorem}

Our argument will show that $P \in B_{k,n}$ is strictly
star-shaped with respect to all points in $\widehat P^n$.  In
particular, all polygons in the orbit are
strictly star-shaped with respect to the collapse point $\widehat P_{\infty}$.
See Corollary \ref{star-extra}.

One might wonder if some version of Glick's formula works for
the $\widehat P_{\infty}$ in general.
I discovered experimentally that this is indeed the case for
$n=3k+1$ and $n=3k+2$.    See \S \ref{Glick} for a discussion
of this and related matters.

Here is a corollary of our results that is just about convex polygons.

\begin{corollary}
  Suppose that $n>3k$ and $P$ is a convex $n$-gon.
  Then the sequence $\{\Delta_k^{\ell}(P)\}$ shrinks to a
  point as $\ell \to \infty$, and each member of this
  sequence if strictly star-shaped with respect to the
  collapse point.
\end{corollary}

\subsection{The Triangulations}

In \S \ref{triang}
we associate to each $k$-bird $P$ a triangulation $\tau_P  \subset
\P$, the
projective plane.
Here $\tau_P$ is an embedded
degree $6$ triangulation of
$P_{-\infty}-P_{\infty}$.
The edges are made from the segments in the $\delta$-diagonals of
$P$ and its iterates for $\delta=1,k,k+1$.

\begin{center}
\resizebox{!}{4.75in}{\includegraphics{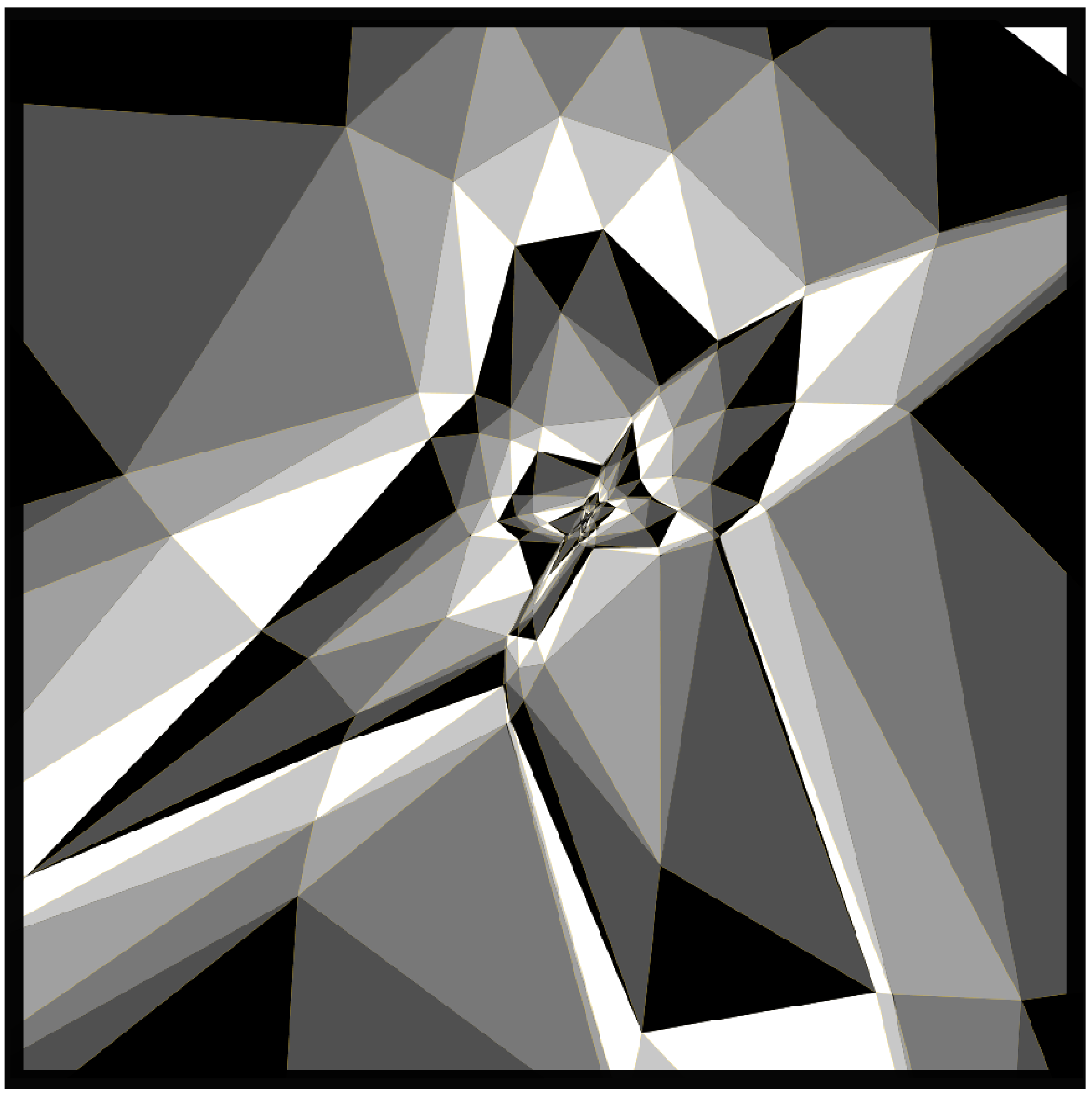}}
\newline
Figure 1.8: The triangulation associated to a member of $B_{5,16}$.
\end{center}

Figure 1.8 shows this tiling associated to a member of $B_{5,16}$.
In this figure, the interface between the big black triangles
and the big white triangles is some $\Delta_5^{\ell}(P)$ for
some smallish value of $\ell$. (I zoomed into the picture a bit
to remove the boundary of the initial $P$.)
The picture is normalized so that
the line $P_{-\infty}$ is the line at infinity.  When I make these
kinds of pictures (and animations), I
normalize so that the ellipse of inertia of $P$ is the
unit disk.

          \subsection{Paper Organization}

This paper is organized as follows.
\begin{itemize}
  \item In \S 2 we prove Theorem \ref{energy}.
  \item In \S 3 we prove Statement 1 of Theorem \ref{main}.
\item In \S 4 we prove Statement 2 of Theorem \ref{main}.
\item In \S 5 we prove a technical result called the Degeneration Lemma,
  which will help with  Statement 3 of Theorem \ref{main}.
\item In \S 6 we prove Statement 3 of Theorem \ref{main}.
\item In \S 7 we introduce the triangulations discussed above.
Our Theorem \ref{soul2} will help with the proof of Theorem \ref{aux}.
\item In \S 8 we prove Theorem \ref{aux}. 
\item In \S 9, an appendix,
  we sketch an alternate proof of Theorem \ref{energy}
  which Anton Izosimov kindly explained.  We also discuss Glick's
  collapse formula and star relabelings of polygons.
  \end{itemize}

\subsection{Visit the Flapping Bird Exhibit}

Our results inject some more geometry into the pentagram zoo.
Our results even have geometric implications for the pentagram map itself.
See \S \ref{relabel000}.
There are different ways to visit the flapping bird exhibit in the zoo.
You could read the proofs here, or you might
just want to to look at some images:
\newline 
         {\bf http://www.math.brown.edu/$\sim$reschwar/BirdGallery\/} \newline
         You can also download and play with the software I wrote:
 \newline 
          {\bf http://www.math.brown.edu/$\sim$reschwar/Java/Bird.TAR\/} \newline
          The software has detailed instructions.
          You can view this paper as a
justification for why the nice images actually exist.

\subsection{Acknowledgements}

I would like to thank
Misha Gekhtman,
Max Glick,
Anton Izosimov,
Boris Khesin,
Valentin Ovsienko,
and Serge Tabachnikov
for many discussions about the pentagram zoo.
I would like to thank Anton, in particular, for
extensive discussions about the material in \S 9.

\newpage

\section{The Energy}

The purpose of this chapter is to prove Theorem \ref{energy}.  The
proof, which is similar to what I do in \cite{SCH1}, is more of a
verification than a conceptual explanation.  My computer program
allows the reader to understand the technical details of the proof
better.  The reader might want to just skim this chapter on the first
reading.  In \S \ref{ANTON} I will sketch an alternate proof, which I
learned from Anton Izosimov.  Izosimov's proof also uses the first two
sections of this chapter.

\subsection{Projective Geometry}

 Let $\P$ denote the real projective plane.
 This is the space of $1$-dimensional subspaces
 of $\R^3$.  The projective plane $\P$ contains $\R^2$ as the
 {\it affine patch\/}.   Here
 $\R^2$ corresponds to vectors of the form
 $(x,y,1)$, which in turn define elements of $\P$.
 
Let $\P^*$ denote the dual projective plane,
namely the space of lines in $\P$.  The elements
in $\P^*$ are naturally equivalent to
$2$-dimensional subspaces of $\R^3$.
The line in $\P$ such a subspace $\Pi$ defines is equal to
the union of all $1$-dimensional subspaces of $\Pi$.

Any invertible linear transformation of $\R^3$ induces a
{\it projective transformation\/} of $\P$, and also of $\P^*$.
These form the {\it projective group\/} $PSL_3(\R)$.
Such maps preserve collinear points and coincident lines.

A {\it duality\/}
from $\P$ to $\P^*$ is an analytic diffeomorphism
$\P \to \P^*$ which maps collinear points to
coincidence lines.  The classic example is the
map which sends each linear subspace of
$\R^3$ to its orthogonal complement.

A {\it PolyPoint\/} is a cyclically ordered list
of points of $\P$.  When there are $n$ such points,
we call this an $n$-{\it Point\/}.
A {\it PolyLine\/} is a cyclically ordered list
of lines in $\P$, which is the same as a cyclically
ordered list of points in $\P^*$.
A projective duality maps PolyLines to PolyPoints,
and {\it vice versa\/}.

Each $n$-Point determines $2^n$ polygons in $\P$
because, for each pair of consecutive points, we may
choose one of two line segments to join them.
As we mentioned in the introduction, we have a canonical
choice for $k$-birds.
Theorem \ref{energy} only involves
PolyPoints, and our proof uses PolyPoints and PolyLines.

Given a $n$-Point $P$, we let $P_j$ be its $j$th point.
We make a similar definition for $n$-Lines.  We always
take indices mod $n$.

\subsection{Factoring the Map}
\label{duality}

Like the pentagram map, the map $\Delta_k$ is the product of
$2$ involutions.  This factorization will be useful here and
in later chapters.

Given a PolyPoint $P$, consisting of points $P_1,...,P_n$,
we define $Q=D_m(P)$ to be the PolyLine  whose successive
lines are $P_{0,m}$, $P_{1,m+1}$, etc.  Here
$P_{0,m}$ denotes the line through $P_0$ and $P_{m}$, etc.
We labed the vertices so that
\begin{equation}
  \label{convention}
  Q_{-m-i} = P_{i,i+m}.
\end{equation}
This is a convenient choice.
We define the action of $D_m$ on PolyLines in the same way,
switching the roles of points and lines.  For PolyLines,
$P_{0,m}$ is the intersection of the line $P_0$ with the line
$P_m$.  The map $D_m$ is an involution which swaps PolyPoints with
PolyLines.
We have the compositions
\begin{equation}
  \label{factor}
  \Delta_k = D_k \circ D_{k+1}, \hskip 30 pt
  \Delta_k^{-1}=D_{k+1} \circ D_k.
\end{equation}

The {\it energy\/} $\chi_k$ makes sense for $n$-Lines as well
as for $n$-Points.  The quantities
$\chi_k \circ D_k(P)$ and $\chi_k \circ D_{k+1}(P)$ can be computed
directly from the PolyPoint $P$.  Figure 2.1 shows schematically the $4$-tuples
associated to $\chi(0,k,Q)$ for $Q=P$ and  $D_k(P)$ and $D_{k+1}(P)$.
In each case, $\chi_k(Q)$ is a product of $n$ cross ratios like these.
If we want to compute the factor of $\chi_k(D_k(P))$ associated
to index $i$ we
subtract  (rather than add)
$i$ from the indices shown in the middle figure. A similar
rule goes for $D_{k+1}(P)$.

\begin{center}
\resizebox{!}{1.7in}{\includegraphics{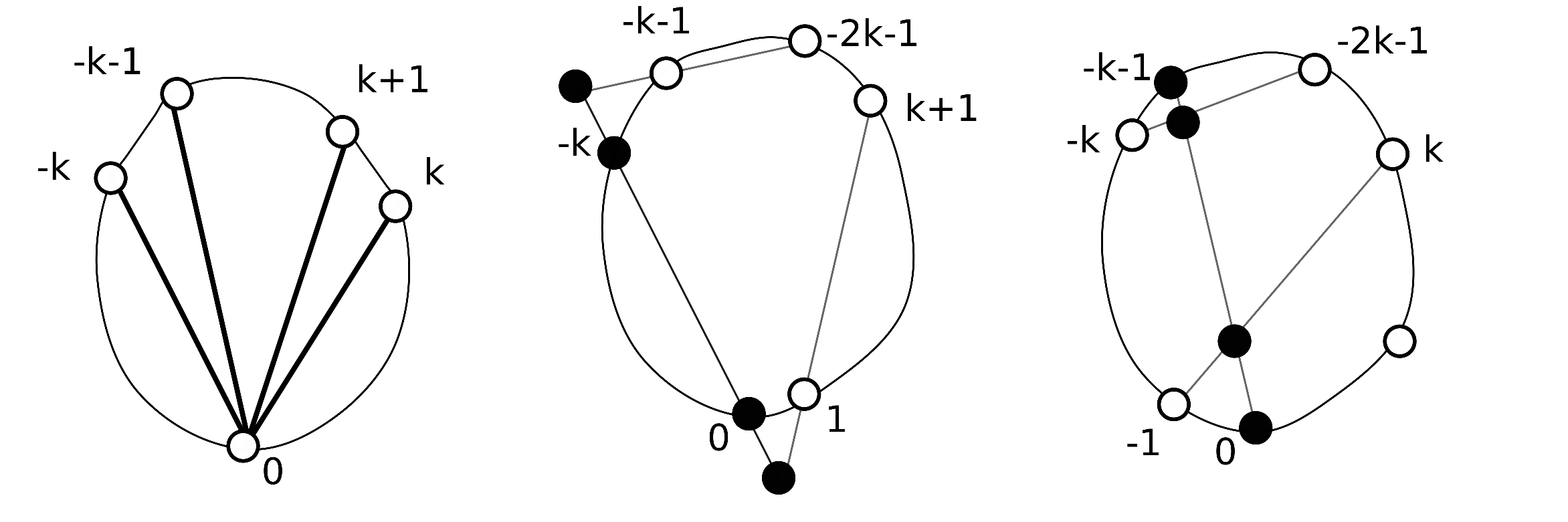}}
\newline
Figure 2.1: Computing the $k$-energy.
\end{center}

Theorem \ref{energy} follows from the next two results.

\begin{lemma}
  \label{energy1}
  $\chi_k \circ D_k= \chi_k$.
\end{lemma}

\begin{lemma}
  \label{energy2}
  $\chi_k \circ D_{k+1}= \chi_k$.
\end{lemma}

These results have almost identical proofs.  We
consider Lemma \ref{energy1} in detail and then explain
the small changes needed for Lemma \ref{energy2}.

\subsection{Proof of the First Result}

We study the ratio
\begin{equation}
  R(P)=\frac{\chi_k \circ D_k(P)}{\chi_k(P)}.
\end{equation}
We want to show that $R(P)$ equals $1$ wherever it is defined.
We certainly have $R(P)=1$ when $P$ is the regular $n$-Point.

Given a PolyPoint $P$ we choose a pair of vertices
$a,b$ with $|a-b|=k$.
We define $P(t)$ to be the
PolyPoint obtained by replacing $P_a$ with
\begin{equation}
  (1-t)P_a + t P_{b}.
\end{equation}
Figure 2.2 shows what we are talking about, in case $k=3$.
We have rotated the picture so that $P_a$ and $P_b$ both
lie on the $X$-axis.
\begin{center}
\resizebox{!}{1.3in}{\includegraphics{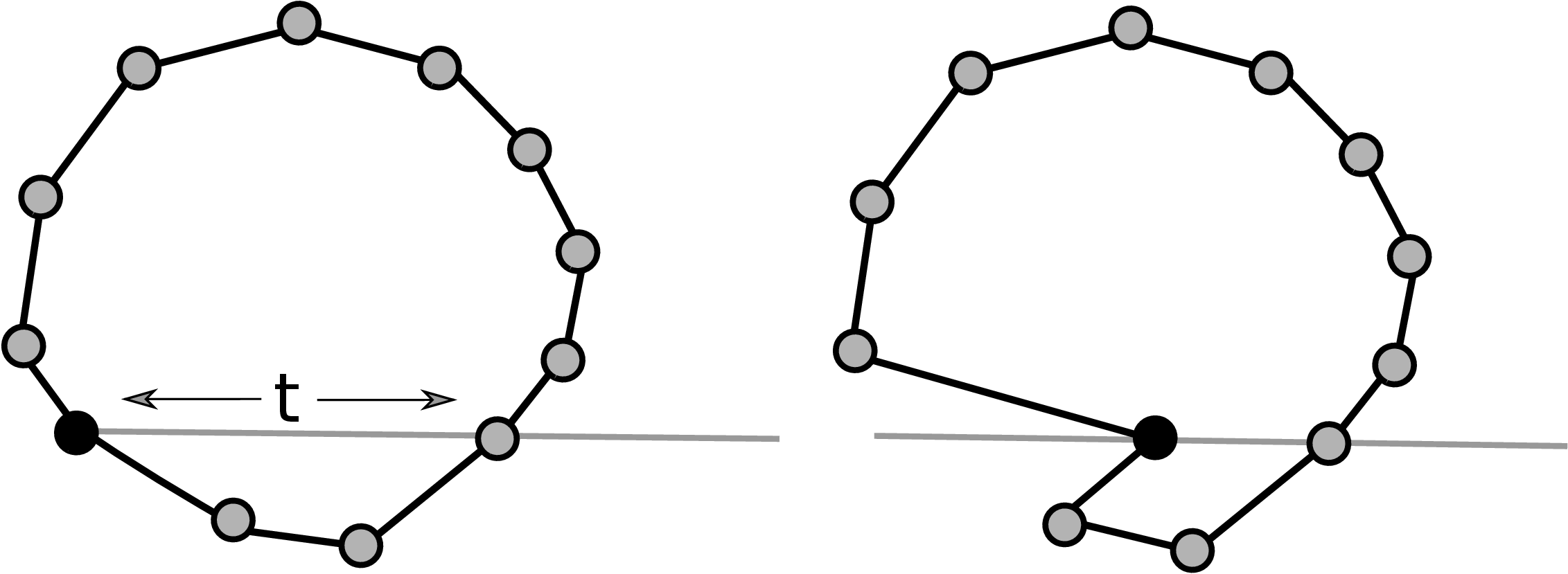}}
\newline
Figure 2.2: Connecting one PolyPoint to another by sliding a point.
\end{center}

The two functions
\begin{equation}
  f(t)=\chi_k(P(t)), \hskip 30 pt
  g(t)=\chi_k \circ D_k(P(t))
\end{equation}
are each rational functions of $t$.  Our notation does not
reflect that $f$ and $g$ depend on $P,a,b$.

A {\it linear fractional transformation\/}
is a map of the form
$$t \to \frac{ \alpha t+\beta }{\gamma t+\delta}, \hskip 30 pt
\alpha,\beta,\gamma,\delta \in \R, \hskip 30 pt \alpha \delta -\beta \gamma \not =0.
$$

\begin{lemma}[Factor I]
  If $n \geq 4k+2$ and $P$ is a generically chosen $n$-Point, then
  $f(t)$ and $g(t)$ are each products of
  $4$ linear fractional transformations.
  The zeros of $f$ and $g$ occur at the same points
  and the poles of $f$ and $g$ occur at the same points.
  Hence $f/g$ is constant.
\end{lemma}

The only reason we choose $n \geq 4k+2$ in the Factor Lemma is so that
the various diagonals involved in the proof do not have common
endpoints. 
The Factor Lemma I works the same way for all $k$ and for all
choices of (large) $n$.
We write $P \leftrightarrow Q$ if we can choose indices $a,b$ and
some $t \in \R$ such that $Q=P(t)$.  The Factor Lemma implies
that when $P,Q$ are generic and $P \leftrightarrow Q$ we have
$R(P)=R(Q)$.  The result for non-generic choices of $P$ follows
from continuity.  Any $n$-Point $Q$ can be included in a finite
chain
$$P_0 \leftrightarrow P_1 \leftrightarrow \cdots \leftrightarrow P_{2n} = Q,$$
where $P_0$ is the regular $n$-Point.
Hence $R(Q)=R(P_0)=1$.
This shows that Lemma \ref{energy1} holds for
$(k,n)$ where $k \geq 2$ and $n \geq 4k+2$.
(The case $k=1$ is a main result of
\cite{SCH1}, and by now has many proofs.)

\begin{lemma}
  If Lemma \ref{energy1} is true for
all large values of $n$, then it is true for
all values of $n$.
\end{lemma}

\startproof
If we are interested in the
result for small values of $n$, we can replace
a given PolyPoint $P$ with its $m$-fold cyclic
cover $mP$.  We have $\chi_k(mP)=\chi_k(P)^m$
and $\chi_k(D_k(mP))=\chi_k(D_k(p))^m$.
Thus, the result for large $n$ implies the
result for small $n$.
\endproof

In view of Equation
\ref{ENERGY} we have
\begin{equation}
  f(t)=f_1(t)...f_n(t), \hskip 30 pt f_j(t)=\chi(j,k,P(t)).
\end{equation}
Thus $f(t)$ is the product of $n$ ``local'' cross ratios.
We call an index $j$ {\it asleep\/} if none of the lines
involved in the cross ratio $f_j(t)$ depend on $t$.
In other words, the lines do not vary at all with $t$.
Otherwise we call $j$ {\it awake\/}.

As we vary $t$, only the diagonals
$P_{0,h}$ change for $h=-k,-k-1,k+1,k$.  From this fact,
it is not surprising that there are precisely $4$ awake
indices.  These indices are
\begin{equation}
  \label{vertices}
  j_0=0,\hskip 10 pt   j_1=k+1,\hskip 10 pt   j_2=-k-1,\hskip 10 pt   j_3=-k.
\end{equation}
The index $k$ is not awake because the diagonal
$P_{0,k}(t)$ does not move with $t$.

We define a {\it chord\/} of $P(t)$ to be a line defined
by a pair of vertices of $P(t)$.
The point $P_0(t)$ moves at linear speed, and
the $4$ lines involved in the calculation of
$f_{c_j}(t)$ are distinct unless
$P_0(t)$ lies in one of the chords of $P(t)$.
Thus $f_{c_j}(t)$ only has zeros and poles at
the corresponding values of $t$.  It turns out that
only the following chords are involved.
{\small
\begin{equation}
  \label{chords}
  \begin{matrix} -k \cr -k-1 \end{matrix}
  \hskip 15 pt
  \begin{matrix} -k \cr k+1 \end{matrix}
  \hskip 15 pt
  \begin{matrix} -k \cr 1 \end{matrix}
  \hskip 15 pt
  \begin{matrix} -k \cr -2k-1 \end{matrix}
  \hskip 15 pt
  \begin{matrix} -k-1 \cr -1 \end{matrix}
  \hskip 15 pt
  \begin{matrix} -k-1 \cr -2k-1 \end{matrix}
  \hskip 15 pt
 \begin{matrix} k+1\cr 1 \end{matrix}
  \hskip 15 pt
  \begin{matrix} k+1 \cr 2k+1\end{matrix}
\end{equation}
\/}
We call these $c_0,...,c_7$.
For instance, $c_0$ is the line through
$P_{-k}$ and $P_{-k-1}$.
Let $t_j$ denote the value of $t$ such that
$P(t_j) \in c_j$.

The PolyPoint $Q(t)=D_k(P(t))$ has the same
structure as $P(t)$.   Up to
projective transformations $Q(t)$ is also obtained from
the regular PolyPoint by moving a single vertex
along one of the $k$-diagonals.
The pattern of zeros and poles
is not precisely the same because the
chords of $Q(t)$ do not correspond to the
chords of $P(t)$ in a completely straightforward
way.  The $k$-diagonals of $Q(t)$ correspond to the
vertices of $P(t)$ and {\it vice versa\/}. The
$(k+1)$ diagonals of $Q(t)$ correspond to the
vertices of $\Delta_k^{-1}(P(t))$. This is what
gives us our quadruples of points in the middle
picture in Figure 2.1.

We now list the pattern of zeros and poles.  We explain our
notation by way of example. The quadruple $(f,2,4,5)$ indicates
that $f_{c_2}$ has a simple zero at $f_4$ and a simple pole at $t_5$.
\begin{equation}
  (f,0,0,1), \hskip 20 pt
  (f,1,6,7), \hskip 20 pt
  (f,2,4,5), \hskip 20 pt
  (f,3,2,3).
\end{equation}
\begin{equation}
  (g,0,6,5), \hskip 20 pt
  (g,1,0,3), \hskip 20 pt
  (g,2,2,1), \hskip 20 pt
  (g,3,4,7).
\end{equation}
Since these functions have holomorphic extensions to $\C$ with no
other zeros and poles, these functions are
linear fractional transformations.  This pattern establishes
the Factor Lemma I.

Checking that the pattern is correct is just a matter of inspection.
We give two example checks.
\begin{itemize}
\item To see why $f_{c_2}$ has a simple zero at $t_4$ we consider
  the quintuple $$(-k-1,-2k-1,-2k-2,0,-1).$$
  At $t_4$ the two
  diagonals $P_{-k-1,0}$ and $P_{-k-1,-1}$ coincide.  In terms of
  the cross ratios of the slopes we are computing
  $\chi(a,b,c,d)$ with $a=b$.  The point $P_0(t)$ is moving with
  linear speed and so the zero is simple.
\item To see why $g_{c_2}$ has a simple pole at $t_1$ we
  consider the $4$ points
  \begin{equation}
    \label{exam1}
    P_{2k+2,k+2}  \cap P_{1,k+1},\hskip 10 pt
  P_{k+1}, \hskip 10 pt
  P_{1}, \hskip 10 pt
  P_{1,k+1} \cap P_{-k,0}.
  \end{equation}
  These are all contained in the $k$-diagonal $P_{1,k+1}$, which
  corresponds to the vertex $(-k-1)$ of $D_k(P)$.
  At $t=t_1$ the three points $P_0(t)$ and $P_{-k}$ and $P_{k+1}$
  are collinear.  This makes the $2$nd and $4$th listed
  point coincided.  In terms of our cross ratio
  $\chi(a,b,c,d)$ we have $b=d$. This gives us a pole.
  The pole is simple because the points come together at
  linear speed.
  \end{itemize}

  The other explanations are similar.  The reader can see
  graphical illustrations of these zeros and poles using
  our program.

\subsection{Proof of the Second Result}

The proof of Lemma \ref{energy2} is essentially
identical to the proof of Lemma \ref{energy1}.
Here are the changes.  The Factor Lemma II
has precisely the same statement as
the Factor Lemma I, except that
\begin{itemize}
 \item  When defining $P(t)$ we use points $P_a$ and $P_b$ with
 $|a-b|=k+1$.
\item We are comparing $P(t)$ with $D_{k+1}(P(t))$.
\end{itemize}
This changes the definition of the functions $f$ and $g$.
With these changes made, the Factor Lemma I is replaced
by the Factor Lemma II, which has an identical statement.
This time the chords involved are as follows.

{\small
\begin{equation}
  \label{chords2}
  \begin{matrix} -k-1\cr -k \end{matrix}
  \hskip 15 pt
  \begin{matrix} -k-1\cr k \end{matrix}
  \hskip 15 pt
  \begin{matrix} -k-1 \cr -1 \end{matrix}
  \hskip 15 pt
  \begin{matrix} -k-1 \cr -2k-1 \end{matrix}
  \hskip 15 pt
  \begin{matrix} -k \cr 1 \end{matrix}
  \hskip 15 pt
  \begin{matrix} -k \cr -2k-1 \end{matrix}
  \hskip 15 pt
 \begin{matrix} k\cr -1 \end{matrix}
  \hskip 15 pt
  \begin{matrix} k\cr 2k+1\end{matrix}
\end{equation}
\/}

This time the $4$ awake indices are:
\begin{equation}
  \label{vertices2}
  j_0=0,\hskip 10 pt   j_1=k,\hskip 10 pt   j_2=-k-1,\hskip 10 pt   j_3=-k.
\end{equation}
Here is the pattern of zeros and poles.
\begin{equation}
  (f,0,1,0), \hskip 20 pt
  (f,1,7,6), \hskip 20 pt
  (f,2,3,2), \hskip 20 pt
  (f,3,5,4).
\end{equation}
\begin{equation}
  (g,0,5,6), \hskip 20 pt
  (g,1,3,0), \hskip 20 pt
  (g,2,7,4), \hskip 20 pt
  (g,3,1,2).
\end{equation}
The pictures in these cases look almost identical to the previous case.
The reader can see these pictures by operating my computer program.
Again, the zeros of $f$ and $g$ are located at the same places,
and likewise the poles of $f$ and $g$ are located
at the same places.  
Hence $f/g$ is constant.  This completes the proof
the Factor Lemma II, which implies
Lemma \ref{energy2}.

\newpage

\section{The Soul of the Bird}

\subsection{Goal of the Chapter}
\label{goal2}

Given a polygon $P \subset \R^2$, let $\widehat P$ be the
closure of the bounded components of $\R^2-P$ and
let $P^I$ be the interior of $\widehat P$. (Eventually
we will see that birds are embedded, so $\widehat P$ will
be a closed topological disk and $P^I$ will be an open
topological disk.)

Suppose now that $P(t)$ for $t \in [0,1]$ is a path
in $B_{n,k}$ starting at the regular $n$-gon $P(0)$.
We can adjust by a
continuous family of projective transformations
so that $P(t)$ is a bounded polygon in $\R^2$ for all
$t \in [0,1]$.
We orient $P(0)$ counter-clockwise around $P^I(0)$.
We extend this orientation choice continuously to $P(t)$.
We let $P_{ab}(t)$ denote the diagonal through
vertices $P_a(t)$ and $P_b(t)$. We orient $P_{a,b}(t)$ so that
it points from $P_a(t)$ to $P_b(t)$.
We take indices mod $n$.

We now recall a definition from the introduction:
When $P$ is embedded, we say that $P$ is
{\it strictly star shaped\/} with respect to
$x \in P^I$ if
each ray emanating from $x$ intersects $P$ exactly once.

\begin{center}
\resizebox{!}{2.53in}{\includegraphics{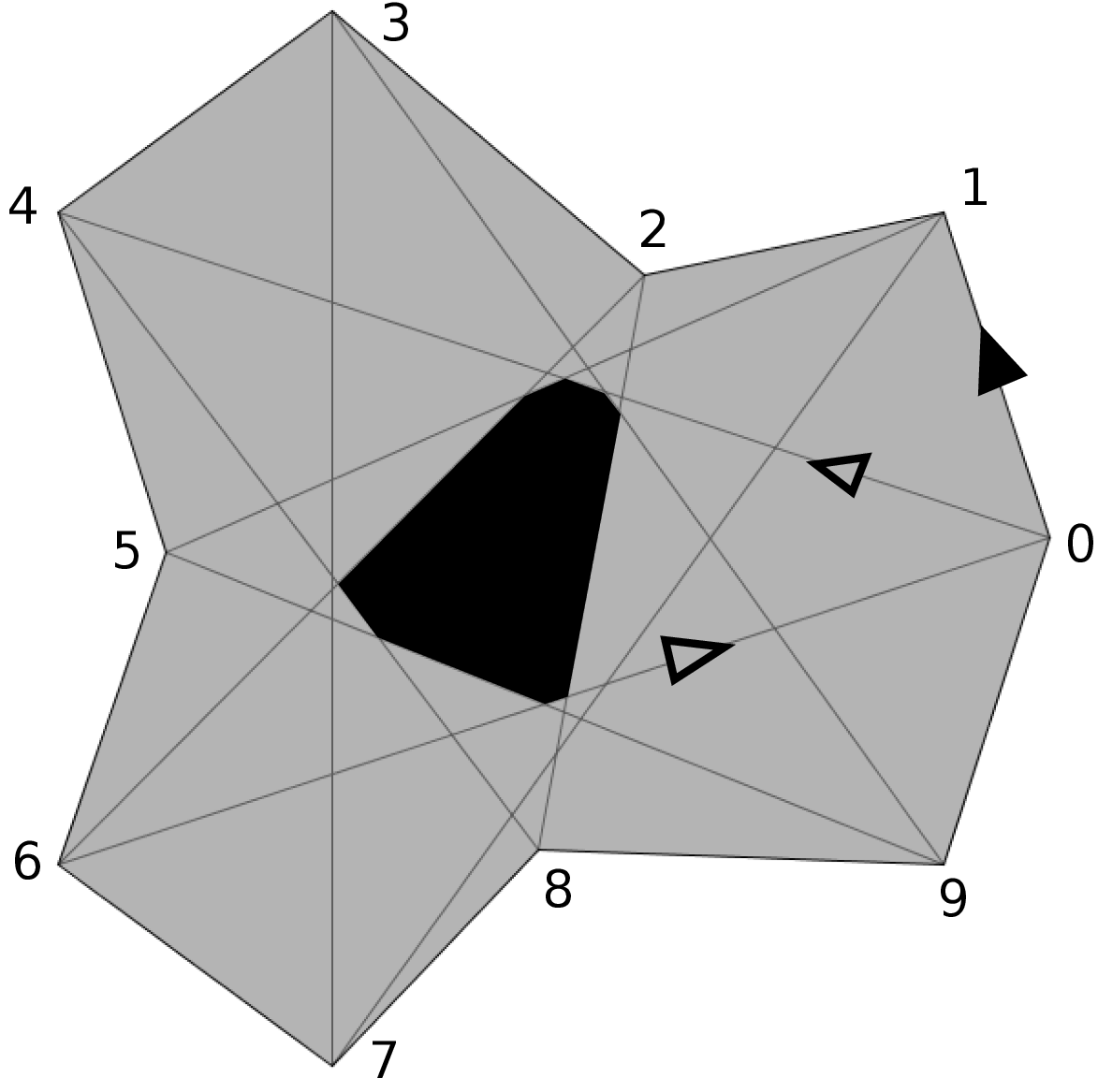}}
\newline
Figure 3.1:  The soul of a $3$-bird
\end{center}

Each such $(k+1)$-diagonal defines an
oriented line that contains it, and also the (closed)
{\it distinguished half plane\/} which lies
to the \underline{left} of the oriented line.
These $n$  half-planes vary continuously with $t$.
The {\it soul\/} of $P(t)$, which we denote $S(t)$, is the
intersection of the distinguished half-planes. Figure
3.1 shows the an example.

The goal of this chapter is to prove the following result.
\begin{theorem}
  \label{soul}
  Let $P$ be a bird and let $S$ be its soul.
  Then:
  \begin{enumerate}
\item $S$ is has non-empty interior.
\item $S \subset P^I$.
\item $P$ is strictly star-shaped with respect to any point in $S$.
  \end{enumerate}
\end{theorem}
Theorem \ref{soul} immediately implies Statement 1 of Theorem \ref{main}.

We are going to give a homotopical proof of
Theorem \ref{soul}.  We say that a value
$t \in [0,1]$ is a {\it good parameter\/}
if Theorem \ref{soul} holds for $P(t)$.  All three conclusions of
Theorem \ref{soul} are open conditions. Finally,
$0$ is a good parameter.  For all these reasons, it suffices to
prove that the set of good parameters is closed.  By truncating
our path at the first supposed failure, we reduce to the
case when Theorem \ref{soul}
holds for all $t \in [0,1)$.

  \subsection{The Proof}

  For ease of notation we set $X=X(1)$ for any object
  $X$ associated to $P(1)$.

\begin{lemma}
  \label{L1}
  If $P$ is any $k$-bird, then
    $P_0$ and $P_{2k+1}$ lie to the left of
    $P_{k,k+1}$.  The same goes if all indices
    are cyclically shifted by the same amount.
\end{lemma}

\startproof
Consider the triangle
with vertices $P_0(t)$ and $P_k(t)$ and $P_{k+1}(t)$.
The $k$-niceness condition implies that this
triangle is non-degenerate for all $t \in [0,1]$.
Since $P_0(t)$ lies to to the left
of $P_{k,k+1}(t)$, the non-degeneracy implies
the same result for $t=1$. The same
argument works for the triple $(2k+1,k,k+1)$.
\endproof

     \begin{lemma}
       $S$ is non-empty and contained in $P^I$.
     \end{lemma}

     \startproof
     By continuity, $S$ is nonempty and contained in $P \cup P^I$.  
     By the $k$-niceness property and continuity, $P_1$ lies strictly to the right
     of $P_{0,k+1}$.  Hence the entire half-open edge
  $[P_0,P_1)$ lies strictly to the right of $P_{0,k+1}$.
  Hence $[P_0,P_1)$ is disjoint from $S$.
      By cyclic relabeling, the same goes for all the
      other half-open edges. Hence $S \cap P=\emptyset$.
      Hence $S \subset P^I$.
     \endproof

     \begin{lemma}
       $P$ is strictly star-shaped with respect to any point of $S$.
     \end{lemma}

     \startproof
     Since $P(t)$ is strictly star-shaped with respect to all points of
     $S(t)$ for $t<1$, this lemma can only fail if there is an edge
     of $P$ whose extending line contains a point $x \in S$.
     We can cyclically relabel so that the edge of $\overline{P_0P_1}$.
     
\begin{center}
\resizebox{!}{.75in}{\includegraphics{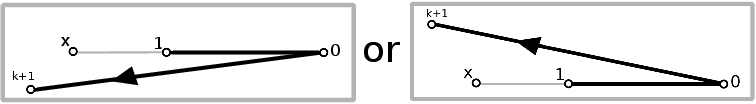}}
\newline
Figure 3.2:  The diagonal $P_{0,k+1}$ does not separate $1$ from $x$.
\end{center}

     Since $x \not \in P$, either $P_1$ lies between $P_0$ and $x$ or
     $P_0$ lies in between $x$ and $P_1$. In the first case, both
     $P_1$ and $x$ lie on the same side of the diagonal $P_{0,k+1}$.
     This is a contradiction: $P_1$ is supposed to lie on the right
     and $x$ is supposed to lie on the left.   In the second case
     we get the same kind of contradiction with respect to the
     diagonal $P_{-k,1}$.
     \endproof

We say that $P$ has {\it opposing $(k+1)$-diagonals\/}
if it has a pair of $(k+1)$-diagonals which lie in the
same line and point in opposite directions.
In this case, the two left half-spaces are on opposite
sides of the common line.

\begin{lemma}
    \label{oppo}
    $P$ does not have opposing $(k+1)$-diagonals.
  \end{lemma}

  \startproof
  We suppose that $P$ has opposing diagonals and
  we derive a contradiction.
  In this case $S$, which is the intersection of
  all the associated left half-planes, must be
  a subset of the line $L$ containing these diagonals.
  But then $P$ intersects $L$ in at least $4$ points,
  none of which lie in $S$.  But then $P$ cannot
  be strictly star-shaped with respect to any
  point of $S$. This is a contradiction.
  \endproof

We call three $(k+1)$-diagonals of $P(t)$ {\it interlaced\/} if the
intersection of their left half-spaces is a triangle.
See Figure 3.3.

\begin{center}
\resizebox{!}{2.2in}{\includegraphics{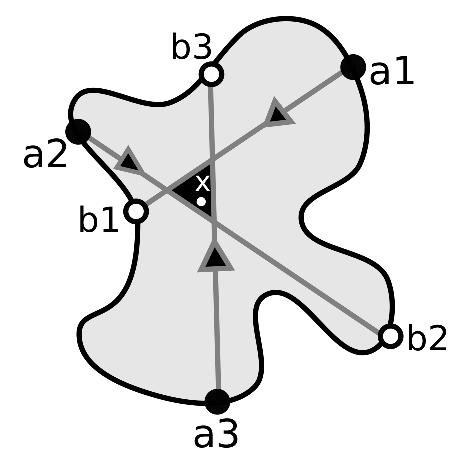}}
\newline
Figure 3.3:  Interlaced diagonals on $P(t)$.
\end{center}

Given interlaced $(k+1)$-diagonals, and a point $x$
in the intersection, the circle of rays emanating
from $x$ encounters the endpoints of the diagonals
in an alternating pattern:
$a_1,b_3,a_2,b_1,a_3,b_2$, where
$a_1,a_2,a_3$ are the tail points and
$b_1,b_2,b_3$ are the head points.
Here $a_1$ names the vertex $P_{a_1}(t)$, etc.

\begin{lemma}
  \label{struct}
  $P(t)$ cannot have interlaced diagonals for $t<1$.
 \end{lemma}

\startproof
Choose $x \in S(t)$. The $n$-gon $P(t)$ is strictly
star-shaped with respect to $x$.
Hence, the vertices of $P$ are encountered in
order (mod $n$) by a family of rays that emanate from $x$ and
rotates around full-circle.  Given the order these
vertices are encountered, we have
$a_{j+1}=a_j+\eta_j$,
where $\eta_j \leq k$.  Here we are taking
the subscripts mod $3$ and the vertex values mod $n$.
This tells us that $n=\eta_1+\eta_2+\eta_3 \leq 3k$.
This contradicts the fact that $n>3k$.
\endproof

It only remains to show that $S$ has non-empty interior.
  A special case of
  Helly's Theorem says the following: If we have a finite number
  of convex subsets of $\R^2$ then they all intersect provided
  that every $3$ of them intersect.  Applying Helly's Theorem to
  the set of interiors of our distinguished half-planes, we
  conclude that we can find $3$ of these open half-planes
  whose triple intersection is empty.  On the other hand,
  the triple intersection of the {\it closed\/} half-planes
  contains $x$. Since $P$ has no opposing diagonals,
  this is only possible if the $3$ associated
  diagonals are interlaced for $t$ sufficiently close to $1$. This contradicts
  Lemma \ref{struct}.   Hence $S$ has non-empty interior.
  
\newpage

\section{The Feathers of the Bird}
\label{featherAAA}

\subsection{Goal of the Chapter}

Recall that $P^I$ is the interior of the region bounded by $P$.
We call the union of black triangles in Figure 4.1 the
{\it feathers\/} of the bird.  the black region in the center is the soul.

\begin{center}
\resizebox{!}{3in}{\includegraphics{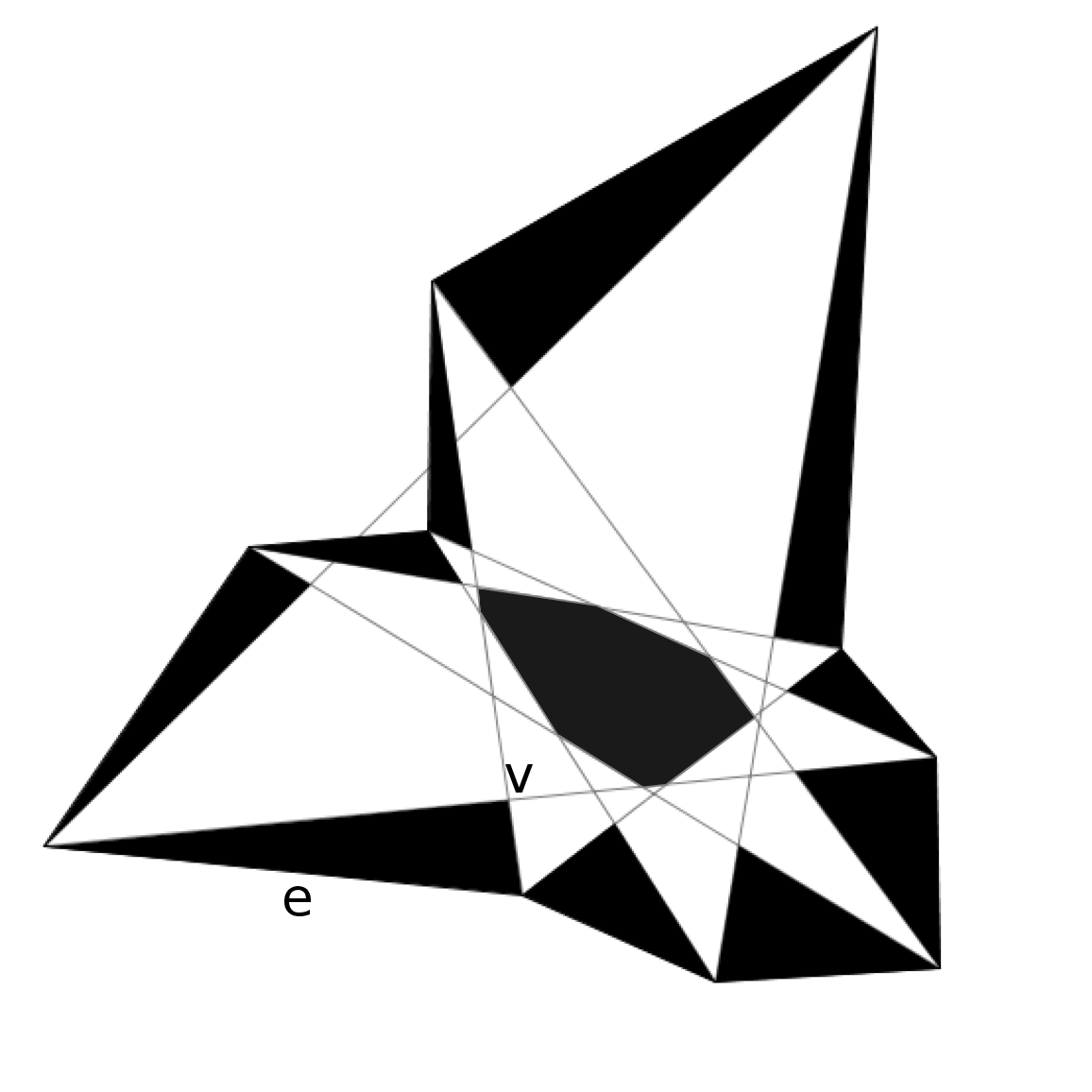}}
\newline
Figure 4.1 The feathers of a $3$-bird.
\end{center}

Each feather $F$ of a $k$-bird $P$ is the convex hull of
its {\it base\/}, an edge $e$ of $P$, and its {\it tip\/},
a vertex of $\Delta_k(P)$.

The goal of this chapter is to prove the following result,
which says that the simple topological picture shown in
Figure 4.1 always holds.

\begin{theorem}
  \label{feather}
  The following is true.
  \begin{enumerate}
  \item  Let $F$ be an feather with base $e$. Then
  $F-\{e\} \subset P^I$.
\item Distinct feathers can only intersect at a vertex of $P$.
    \item The line segment connecting
    two consecutive feather tips lies in $P^I$.
    \end{enumerate}
\end{theorem}

When we apply $\Delta_k$ to $P$ we are just specifying the
points of $\Delta_k(P)$.  We define the
{\it polygon\/} $\Delta_k(P)$ so that the edges
are the bounded segments connecting the consecutive tips of
the feathers of $P$. With this definion, Statement 2 of
Theorem \ref{main} follows immediately from
Theorem \ref{feather}.

\subsection{The Proof}

There is one crucial idea in the proof of
Theorem \ref{feather}:  The soul of $P$ lies in the
sector $F^*$ opposite any of its feathers $F$.  See
Figure 4.2.

\begin{center}
\resizebox{!}{1.2in}{\includegraphics{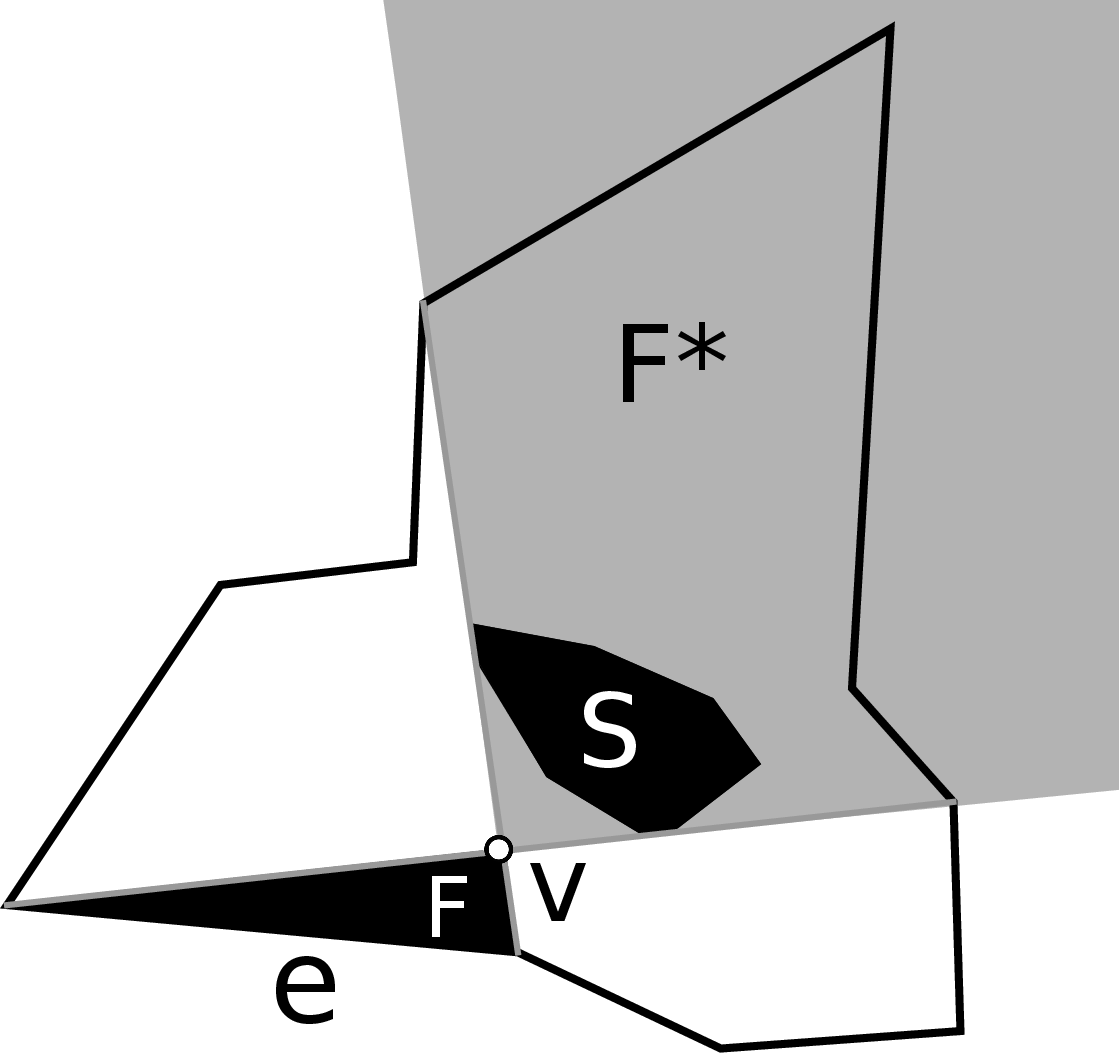}}
\newline
Figure 4.2 The soul lies in the sectors opposite the feathers.
\end{center}

We will give a homotopical proof of Theorem \ref{feather}.
By truncating our path of birds, we can assume that
Theorem \ref{feather} holds for all $t \in [0,1)$.
We set $P=P(1)$, etc.
\newline
\newline
{\bf Statement 1:\/}
  Figure 4.3 shows the $2$ ways that Statement 1 could
  fail:
  \begin{enumerate}
  \item The tip $v$ of the feather $F$ could coincide with some $p \in P$.
  \item Some $p \in P$ could lie in the interior point of $\partial F-e$.
  \end{enumerate}

\begin{center}
\resizebox{!}{2.1in}{\includegraphics{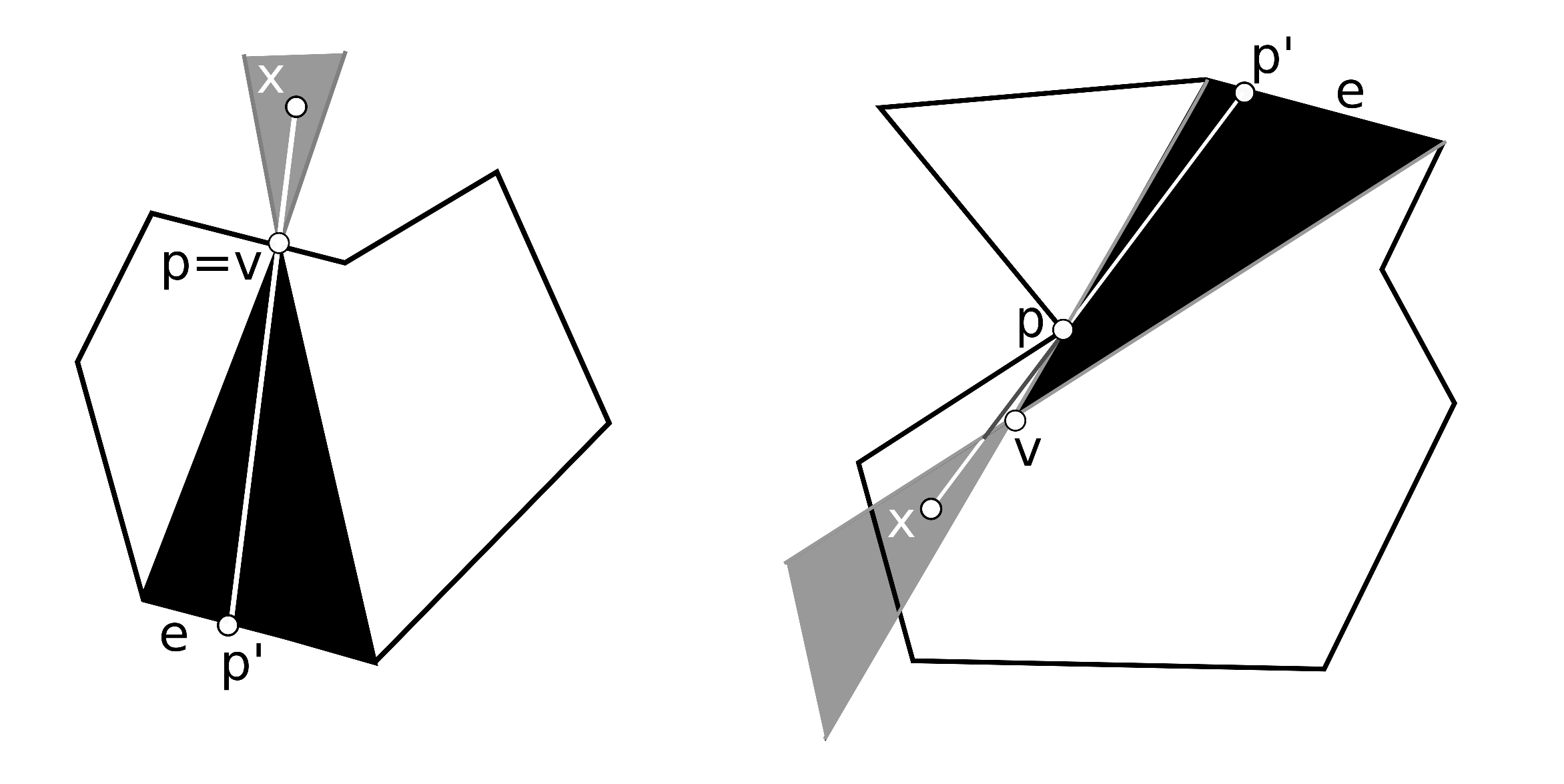}}
\newline
Figure 4.3: Case 1 (left) and Case 2 (right).
\end{center}

  For all $x \in F^*$, the ray $\overrightarrow{xp}$
  intersects $P$ both at $p$ and at a point $p' \in e$. This
  contradicts the fact that for any $x \in S \subset F^*$, the
  polygon $P$ is strictly star-shaped with respect to $x$.
  This establishes Statement 1 of Theorem \ref{feather}.
  \newline
  \newline
  {\bf Statement 2:\/}
  Let $F_1$ and $F_2$ be two feathers of $P$,
  having bases $e_1$ and $e_2$.  For Statement 2,
  it suffices to prove that
    $F_1-e_1$ and $F_2-e_2$ are disjoint.

  The same homotopical argument as for Statement 1 reduces us to the
case when $F_1$ and $F_2$ have disjoint interiors but
$\partial F_1-e_1$ and
$\partial F_2 -e_2$ share a common point $x$.
If $\partial F_1$ and $\partial F_2$ share an entire line segment then,
thanks to the fact that all the feathers are oriented
the same way, we would have two $(k+1)$ diagonals of
$P$ lying in the same line and having opposite
orientation. Lemma \ref{oppo} rules this out.

In particular $x$ must be the tip of at least one
feather.
Figure 4.4 shows the case when $x=v_1$, the tip of $F_1$,
but $x \not = v_2$.  The case when $x=v_1=v_2$ has
a similar treatment.

\begin{center}
\resizebox{!}{4in}{\includegraphics{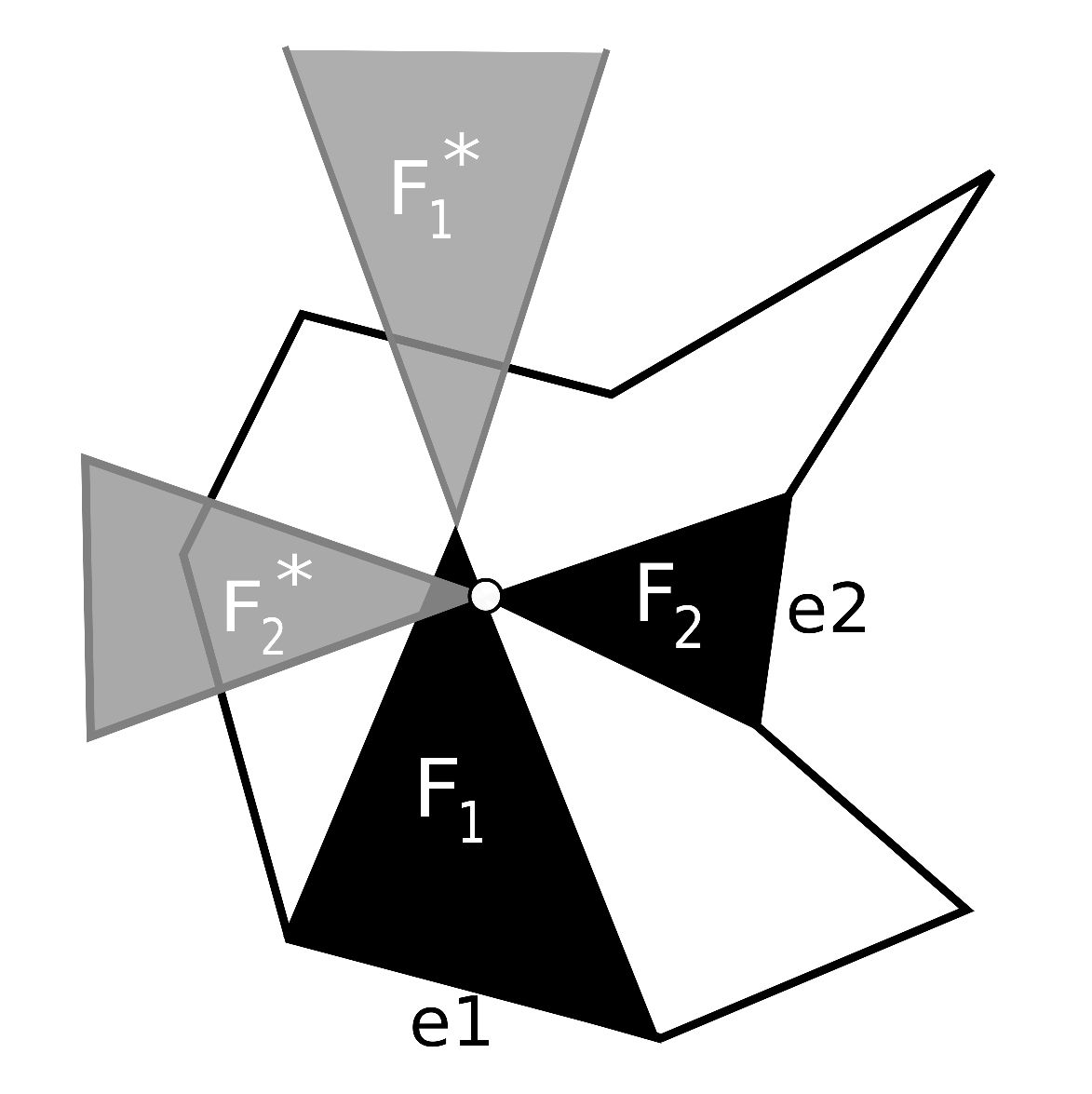}}
\newline
Figure 4.4: Opposiing sectors are disjoint
\end{center}

In this case, the two sectors $F_1^*$ and $F_2^*$
are either disjoint or intersect in a single point.
This contradicts the fact that $S \subset F_1^* \subset F_2^*$
has non-empty interior.  This contradiction
establishes Statement 2 of Theorem \ref{feather}.
\newline
\newline
{\bf Statement 3:\/}
Recall that $\widehat P=P \cup P^I$.
Let $F_1$ and $F_2$ be consecutive feathers with
bases $e_1$ and $e_2$ respeectively.
Let $f$ be the edge connecting the tips of $F_1$ and $F_2$.
Our homotopy idea reduces us to the case when
$f \subset \widehat P$ and $f \cap P \not = \emptyset$.
Figure 4.5 shows the situation.

\begin{center}
\resizebox{!}{3.5in}{\includegraphics{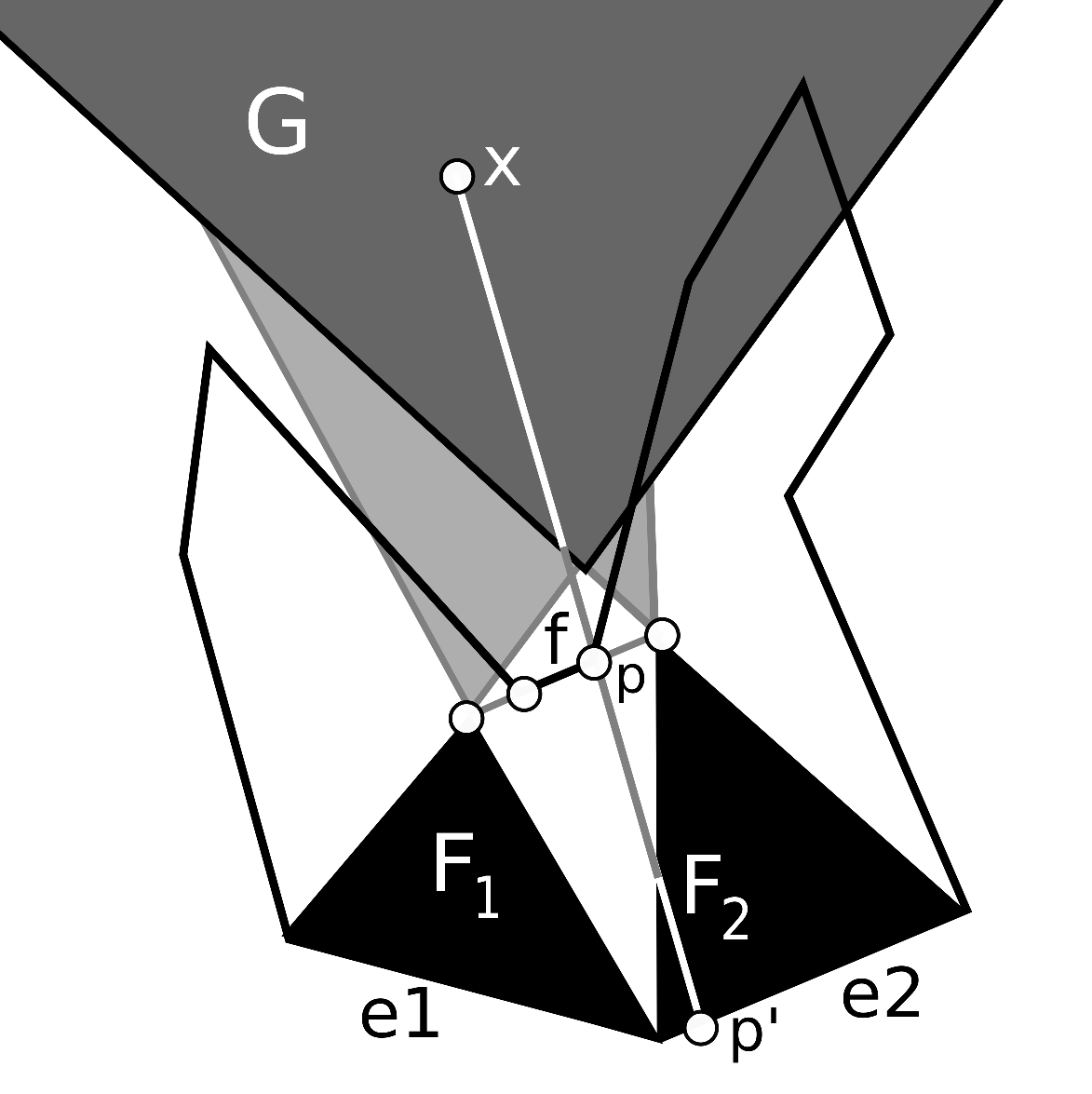}}
\newline
Figure 4.5: The problem a common boundary point
\end{center}

Note that $f \cap P$ must be strictly contained in
the interior of $f$ because (by Statement 1 of
Theorem \ref{feather}) the endpoints of $f$ lie
in $P^I$.
But then, $f \cap P$ is disjoint from $F_1^* \cap F_2^*$,
which is in turn contained in the shaded region $G$.
For any $x \in G$ and each vertex $p$ of $f$, the
ray the ray $\overrightarrow{xp}$ also intersects
$P$ at a point $p' \in e_1 \cup e_2$. This gives the
same contradiction as above when we take $x \in S \subset F_1^* \cap F_2^* \subset G$.
This completes the proof of Statement 3 of
Theorem \ref{feather}.

\newpage

\section{The Degeneration of Birds}
\label{DEGEN}

\subsection{Statement of Result}

Let $B_{k,n}$ denote the space
of $n$-gons which are $k$-birds.
Let $\chi_k$ denote the $k$-energy.
With the value of $k$ fixed in the background, we
say that a {\it degenerating path\/} is
a path $Q(t)$ of $n$-gons such that
\begin{enumerate}
  \item $Q(t)$ is planar for all $t \in [0,1]$. 
  \item All vertices of $Q(t)$ are distinct for all $t \in [0,1]$.
\item $Q(t) \in B_{k,n}$ for all $t \in [0,1)$ but $Q(1) \not \in B_{k,n}$.
\item $\chi_k(Q(t))>\epsilon_0>0$ for all $t \in [0,1]$.
\end{enumerate}

In this chapter we will prove the following result,
which will help us prove that
$\Delta_k(B_{k,n}) \subset B_{k,n}$ in the next chapter.
The reader should probably just use the
statement as a black box on the
first reading.

\begin{lemma}[Degeneration]
  If $Q(\cdot)$ is a degenerating path, then
  all but at most one vertex of $Q(1)$ lies in a
  line segment.
\end{lemma}

\noindent
{\bf Remark:\/} Our proof only uses the fact that
$Q$ has nontrivial edges,
nontrivial $k$-diagonals, and nontrivial $(k+1)$-diagonals.  Some
of the other vertices could coincide and it would
not matter.  Also, the same proof works if, instead of a
continuous path, we have a convergent
sequence $\{Q(t_n)\}$ with $t_n \to 1$ and a limiting
polygon $Q(1)=\lim Q(t_n)$.
\newline

\noindent
{\bf Example:\/}
Let us give an example for the case $k=1$ and $n=5$.
Figure 5.0 shows a picture of a pentagon $Q(t)$ for $t=1-s$.

\begin{center}
\resizebox{!}{1.1in}{\includegraphics{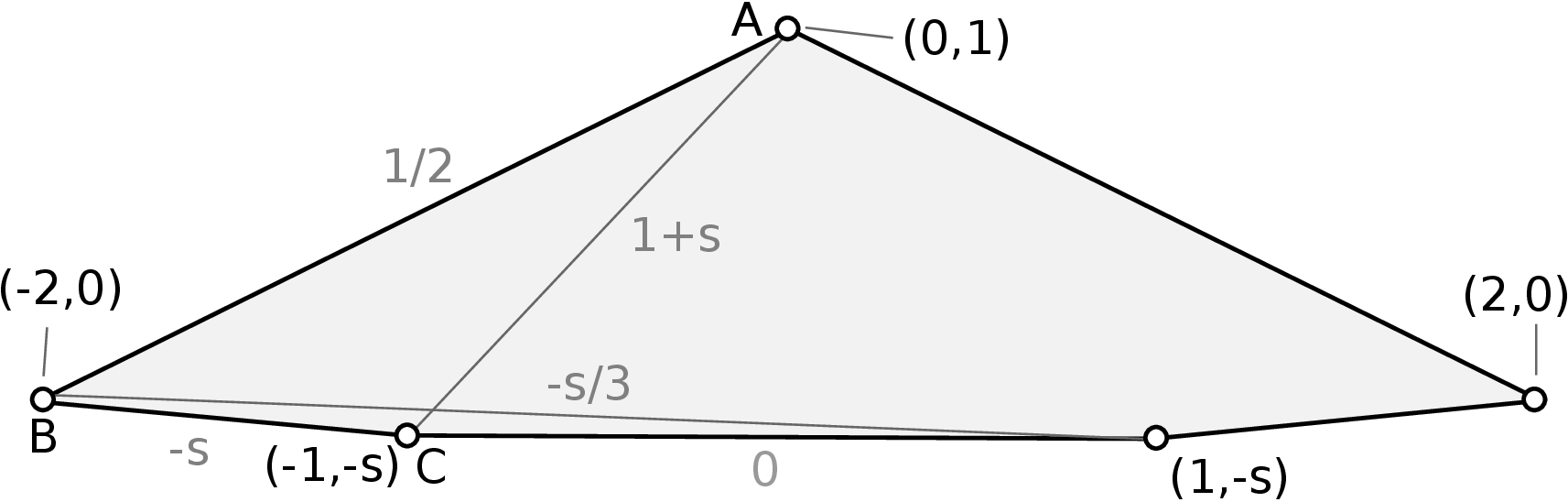}}
\newline
Figure 5.0: A degenerating path in the case $k=1$ and $n=5$.
\end{center}

Here $s$ ranges from $1$ to $0$ as $t$ ranges from
$0$ to $1$.    We have labeled some of the slopes to
facility the calculation (which we leave to the reader) that
$\chi_1(Q(t))$ remains uniformly bounded away from $0$.

        \subsection{Distinguished Diagonals}
        \label{RULES}

        We orient $Q(t)$ so that it
        goes counter-clockwise around the region it bounds.
        We orient the diagonal $Q_{ab}$ so that it points from $Q_a$ to $Q_b$.
        For $t<1$ the vertices
        $Q_1(t)$ and $Q_k(t)$ lie to the right of the diagonal $Q_{0,k+1}(t)$,
        in the sense that a person walking along this diagonal
        according to its orientation would see that points in the right.
        This has the same proof as Lemma \ref{L1}.
        The same rule holds for all cyclic relabelings of these points.
        The rule holds when $t<1$. Taking a
        limit, we get a weak version of the rule:  Each of
        $Q_1(1)$ and $Q_k(1)$ either lies to the right of the
        diagonal $Q_{0,k+1}(1)$ or on it. The same goes
        for cyclic relabeings. We call this the
        {\it Right Hand Rule\/}.

        Say that a {\it distinguished diagonal\/} of $Q(t)$ is
        either a $k$-diagonal or a $(k+1)$-diagonal.
        There
        are $2n$ of these, and they come in a natural cyclic order:
        \begin{equation}
        Q_{0,k}(t) \hskip 10 pt
        Q_{0,k+1}(t), \hskip  10 pt
        Q_{1,k+1}(t), \hskip 10 pt
        Q_{1,k+2}(t), ...
        \end{equation}
        The pattern alternates between $k$ and $(k+1)$-diagonals.
        We say that a {\it diagonal chain\/} is a consecutive list of
        these.

    We say that one oriented segment $L_2$ lies {\it ahead\/}
    of another one $L_1$ if we can rotate $L_1$ by
    $\theta \in (0,\pi)$ radians counter-clockwise to arrive at a segment parallel to $L_2$,
     In this case we
    write $L_1 \prec L_2$.  We have
        \begin{equation}
      \label{TURN}
      Q_{0,k+1}(t) \prec Q_{1,k+1}(t) \prec Q_{1,k+2}(t) \prec Q_{2,k+2}(t).
    \end{equation}
\begin{center}
\resizebox{!}{1.5in}{\includegraphics{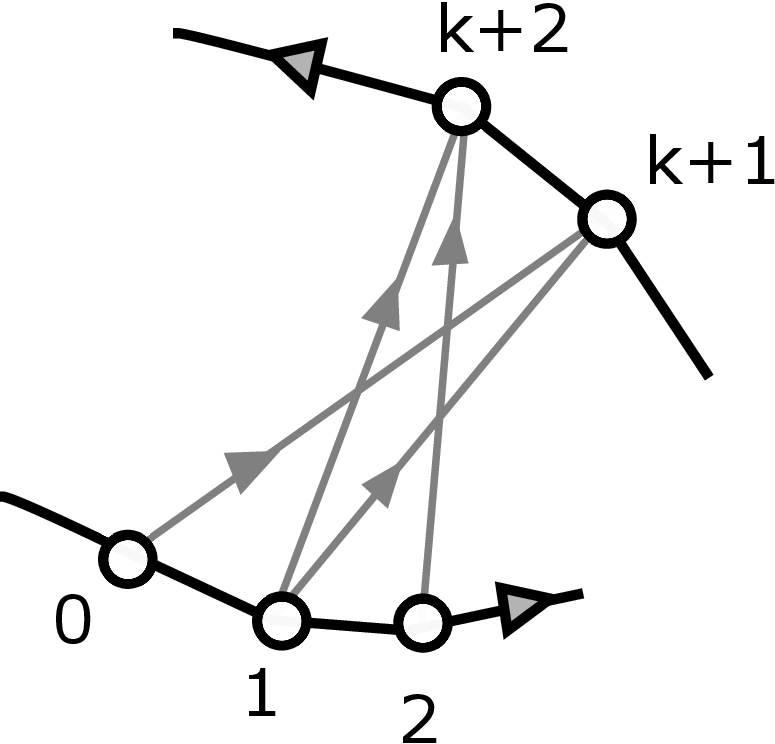}}
\newline
Figure 5.1: The turning rule
\end{center}
This certainly holds when $t=0$.
By continuity and the Right Hand Rule, this holds
for all $t<1$.
      Taking a limit, we see that the $k$-diagonals of
      $Q(1)$ weakly turn counter-clockwise in the sense
      that either $L_1 \prec L_2$ for consecutive diagonals
      or else $L_1$ and $L_2$ lie in the same line and
      point in the same direction.  Moreover, the total
      turning is $2\pi$.  If we start with one distinguished
      diagonal and move through the cycle then the turning
      angle increases by jumps in $[0,\pi]$ until it reaches $2\pi$.
      We call these observations {\it the Turning Rule\/}.

      \subsection{Collapsed Diagonals}
        \label{SUB}
        
        Figure 5.2 shows the distinguished diagonals incident to $Q_0$.
        We always take indices mod $n$. Thus $-k-1= n-k-1$ mod $n$.
        
     \begin{center}
     \resizebox{!}{1.7in}{\includegraphics{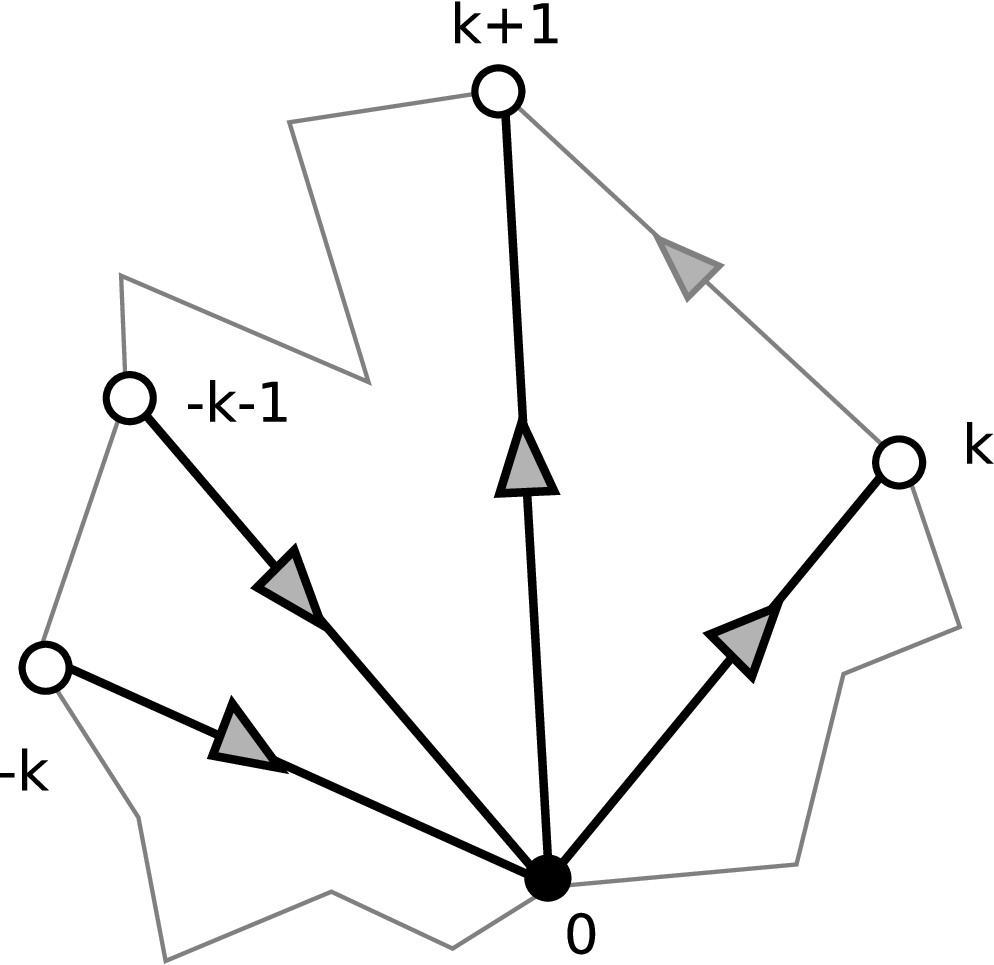}}
     \newline
     Figure 5.2: The $4$ distinguished diagonals incident to $Q_0(t)$.
     \end{center}

        We say that $Q$ has
        {\it collapsed diagonals\/} at a vertex if $Q$ if the
        $4$ distinguished diagonals incident to $Q_k$
        do not all lie on distinct lines.
        We set $Q=Q(1)$. 
       We set $X=X(1)$ for each object $X$ associated to $Q(1)$.

        Since $Q$ is planar but not $k$-nice, $Q$
        must have collapsed diagonals at some vertex.
        We relabel so that the collapsed diagonals are at $Q_0$.
     
     \begin{lemma}
       \label{dichotomy}
       If $Q$ has collapsed diagonals at $Q_0$ then
       $Q_{-k-1,0}$ and $Q_{0,k+1}$
       point in opposite directions or
       $Q_{-k,0}$ and $Q_{0,k}$ point in the same direction.
     \end{lemma}

     \startproof
     Associated to each diagonal incident to $Q_0$ is the ray which
     starts at $Q_0$ and goes in the direction of the other endpoint of
     the diagonal. (Warning: The ray may have the opposite orientation than
     the diagonal it corresponds to.)
     If the angle between any of the rays
     tends to $\pi$ as $t \to 1$ then the angle
     between the outer two rays tends to $\pi$.
     In this case $Q_{-k,0}$ and $Q_{0,k}$ point in the same directions.
     If the angle between non-adjacent rays tends to $0$ then
     $Q_{-k-1,0}$ and $Q_{0,k+1}$ are squeezed together and
     point in opposite directions.

     Suppose that the angle between adjacent rays tends to $0$.
     If the two adjacent rays are the middle ones, we have the
     case just considered.  Otherwise, either the angle between
     the two left rays tends to $0$ or the angle between the
     two right rays tends to $0$.  In either case, the
     uniform lower bound on the cross ratio forces a third
     diagonal to point either in the same or the opposite
     direction as these adjacent diagonals when $t=1$.
     Any situation like this leads to a case we have already considered.
     \endproof

     \subsection{The Case of Aligned Diagonals}
     \label{ALIGN}
     \label{ALIGNED}

     We say that $Q$ has {\it aligned diagonals\/} at the vertex $Q_0$
     if
     $Q_{-k,0}$ and $Q_{0,k}$ are parallel.  This is the second
     option in Lemma \ref{dichotomy}.  We make the same
     kind of definition at other vertices, with the indices shifted
     in the obvious way,.

     \begin{lemma}
       \label{conseq}
       \label{degen}
       Suppose $Q$ does not lie in a single line.
       Suppose also that $Q$
       has aligned diagonals at $Q_0$.  Then
       the diagonals $Q_{-k,0},Q_{-k,1},...,Q_{-1,k},Q_{0,k}$ all
       are parallel and (hence) the $2k+1$ points
        $Q_{-k},...,Q_0,...,Q_k$ are contained in the line defined by
        these diagonals.
     \end{lemma}

     \startproof 
    These two diagonals define a short chain of diagonals,
    which starts with the first listed diagonal and ends with the
    second one.
    They also define a long chain, which starts with the second and ends with the first.
     The total turning of the diagonals is $2\pi$, so one of the
     two chains defined by our diagonals turns $2\pi$ and the
     other turns $0$.   Suppose first that the long chain has $0$
     turning.  This chain involves all points
     of $Q$, and forces all points of $Q$ to be on the same line.
     So, the short chain must consist of parallel diagonals.
     \endproof

     All we use in the rest of the proof is that
     $Q_{-k},...,Q_k$ are all contained in a line $L$.
     By shifting our indices, we can assume
     that $Q_{k+1} \not \in L$.   This relabeling trick comes
     with a cost.   Now we cannot say whether the
     points $Q_{-k}....Q_k$ come in order on $L$.    We now
     regain this control.
    
    \begin{lemma}
      The length $2k$-diagonal chain
      $Q_{-k,0} \to ... \to Q_{0,k}$ consists entirely of parallel diagonals.
      There is no turning at all.
    \end{lemma}

    \startproof
    The diagonals $Q_{-k,0}$ and $Q_{0,k}$.
    are either parallel or anti-parallel.
    If they are anti-parallel, then the angle
    between the corresponding rays incident $Q_{0}(t)$ tends to $0$
    as $t \to 1$.  But these are the outer two rays.  This forces the angle between all
    $4$ rays incident to $Q_0(t)$ to tend to $0$.  The whole picture just folds up like
    a fan.  But one or these rays corresponds to $Q_{0,k+1}(t)$.  This
    picture forces $Q_{k+1} \in L$.  But this is not the case.

    Now we know that $Q_{-k,0}$ and $Q_{0,k}$ are parallel.
    All the diagonals in our chain are either parallel or anti-parallel to the
    first and last ones in the chain.  If we ever get an anti-parallel pair,
    then the diagonals in the chain must turn $2\pi$ around.  But then
    none of the other distinguished diagonals outside our chain turns
    at all.  As in Lemma \ref{degen}, this gives
    $Q \subset L$, which is false.
    \endproof

    We rotate the picture so that $L$ coincides with the $X$-axis and so that
    $Q_{0,k}$ points in the positive direction.   Since we are already using the
    words {\it left\/} and {\it right\/} for another purpose, we say that
    $p \in L$ is {\it forward of\/} of $q \in L$ if $p$ has larger $X$-coordinate.
    Likewise we say that $q$ is {\it backwards of\/} $p$ in this situation.
    We say that $Q_{0,k}$ points {\it forwards\/}.  We have
    established that $Q_{-k,0},...,Q_{0,k}$ all point forwards.
        
    \begin{lemma}
      $Q_{k+2} \in L$ and both
      $Q_{1,k+2}$ and $Q_{2,k+2}$ point backwards.
        \end{lemma}

        \startproof
        We have arranged that $Q_{k+1} \not \in L$.
    Let us first justify the fact that $Q_{k+1}$ lies above $L$.
    This follows from Right Hand Rule applied to
    $Q_{0,k+1}$ and $Q_k$ and the fact that $Q_{0,k}$ points forwards.
    Since $Q_{-k}, Q_{-k+1}, Q_1$ are collinear,
    $Q$ has collapsed diagonals at $Q_1$.  But $Q$ cannot have aligned
    diagonals because $Q_{1,k+1}$ is not parallel to $Q_{-k,1}$. Hence
    $Q$ has folded diagonals at $1$.    This means that the diagonals
    $Q_{-k,1}$ and
    $Q_{1,k+2}$ point in opposite directions.  This forces $Q_{k+2}
    \in L$ and morever we can say that $Q_{1,k+2}$ points backwards.

    We have $Q_2 \in L$ because $2 \leq k$.
    We want to see that $Q_{2,k+2}$ points forwards and they
    Suppose not.
   We consider the  $3$ distinguished diagonals
  $$Q_{0,k}, \hskip 10 pt Q_{1,k+2}, \hskip 10 pt Q_{2,k+2}.$$
  These diagonals
  respectively point forwards, backwards, forwards and they all
  point one direction or the other along $L$.
  But then, in going from $Q_{0,k}$ to $Q_{2,k+2}$, the diagonals have
  already turned $2\pi$.  Since the total turn is $2\pi$, the diagonals
    $Q_{2,k+2},\ Q_{3,k+3},...,Q_{n,n+k}$ are all parallel.  But then
  $Q_2,...,Q_{n} \in L$. This contradicts the fact that $Q_{k+1} \not \in L$.
  \endproof

\begin{lemma}
  \label{either}
  For at least one of the two indices $j \in \{2k+2,2k+3\}$ we have
  $Q_j \in L$ and $Q_{k+2,j}$ points forwards.
\end{lemma}

\startproof
Since $Q_1,Q_2,Q_{k+2}$ are collinear, $Q$ has collapsed diagonals at $Q_{k+2}$.
So, by Lemma \ref{dichotomy},
we either have folded diagonals at $Q_{k+2}$ or aligned diagonals
at $Q_{k+2}$.  The aligned case gives $Q_{2k+2} \in L$ and the folded
case gives $Q_{2k+3} \in L$.
We need to work out the direction of pointing in each case.

Consider the aligned case.
Suppose $Q_{k+2,2k+2}$ points backwards, as shown in Figure 5.3.
\begin{center}
\resizebox{!}{1in}{\includegraphics{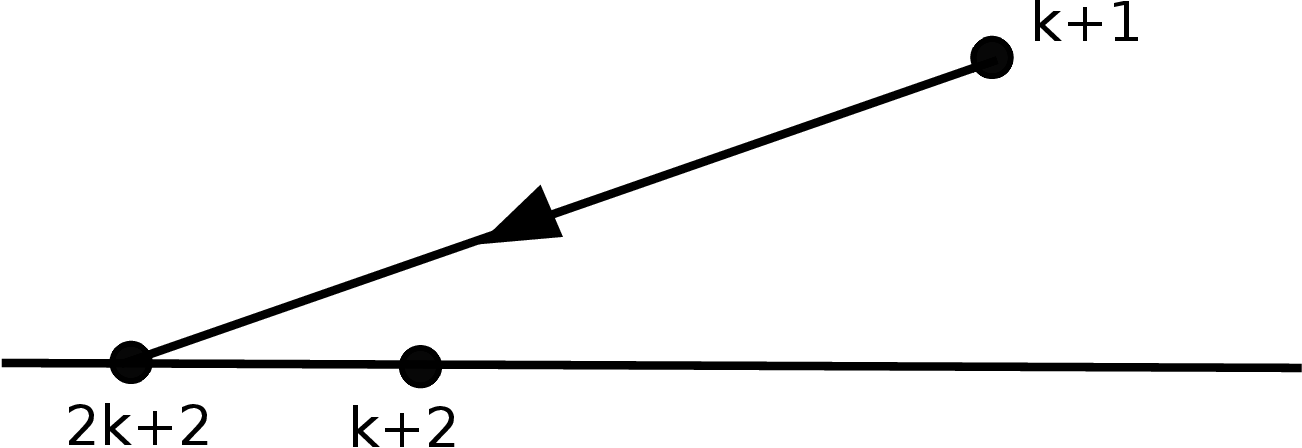}}
\newline
Figure 5.3:  Violation of the Right Hand Rule
\end{center}

This violates the Right Hand Rule for $Q_{k+2}$ and $Q_{k+1,2k+2}$
because $Q_{k+1}$ lies above $L$.  

Consider the folded case.  Since $Q_{k+2,2k+3}$ and
$Q_{1,k+2}$ point in opposite directions, and
$Q_{1,k+2}$ points backwards (by the previous lemma), $Q_{k+2,2k+3}$ points forwards.
\endproof

        Let $j \in \{2k+2,2k+3\}$ be the index from Lemma \ref{either}.
        Consider the $3$ diagonals
        $$Q_{0,k}, \hskip 10 pt Q_{1,k+1}, \hskip 10 pt Q_{k+2,j}.$$
        These diagonals are all parallel to $L$ and respectively
        point in the forwards, backwards, forwards direction.
        This means that
        the diagonals in the chain  $Q_{0,k} \to ... \to Q_{k+2,j}$
        have already
        turned $2\pi$ radians.  But this means that the diagonals
        $$Q_{k+2,2k+3},\hskip 10 pt Q_{k+3,2k+3}, \hskip 10 pt
        Q_{k+3,2k+4},\hskip 10 pt ... \hskip 10 pt Q_{0,k}=Q_{n,n+k}$$
        are all parallel and point forwards along $L$.
        Hence $Q_{k+2}, Q_{k+3}, ... , Q_{n} \in L$.
        Hence all points but $Q_{k+1}$ lie in $L$.

        \subsection{The Case of Double Folded Diagonals}

        We fix a value of $k$.
        Say that two indices $a,b \in \Z/n$ are {\it far\/}
        if their distance is at least $k$ in $\Z/n$.
        We say that $Q$ has {\it far folded diagonals\/}
        if $Q$ has folded diagonals at $Q_a$ and
        $Q$ has folded diagonals at $b$ and $a,b$ are far.

        In this case we have two parallel diagonals
        $Q_{a,a+k+1}$ and $Q_{b,b+k+1}$.  As in the
        proof of Lemma \ref{degen}, one of the two
        diagonal chains defined by these diagonals
        consists of parallel diagonals.  The far condition
        guarantees that at least $2k+1$ consecutive points
        are involved in each chain. But then we get
        $2k+1$ collinear points.   So, if $Q$ has
        far folded diagonals, then the same proof as
        in the previous section shows that the conclusion
        of the Degeneration Lemma holds for $Q$.

     \subsection{Good Folded Diagonals}

     We say that the folded diagonals
     $Q_{-k-1,0}$ and $Q_{0,k+1}$ are {\it good\/}
     if all the points
     $Q_{k+1}, Q_{k+2},...,Q_{n-k-1}$ are collinear.
     This notion is empty when $k=2$ and $n=7$ but
     otherwise it has content.  In this section
     we treat the case when $Q$ has a pair of
     good folded diagonals.    We start by
     discussing an auxiliary notion.

        We say that $Q$ has
        {\it backtracked edges\/} at $Q_a$
        if the angle between the
        edges $Q_{a,a+1}$ and $Q_{a,a-1}$ is either
        $0$ or $2\pi$.

        \begin{lemma}
          \label{nofold}
          If $Q$ has backtracked edges at $Q_a$
          then $Q$ has folded diagonals at $Q_a$.
        \end{lemma}

        \startproof
        For $t \in [0,1)$, the edges of $Q$ emanating from $a$ divide the plane
        into $4$ sectors, and one of these sectors, $C(t)$ contains
        all the distinguished diagonals emanating from $Q_a(t)$.   The
        sector $C(t)$ is the one which locally intersects $Q(t)$ near
        $Q_a(t)$.
        The angle of $C(t)$ tends to $0$ as $t \to 1$, forcing all the
        distinguished diagonals emanating from $Q_a(t)$ to squeeze
        down as $t \to 1$.  This gives us the folded diagonals.
        \endproof

        We will use Lemma \ref{nofold} in our analysis of good folded
        edges.  Now we get to it.
      We rotate so that our two diagonals are
     in the line $L$, which is the $X$-axis.
     We normalize so that $Q_0$ is the origin, and
     $Q_{k+1}$ and $Q_{-k-1}$ are forward of $Q_0$.
     
     \begin{lemma}
       \label{goodfold}
       If $n>3k+1$ and $Q_{-k-1,0}, Q_{0,k+1}$ are good folded
       diagonals, then  the Degeneration Lemma is true for $Q$.
\end{lemma}

\startproof
Suppose first that $Q_1 \in L$.  Then $Q$ has
folded diagonals at $Q_{k+1}$.
When $n>3k+1$ the indices
$(k+1)$ and $(-k-1)$ are $k$-far.
This gives $Q$ far folded diagonals,
a case we have already treated.

To finish our proof, we
show that $Q_1 \in L$.
We explore some of the other points.
We know that $Q_{k+1},...,Q_{n-k-1} \in L$.
We can relabel dihedrally so that $Q_{n-k-1}$ is
forwards of $Q_{k+1}$.
We claim that $Q_{k+2}$ is forwards of $Q_{k+1}$.
Suppose not.
Then there is some index
$a \in \{k+2,...,-k-2\}$ such that
$Q_a$ is backwards of $Q_{a \pm 1}$.
What is going on is that our points
would start by moving backwards on $L$ and
eventually they have to turn around.  The index
$a$ is the turn-around index.
This gives us backtracked edges at $Q_a$.
By Lemma \ref{nofold}, we have
folded diagonals at $Q_a$.  But $a$
and $0$ are $k$-far indices.
This gives $Q$ far-folded diagonals.

The only way out of the contradiction is
that $Q_{k+2}$ is forwards of $Q_{k+1}$.

\begin{center}
\resizebox{!}{.8in}{\includegraphics{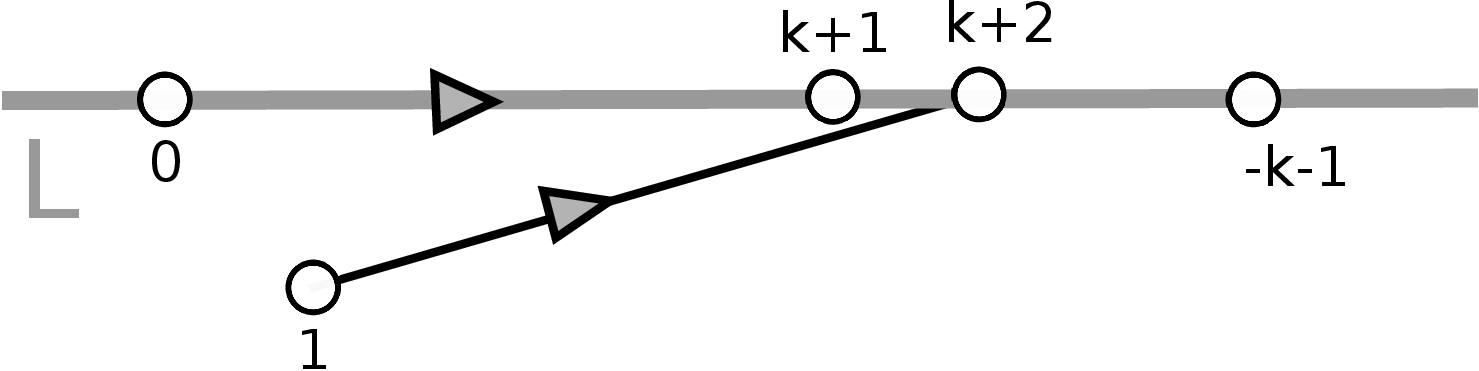}}
\newline
Figure 5.4: A contradiction involving $Q_1$.
\end{center}

Suppose $Q_1 \not \in L$.
by the Right Hand Rule applied to the
diagonal $Q_{0,k+1}$, the point
$Q_1$ lies beneath $L$, as shown in
Figure 5.4. But then $Q_{k+1}$ lies to the left of the diagonal
$Q_{1,k+2}$. This violates the Right Hand Rule.
Now we know that $Q_1 \in L$.
\endproof

\begin{lemma}
  \label{special00}
       Suppose $n=3k+1$ and $k>2$.  If $Q_{-k-1,0}, Q_{0,k+1}$ are good folded
       diagonals, then  the Degeneration Lemma is true for $Q$.
     \end{lemma}

     \startproof
     The same argument as in Lemma \ref{goodfold}
establishes the key containment
$Q_1 \in L$.  (We need $k>2$ for this.)
From here, as in Lemma \ref{goodfold}, we deduce that
$Q_{-k-1,0}$ and $Q_{k+1,2k+2}$ are parallel.
This time the conclusion we get from this
is not as good.
We get a run of $k$ points in $L$, and
then a point not necessarily in $L$, and then
a run of $k$ additional points in $L$.

The points are $Q_{k+1},...,Q_{2k+1},...,Q_0$
with the point $Q_{-k}$ omitted.
But then $Q$ has folded diagonals at each of
these points except the outer two,
$Q_{k+1}$ and $Q_0$.   But then
For each such index $h$, we see that both
$Q_{h \pm (k+1)}$ belong to $L$.  This gives us
all but one point in $L$.

It is instructive to consider an example,
say $k=4$ and $n=13$.  In this case, our
initial run of points in $L$ is
$Q_5,Q_6,Q_7,Q_8,Q_{10},Q_{11},Q_{12},Q_{13}.$
The folded diagonals at $Q_6, Q_7, Q_8$ respectively
give $Q_1,Q_2,Q_3 \in L$.  The folded diagonals at
$Q_{10},Q_{11},Q_{12}$ respectively give
$Q_2,Q_3,Q_4 \in L$.
\endproof

Finally we consider the case $k=2$ and $n=7$.
In this case all we know is that
$Q_0, Q_3,Q_4 \in L$ with $Q_3,Q_4$ forwards of $Q_0$.
We can dihedrally relabel to that $Q_4$ is forwards of $Q_3$.
Here $Q_3=Q_{k+1}$ and $Q_4=Q_{k+2}$.  So, now we
can run the same argument as in Lemma \ref{special00}
to conclude that $Q_1 \in L$.   Now we proceed as
in the proof of Lemma \ref{special00}.

\subsection{Properties of the Soul}
\label{souldef}

        We define $S=S(1)$ to be the set of all accumulation points
        of sequences $\{p(t_n)\}$ where $p(t_n) \in S(t_n)$ and
        $t_n \to 1$.  Here $S(t_n)$ is the soul of $P(t_n)$.
        We have one more case to analyze, namely ungood folded
        diagonals. 
        To make our argument go smoothly, we first prove some
        properties about $S$.

        \begin{lemma}
          Suppose that $Q$ has folded diagonals at $Q_0$.
          If the Degeneration Lemma is false for $Q$, then $S$ is
          contained
          in the line segment joining $Q_0$ to $Q_{k+1}$
          \end{lemma}

          \startproof 
        Here is a general statement about $S$.
       Since $S(t)$ is non-empty and closed for all $t<1$, we see by compactness
        that $S$ is also a non-empty closed subset of the closed region
        bounded by $Q$.
        By continuity  $S$ lies to the left of all the closed half-planes
          defined by the oriented $(k+1)$ diagonals (or in their boundaries).
         Since $S$ lies
        to the left of (or on) each $(k+1)$ diagonal,
        $S$ is a subset of the line $L$ common to the
        folded diagonals and indeed $S$ lies to one
        side of the fold point $Q_0$.
        From the way we have normalized, $S$ lies in the
        $X$-axis forward of $Q_0$.  (The point $Q_0$ might
        be an endpoint of $S$.)
        
        If $S$ contains points of $L$ that lie forward of
        $Q_{k+1}$ then either the diagonal
        $Q_{k+1,2k+2}$ points along the positive
        $X$-axis or into the lower half-plane.
        In the former cases, the
        diagonals $Q_{0,k+1},Q_{k+1,2k+2}$ are
        parallel and we get at least $2k+1$ collinear
        points and so the Degeneration Lemma holds for $Q$.

        If $Q_{k+1,2k+2}$ points into the negative half-plane,
        then the diagonal $Q_{0,k+1}$ turns more than
        $\pi$ degrees before reaching $Q_{k+1,2k+2}$.
        But then the diagonals
        in the chain $Q_{-k-1,0} \to ... \to Q_{0,k+1} ... \to
        Q_{k+1,2k+2}$
        turn more than $2\pi$ degrees.  This is a contradiction.
        \endproof

        \noindent
        {\bf Remark:\/} The same argument works with $Q_{-k-1}$ in
        place of $Q_{k+1}$.
        
        \begin{lemma}
          \label{soul1}
          If the Degeneration Lemma is false for $Q$ then $S$ cannot
          intersect $Q$ in the interior of an edge of $Q$.
         \end{lemma}

         \startproof
         Suppose this happens.  We relabel so that the edge is
         $Q_{0,1}$.
         By the Right Hand Rule, the point $Q_1$ is not on the left of
         the diagonal $Q_{0,k+1}$.  At the same time, $S$ is not on
         the
         right of the diagonal.  The only possibility is that
         $Q_1, Q_0, Q_{k+1}$ are collinear.   Likewise
         $Q_{-k},Q_0,Q_1$ are collinear.  Furtheremore,
         the $(k+1)$-diagonals $Q_{-k,1}$ and $Q_{0,k+1}$ are
         parallel.
         Figure 5.5 shows the situation
         for $Q(t)$ and $S(t)$ when $t$ is very near $1$.

       \begin{center}
       \resizebox{!}{.6in}{\includegraphics{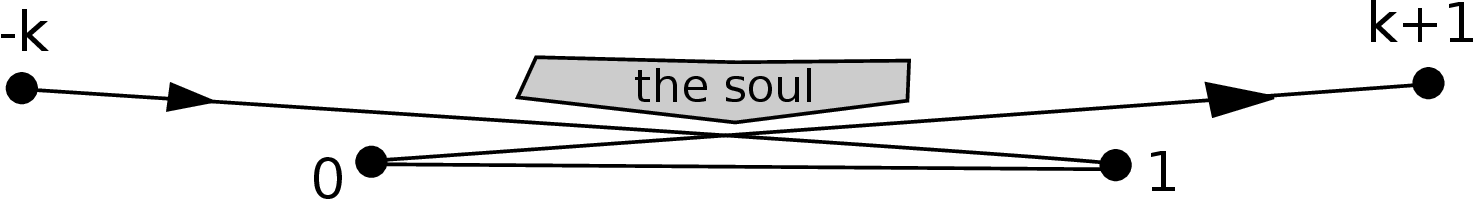}}
       \newline
       Figure 5.5:  The relevant points and lines.
     \end{center}

       But now we have two $(k+1)$-diagonals that are parallel and
       which start at indices that are $k$ apart in $\Z/n$.  This
       gives us $2k+1$ consecutive collinear points on the line
       containing our edge.   We know how to finish the Degeneration
       Lemma in this case.   The only way out is that $S$ cannot
       intersect $Q$ in the interior of an edge of $Q$.
       \endproof
     
       \begin{lemma}
         If the Degeneration Lemma is false for $Q$, then $S$ cannot
         contain a vertex of $Q$.
         \end{lemma}

        \startproof
        We relabel so that $Q_0 \in S$.
        The same analysis as in the previous
        lemma shows that   $Q_1, Q_0, Q_{-k}$
        are collinear.  Figure 5.6. shows the situation
        for $t$ near $1$. 
        At the same time, the points $Q_{-1}, Q_0, Q_k$ are collinear.

       \begin{center}
       \resizebox{!}{.6in}{\includegraphics{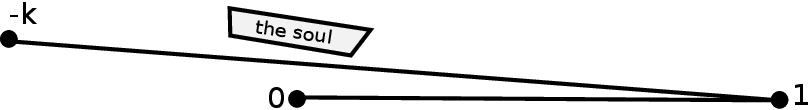}}
       \newline
       Figure 5.6:  The relevant points and lines
        \end{center}
           
        To avoid a case of the Degeneration Lemma we have already
        done, $Q$ must have folded diagonals at $Q_{-k}$. Likewise
        $Q$ must have folded diagonals at $Q_k$.
       But then $Q$ has far folded diagonals, and the
       Degeneration Lemma holds for $Q$.
       \endproof

        Now let us bring back our assumptions:  $Q$ has folded
        diagonals
        at $Q_0$ and the points $Q_0,Q_{k+1},Q_{-k-1}$ all lie in
        the $X$-axis in the forward order listed.

        \begin{corollary}
          \label{pointX}
          If the Degeneration Lemma is false for $Q$ then $S$ lies in
          the
          open interval bounded by $Q_0$ and $Q_{k+1}$ and no point of
          $S$ lies in $Q$.  In particular, $S$ contains a point $x$,
          forwards of $Q_0$ and
          backwards of both $Q_{k+1}$ and $Q_{-k-1}$, that
          is disjoint from $Q$.
          \end{corollary}

     \subsection{Ungood Folded Diagonals}
     \label{ungood}

     The only case left is when $Q$ does not have $2k+1$ consecutive
     collinear points, and when all folded diagonals of $Q$ are
     ungood.    Without loss of generality, we will consider the case
     when $Q$ has ungood folded diagonals at $Q_0$. 
     We normalize as in the previous section, so that     
    $Q_{0}, Q_{k+1}, Q_{-k-1}$ lie in forward order on
    $L$, which is the $X$-axis.  Let $x$ be a point from
    Corollary \ref{pointX}.

We call an edge of $Q$ {\it escaping\/} if $e \cap L$ is a single point.
We call two different edges of $Q$, in the labeled sense,
{\it twinned\/} if they are both escaping and
if they intersect in an open interval.  Even if two distinctly labeled edges
of $Q$ coincide, we consider them different as labeled edges.

\begin{lemma}
  $Q$ cannot have twinned escaping edges.
\end{lemma}

\startproof
Consider $Q(t)$ for $t$ near $1$.  This polygon is
strictly star shaped with respect to a point $x(t)$ near $x$.

\begin{center}
\resizebox{!}{1.2in}{\includegraphics{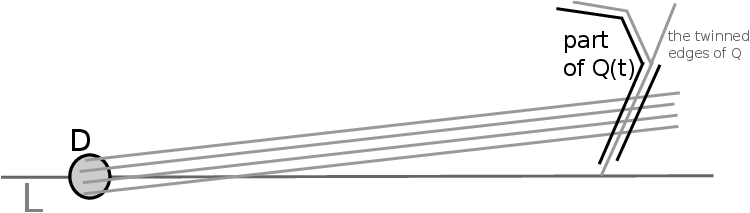}}
\newline
Figure 5.7: Rays intersecting the twinned segments
\end{center}

There is a disk $D$ about $x$ such that every $p \in D$ contains
a ray which intersects the twinned edges in the
middle third portion of their intersection.
Figure 5.7 shows what we mean.
Once $t$ is sufficiently near $1$, the soul $S(t)$ will
intersect $D$, and for all points $p \in D$ there will be
a ray which intersects
$Q(t)$ twice.  This contradicts the fact that
$Q(t)$ is strictly star-shaped with respect to all points of $S(t)$.
\endproof

We say that an escape edge {\it rises above\/} $L$ if it
intersects the upper half plane in a segment.

\begin{lemma}
  \label{rise}
  $Q$ cannot have two escape edges which rise above $L$ and
  intersect $Q$ on the same side of the point $x$.
\end{lemma}

\startproof
This situation is similar to the previous proof.
In this case, there is a small disk $D$ about $x$
such that every point $p \in D$ has a ray which
intersects both rising escape edges transversely,
and in the middle third of each of the two
subsegments of these escape edges that lie above $L$.
Figure 5.8 shows this situation.

\begin{center}
\resizebox{!}{1.4in}{\includegraphics{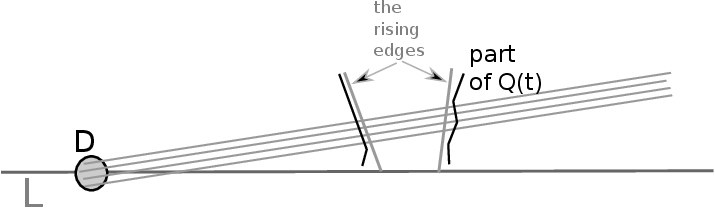}}
\newline
Figuren 5.8: Rays intersecting the rising segments.
\end{center}

In this case, some part of $Q(t)$ closely shadows
our two escape edges for $t$ near $1$.  But then,
once $t$ is sufficiently near $1$, each ray we
have been talking about intersects
$Q(t)$ at least twice, once by each escaping edge.
This gives the same contradiction as in the
previous lemma.
\endproof

We define falling escape segments the same way.
The same statement as in Lemma \ref{rise} works for
falling escape segments.  Since $x \not \in Q$
we conclude that $Q$ can have at most $4$ escaping
segments total.

But $Q=Q_+ \cup Q_-$, where $Q_{\pm}$ is an arc
of $Q$ that starts at $Q_{k+1}$ and ends
at $Q_{-k-1}$.  Since both these arcs
start and end on $L$, and since both
do not remain entirely on $L$, we see that
each arc has at least $2$ escape edges, 
and none of these are twinned.
This means that both $Q_+$ and $Q_-$ have
exactly two escape edges.

Now for the moment of truth:  Consider $Q_+$.
Since $Q_+$ just has $2$ escape edges, they
both have to be either rising or falling.
Also, since $Q_+$ starts and ends on the
same side of $x$, and cannot intersect $x$,
both the escape edges for
$Q_+$ are on the same side of $x$.   This is a
contradiction.  The same argument would
work for $Q_-$ but we don't need to make it.

\newpage

\section{The Persistence of Birds}

In this chapter we prove Statement 3 of
Theorem \ref{main}, namely the fact that
$\Delta_k(B_{n,k})=B_{n,k}$.
First we use the Degeneration Lemma to prove that
$\Delta_k(B_{n,k}) \subset B_{n,k}$.  Then we
deduce the opposite containment from
projective duality and from the factoring of
$\Delta_k$ given in \S \ref{duality}.

\subsection{Containment}

Suppose for the sake of contradiction that
there is some $P \in B_{k,n}$ such that
$\Delta(P) \not \in B_{k,n}$.
Recall that there is a continuous path
$P(t)$ for $t \in [0,1]$ such that
$P(0)$ is the regular $n$-gon.

Define $Q(t)=\Delta_k(P(t))$.
There is some $t_0 \in [0,1]$ such that
$Q(t_0) \not \in B_{k,n}$.  We can truncate
our path so that $t_0=1$.   In other words,
$Q(t) \in B_{n,k}$ for $t \in [0,1)$ but
  $Q(1) \not \in B_{k,n}$.

  \begin{lemma}
    $Q(\cdot)$ is a degenerating path.
  \end{lemma}

  \startproof
  Note that $Q(\cdot)$ is planar and hence satisfies Property 1 for
  degenerating paths.  Let $P=P(1)$ and $Q=Q(1)$.
  If $Q$ doe not have all distinct vertices then
  two different feathers of $P$ intersect at a point
  which (by Statement 2 of Theorem \ref{main}) lies in
  $P^I$.    This contradicts
  Statement 2 of Theorem \ref{feather}.
  Hence $Q(\cdot)$ satisfies Property 2 for
  degenerating paths.
  By construction, $Q(t) \in B_{n,k}$ for
  all $t \in [0,1)$.  Hence $Q(\cdot)$ satisfies
  Property 3.
  The energy $\chi_k$ is
  well-defined and continuous on $B_{k,n}$.
  Hence, by compactness, $\chi_k(P(t))>\epsilon_0$
  for some $\epsilon_0>0$ and all $t \in [0,1]$.
Now for the crucial step:  We have already proved that
$\chi_k \circ \Delta_k=\chi_k$.
  Hence
  $\chi_k(Q(t))>\epsilon_0$ for all $t \in [0,1]$.  That is,
  $Q(\cdot)$ satisfies Property 4 for
  degenerating paths.
  \endproof

Now we apply the Degeneration Lemma to $Q(\cdot)$.  
  We conclude that all but at most $1$ vertex of $Q(1)$ lies in a line $L$.
  Stating this in terms of $P(1)$, we can say that all but at most one of the
  feathers of $P(1)$ have their tips in a single line $L$. Call an edge of $P(1)$ {\it ordinary\/} if the
    feather associated to it has its tip in $L$.  We call the
    remaining edge, if there is one, {\it special\/}.
    Thus, all but at most one edge of $P$ is ordinary.
  
Let $S(t)$ be the soul of $P(t)$.  We know that $S(1)$ has
non-empty interior by Theorem \ref{soul}.   For ease of
notation we set $P=P(1)$ and $S=S(1)$.

\begin{lemma}
  \label{one}
      $P$ cannot have ordinary edges $e_1$ and $e_2$ that
      lie on opposite sides of $L$ and are disjoint from $L$.
    \end{lemma}

    \startproof
    Suppose this happens.  Figure 6.1 shows the
    situation.
    
       \begin{center}
       \resizebox{!}{.81in}{\includegraphics{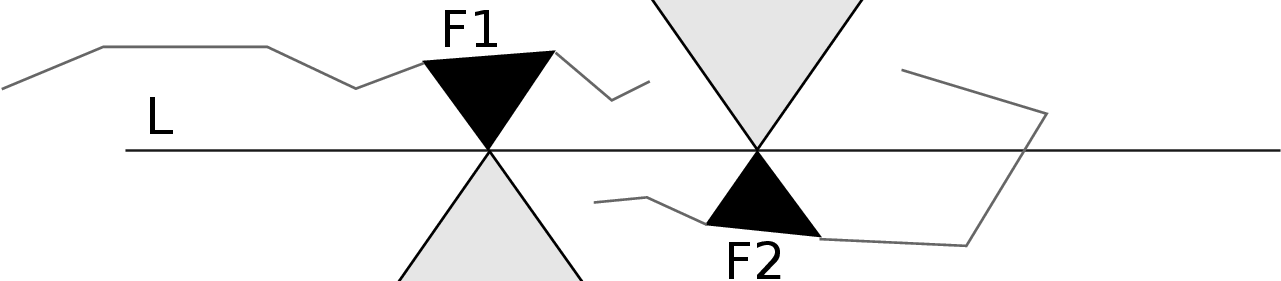}}
       \newline
       Figure 6.1:  Two feathers on opposite sides of $L$.
     \end{center}
    
    Let $F_1$ and $F_2$ be the two associated feathers.
    Then the opposite sector $F_1^{*}$ lies above $L$, and the opposite
    sector $F_2^{*}$ lies below $L$
    and the two tips are distinct.  But then $S(1)$, which must
    lie in the intersection of these sectors, is empty.
    \endproof

    \begin{lemma}
      \label{two}
      $P$ cannot have more than $2$ ordinary edges which intersect $L$.
    \end{lemma}

    \startproof
    Note that an ordinary edge
    cannot lie in $L$ because then the tip would not.
    So, an ordinary edge that intersects $L$ does so either
    at a single vertex or at an interior point.
    As we trace along $L$ in one direction or the other we encounter
    the first intersecting edge and then the last one and then some
    other intersecting edge.  Let $F_1.F_2.F_3$ be the two feathers,
    as shown in Figure 6.3.   Let $e_j$ be the edge of $F_j$ that
    belongs  to $P$.  Let $v_j$ be the tip of $F_j$.  (Figure 6.3
    shows the case when
    $e_j \cap L$ is an interior point of $e_j$ for each $j=1,2,3$,
   but the same argument
    would
    work if some of these intersection points were vertices.)
    
       \begin{center}
       \resizebox{!}{.8in}{\includegraphics{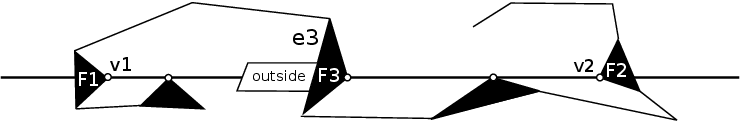}}
       \newline
       Figure 6.2:  Three or more crossing edges
     \end{center}
    
     One of the two arcs $\alpha$ of $Q$ joining $v_1$ to $v_2$ stays in $L$,
     namely the one avoiding the one point of $Q$ not on $L$.
     However, $\alpha$  passes right through $F_3$ and in particular
     crosses $e_3$ transversely. However, one side of $F_3$ is
     outside $P$.  Hence $\alpha$ is not contained in $P^I$,
     the interior of the region bounded by $P$.
     This contradicts Statement 2 of Theorem \ref{main}, which says
    that $Q \subset P^I$.
     \endproof

     The line $L$ divides the plane into two open half-planes,
     which we call the {\it sides\/} of $L$.
     Lemma \ref{one} says that  $P$ cannot have ordinary edges contained in
     opposite sides of $L$.   Lemma \ref{two} says that
     at most $2$ ordinary edges can intersect $L$.
 Hence, all but at most $2$ of the ordinary edges of
 $P$ lie on one side of $L$. Call this the
 {\it abundant side\/}  of $L$.   Call the other
 side the {\it barren side\/}.   The barren side contains
 no ordinary edges at all, and perhaps the special edge.
 In particular, at most two vertices of $P$ lie in the barren side.

       \begin{center}
       \resizebox{!}{1.3in}{\includegraphics{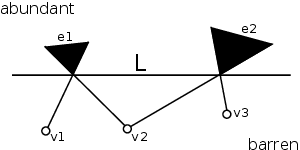}}
       \newline
       Figure 6.3:  Following the diagonals bounding a feather
     \end{center}
    
 At the same time, each ordinary edge on the abundant side
 contributes two vertices to the barren side:  We just follow
 the diagonals comprising the corresponding feather. These
 diagonals cross $L$ from the abundant side into the
 barren side.  Two different ordinary edges contribute
 at least $3$ distinct vertices to the barren side.
 This is a contradiction.

 We have ruled out all possible behavior for $P=P(1)$ assuming
 that $Q=Q(1)$ is degenerate.  Hence, $Q(1)$ is not degenerate.
 This means that $Q(1)$ is a bird.  This completes the
 proof that
 \begin{equation}
   \label{halfresult}
   \Delta_k(B_{k,n}) \subset B_{k,n}.
 \end{equation}

 \subsection{Equality}
 \label{equality}

 We use the notation from \S \ref{duality}.
 Equation \ref{factor} implies that
\begin{equation}
  \label{factor2}
  \Delta_k^{-1}=D_{k+1} \circ \Delta_k \circ D_{k+1}.
\end{equation}
So far, Equation \ref{factor2} makes sense in terms of
PolyPoints and PolyLines.

Below we will explain how
to interpret $D_{k+1}$ as a map from polygons
in $\P$ to polygons in $\P^*$ and also as a map
from polygons in $\P^*$ to polygons in $\P$.
 Since the dual projective plane
 $\P^*$ is an isomorphic copy of $\P$, it makes sense to define
$B_{k.n}^*$.  This space is just the image of $B_{k,n}$ under any projective
duality.
Below we will prove
\begin{theorem}
  \label{dual}
  $D_{k+1}(B_{k,n}) \subset B^*_{k,n}$.
\end{theorem}
It then follows from projective duality that 
$D_{k+1}(B^*_{k,n}) \subset B_{k,n}$.
Combining these equations with Equation \ref{factor2} we see that
$\Delta_k^{-1}(B_{n,k}) \subset B_{n,k}$.
This combines with Equation \ref{halfresult} to finish the proof of
Theorem \ref{main}.

Now we prove Theorem \ref{dual}.

\begin{lemma}
  \label{enhance}
  If $P \in B_{k,n}$, then we can enhance $D_{k+1}(P)$ in such
  a way that $D_{k+1}(P)$ is a planar polygon in $\P^*$.
  The enhancement varies continuously.
\end{lemma}

\startproof
A {\it polygon\/} is a PolyPoint together with additional
data specifying an edge in $\P$ joining each consecutive
pair of points.  Dually, we get a
polygon in $\P^*$  from a PolyLine by specifying, for each
pair of consecutive lines $L_j,L_{j+1}$, an arc of the pencil of lines
through the intersection point which connects $L_j$ to $L_{j+1}$.

Specifying an enhancement of $D_{k+1}(P)$ is the same as
specifing, for each consecutive pair $L_1,L_2$ of
$(k+1)$ diagonals of $P$, an arc of the pencil through their
intersection that connects $L_1,L_2$.  There are two
possible arcs. One of them avoids the interior of the soul of $P$ and the
other one sweeps through the soul of $P$.  We choose
the arc that avoids the soul interior.  Figure 6.4 shows that we
mean for a concrete example.

\begin{center}
\resizebox{!}{1.8in}{\includegraphics{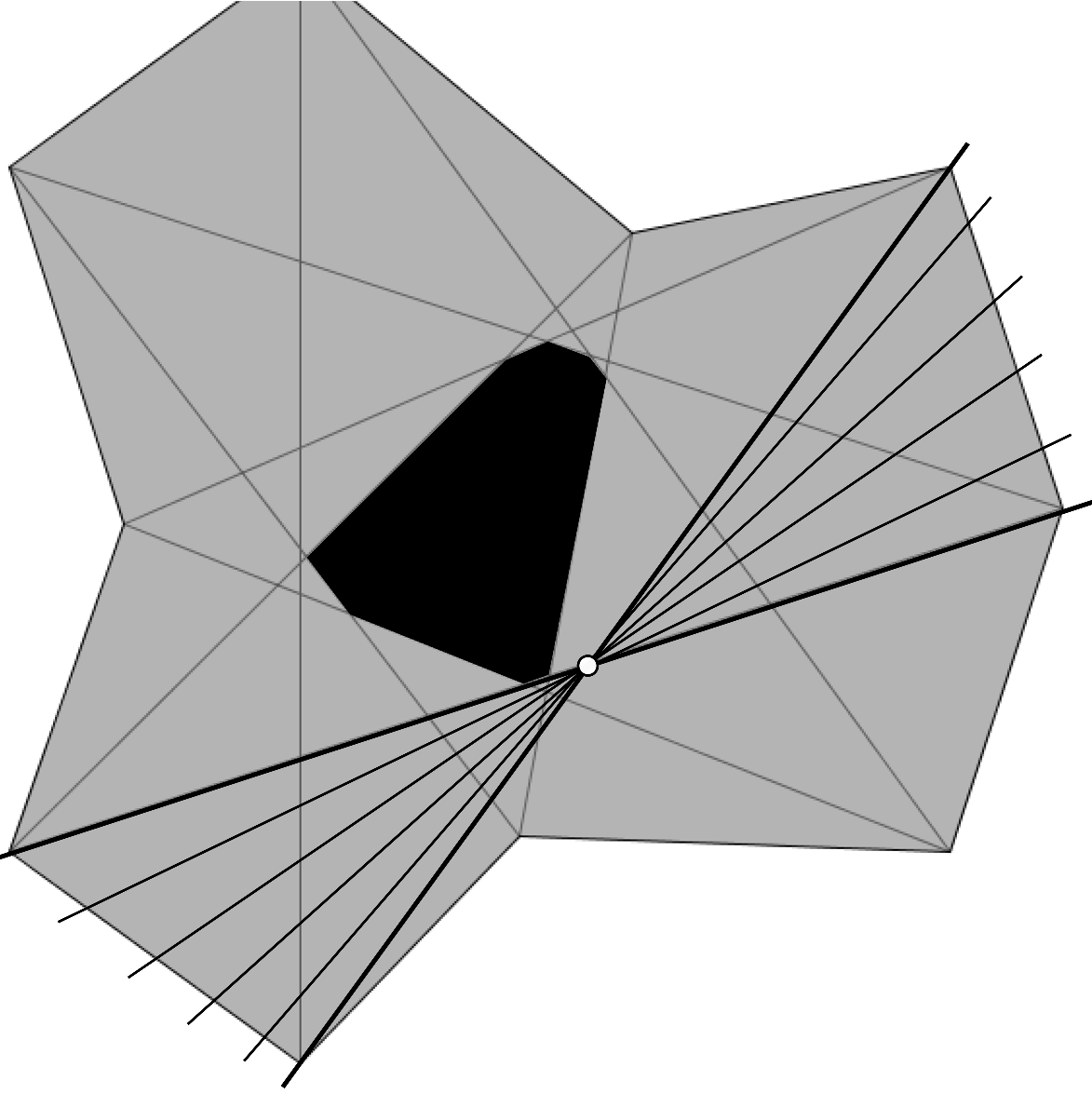}}
\newline
Figure 6.4:  Enhancing a PolyLine to a polygon: Avoid the soul.
\end{center}

Since the soul of $P$ has non-empty interior, there
exists a point $x \in P$ which is disjoint from all
these pencil-arcs.  Applying duality, this exactly says
that there is some line in $\P^*$ which is disjoint
from all the edges of our enhanced $D_{k+1}(P)$.
Hence, this enhancement makes $D_{k+1}(P)$ planar.
Our choice also varies continuously on
$B_{n,k}$.
\endproof

\begin{lemma}
$D_{k+1}$ maps a member of $B_{k,n}$ to an
  $n$-gon which is $k$-nice.
\end{lemma}

\startproof
Let $Q=D_{k+1}(P)$.  A $(k+1)$-diagonal of $Q$ is just a vertex of $P$.
A $k$ diagonal of $Q$ is a vertex of $\Delta_k(p)$.   Thus, to check the
$k$-nice property for $Q$ we need to take $n$-collections of
$4$-tuples of points and check that they are distinct. In each case,
the points are collinear because the lines of $Q$ are coincident.

\begin{center}
\resizebox{!}{1.21in}{\includegraphics{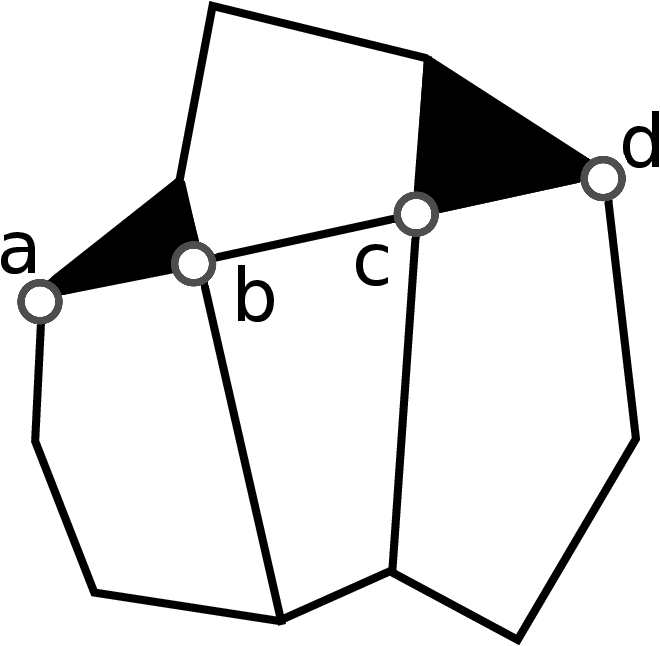}}
\newline
Figure 6.5 One of the $n$ different $4$-tuples we need to check.
\end{center}

Once we make this specification, there is really combinatorially
only possibility for which collections we need to check.
Figure 6.5 shows one such $4$-tuple, $a,b,c,d$.
The shaded triangles are the two feathers of $P$ whose tips
are $b,c$. But $a,b,c,d$ are distinct vertices of $P \cup \Delta_k(P)$ and
so they are distinct.  That is all there is to it.
\endproof

To show that $Q=D_{k+1}(P)$ is a $k$-bird, we consider a
continuous path $P(t)$ from the regular $n$-gon $P(0)$ to
$P=P(1)$. We set $Q(t)=P(t)$.  By construction,
$Q(0)$ is a copy of the regular $n$-gon in $\P^*$, and
$Q(t)$ is $k$-nice for all $t \in [0,1]$, and
$Q(t)$ is a planar polygon for all $t \in [0,1]$.
By definition $Q=Q(1)$ is a $k$-bird.
This completes the proof of
Theorem \ref{dual}.

\newpage

\section{The Triangulation}

\subsection{Basic Definition}
\label{triang}

In this section we gather together the results we have
proved so far and explain how we construct
the triangulation $\tau_P$ associated to a bird
$P \in B_{k,n}$.

Since $\Delta_k(B_{k,n}) \subset B_{k,n}$, we know that
$\Delta_k(P)$ is also a $k$-bird.
Combining this with Theorem \ref{soul}
and Theorem \ref{feather} we can say
that $\Delta_k(P)$ is one embedded $n$-gon
contained in $P^I$, the interior of the region bounded
by the embedded $P$.  The region between
$P$ and $\Delta_k(P)$ is a topological annulus.
Moreover,
$\Delta_k(P)$ is obtained from $P$ by connecting
the tips of the feathers of $P$. The left side
Figure 7.1 shows how this region is triangulated.  The black
triangles are the feathers of $P$ and
each of the white triangles is made from an
edge of $\Delta_k(P)$ and two edges of
adjacent feathers.

\begin{center}
\resizebox{!}{2in}{\includegraphics{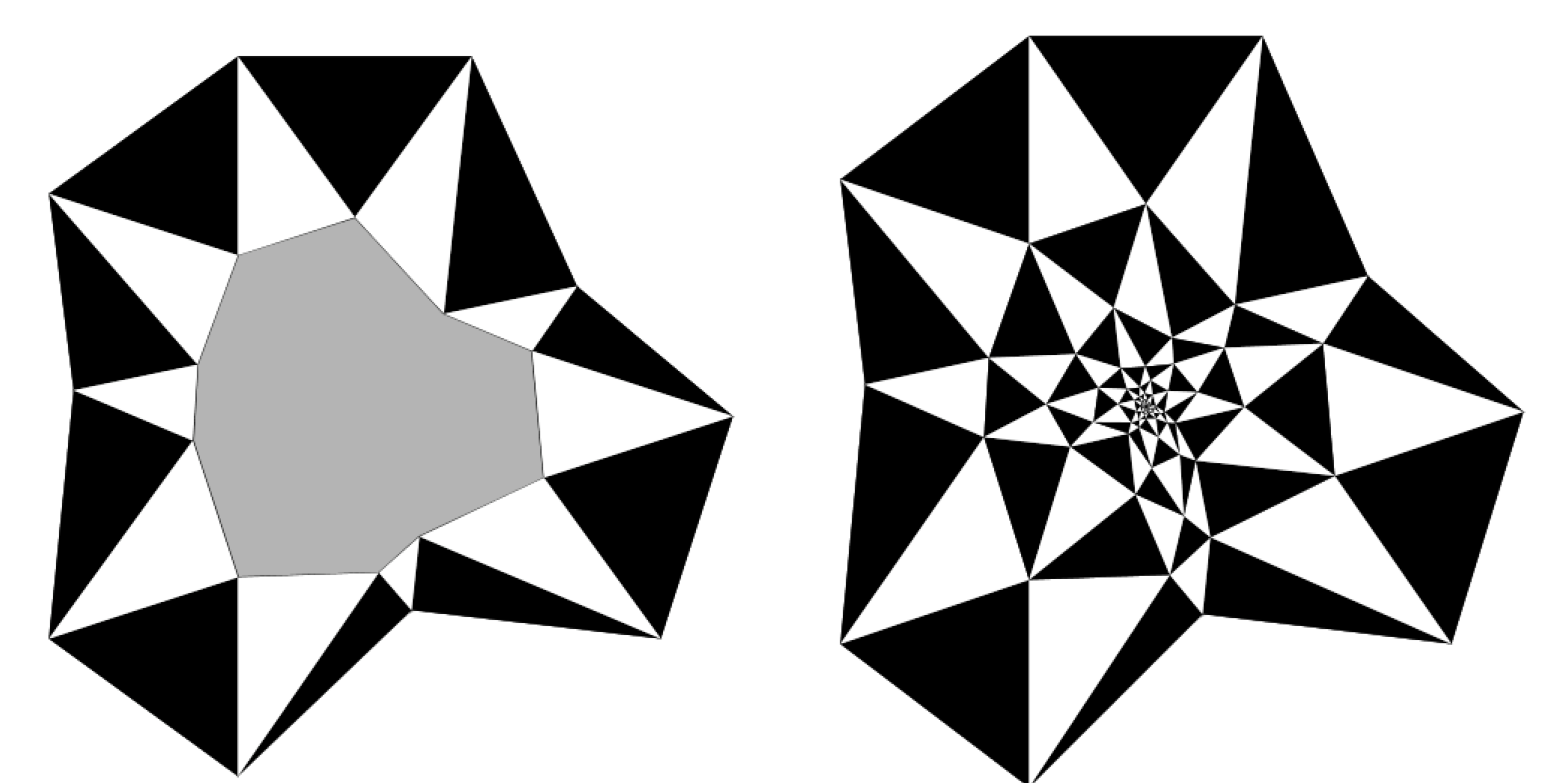}}
\newline
Figure 7.1: The triangulation of the annulus
\end{center}

\begin{lemma}
  For every member $P \in B_{k,n}$, the associated
  $2n$ triangles have pairwise
disjoint interiors, and thus triangulate the annular
region between $P$ and $\Delta_k(P)$.
\end{lemma}

\startproof
As usual, we make a homotopical argument. If this result is
false for some $P$, then we can look at path which starts
at the regular $n$-gon (for which it is true) and stop at
the first place where it fails.  Theorem \ref{feather} tells
us that nothing goes wrong with the feathers of $P$.
The only thing that can go wrong is $\Delta_k(P)$ fails
to be an embedded polygon.  Since this does not happen,
we see that in fact there is no counter-example at all.
\endproof

We can now iterate, and produce $2n$ triangles
between $\Delta_k(P)$ and $\Delta_k^2(P)$, etc.
The right side of Figure 7.1 shows the result of
doing this many times.
The fact that $\Delta_k(B_{k,n})=B_{k,n}$ allows us
to extend outward as well.  When we iterate forever
in both directions, we get an infinite triangulation of
a (topological) cylinder that has degree $6$ everywhere.
This is what Figure 1.6 is showing.  We call
this bi-infinite triangulation $\tau_P$.

\subsection{Some Structural Results}
\label{struct0}

The following result will help with the proof of
Theorem \ref{aux}.

\begin{theorem}
  \label{soul2}
  Let $P \in B_{n,k}$. Let $S$ be the
  soul of $B$.  Then for $\ell \geq n$ we have
  $\Delta^{\ell}_k(P) \subset S$.
\end{theorem}

\startproof
We first note the existence of
certain infinite polygonal arcs in $\tau_P$.
We start at a vertex of $P$ and then move inward
to a vertex of $\Delta_k(P)$ along one of the
edges.  We then continue through this vertex
so that $3$ triangles are on our left and $3$
on our right.  Figure 7.2 below shows the two paths
like this that emanate from the same vertex of $P$.

\begin{center}
\resizebox{!}{3.6in}{\includegraphics{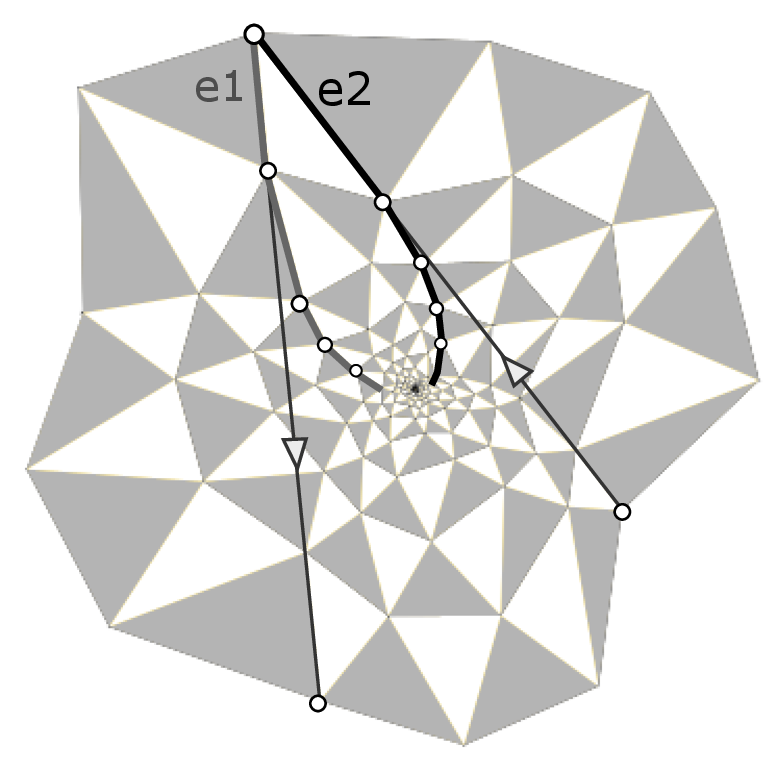}}
\newline
Figure 7.2: The spiral paths.
\end{center}

The usual homotopical argument establishes the fact
that the spiral paths are locally convex.  One can
understand their combinatrics, and how they relate to the
polygons in the orbit, just by looking at the case of
the regular $n$-gon.  We call the two spiral paths in
Figure 7.2 {\it partners\/}.  In the regular $n$-gon the
partners intersect infinitely often. So this is true in general.
Each spiral path has an initial segment joining the initial
endpoint on $P$ to the first intersection point with the
partner.  We define a {\it petal\/} to be the region bounded
by the initial paths of the two partners.

It is convenient to write $P^{\ell}=\Delta_k^{\ell}(P)$.
In the regular case, $P^{\ell}$ is contained in
the petal for $\ell>n-1$..  Hence, the same goes in the
general case.  Because the initial segments are locally
convex, the petal lies to the left of the lines extending
the edges $e_1$ and $e_2$ when these edges are
oriented according to the $(k+1)$-diagonals of $P$.
But this argument works for every pair of
partner spiral paths which start at a vertex of $P$.  We
conclude that for $\ell \geq n$, the polygon
$P^{\ell}$ lies to the left of all the
$(k+1)$-diagonals of $P$.
But the soul of $P$ is exactly the intersection of all these left half planes.
\endproof

Theorem \ref{soul2} in turn gives us
information about the nesting properties of
birds within an orbit.
Let $S_{\ell}$ denote the soul of $P^{\ell}$.   Let
\begin{equation}
  S_{\infty}=\bigcap_{\ell \in \Z}  S_{\ell}, \hskip 30 pt
  S_{-\infty}=\bigcup_{\ell \in \Z} S_{\ell}.
\end{equation}
It follows from Theorem \ref{soul2} that
$\widehat P_{\infty}=S_{\infty}$ and
$\widehat P_{-\infty}=S_{-\infty}$, because
\begin{equation}
S_{\ell+n} \subset  P^{\ell+n} \subset S_{\ell} \subset P^{\ell}.
\end{equation}
Hence these sets are all convex
subsets of an affine plane.

\begin{corollary}
  \label{star-extra}
  Any $P \in B_{k,n}$ is strictly star-shaped with
  respect to all points in the convex hull of $\Delta_k^n(P)$.
\end{corollary}

\startproof
    Since $P^{\ell+n} \subset S_{\ell}$, and
    $P^{\ell}$ is strictly star shaped with
    respect to all points of $S^{\ell}$, we
    see that $P^{\ell}$ is strictly star shaped
    with respect to all points of
    $P^{\ell +n\/}$.  Since $S_{\ell}$ is convex, we can
    say more strongly that $P^{\ell}$ is strictly
    star-shaped with respect to all points of
    the convex hull of $P^{\ell+n}$.
    Now we just set $\ell=0$ and recall the meaning
    of our notation, we get the
    exact statement of the result.
    \endproof

    An immediate corollary is that
    $P$ is strictly star-shaped with
    respect to $\widehat P_{\infty}$.
    (Theorem \ref{aux} says that this is a
    single point.)

\newpage

\section{Nesting Properties of Birds}

\subsection{Duality}

In this chapter we prove Theorem \ref{aux}.
In this first section we show how Statement 1 of
Theorem \ref{aux} implies Statement 2.    We want to prove
that the ``backwards union''
$\widehat P_{-\infty}$ is an affine plane.
Here $P \in B_{n,k}$ is a $k$-bird.
  
    We take $\ell \geq 0$ and consider
    $P^{-\ell}=\Delta_k^{-\ell}(P)$.
    Since $P^{-\ell}$ is planar, there is a closed set
    $\Lambda_{\ell}$
    of lines in $\P$ which miss $P^{-\ell}$.   These sets of
    lines are nested: $\Lambda_1 \supset \Lambda_2 \supset \Lambda_3
    ...$.
    The intersection is non-empty and contains some line $L$.
    We can normalize so that $L$ is the line at infinity.  Thus all
    $P^{-\ell}$ lie in $\R^2$.   We want to see that
    $\widehat P_{-\infty}=\R^2$.

Let 
$D_{k+1}$ be the map from
\S \ref{duality} and \S \ref{equality}.
   From Equation \ref{factor} we see that
   $D_{k+1}$ conjugates $\Delta_k$ to
   $\Delta_k^{-1}$.    With Theorem \ref{dual} in mind,
   define the following ``dual'' $k$-birds:
  \begin{equation}
    \Pi^{\ell}=\Delta_k^{\ell}(D_{k+1}(P))= D_{k+1}(P^{-\ell}).
  \end{equation}
  From Statement 1 of Theorem \ref{aux}, the sequence of $k$-birds
 $\{\Pi^{\ell}\}$
    shrinks to a point
    in the dual plane $\P^*$.
    The vertices of $\Pi^{\ell}$ are the $(k+1)$-diagonals of
    $P^{-\ell}$.
    Because the vertices of $\Pi^{\ell}$ shrink to
    a single point, all the $(k+1)$-diagonals
    of $P^{-\ell}$ converge to a single line $L'$.

    \begin{lemma}
      $L'$ is the line at infinity.
    \end{lemma}

    \startproof
    Suppose not.  When $\ell$ is large, all the $(k+1)$-diagonals
    point nearly in the same direction as $L'$.  In particular,
    this is true of the subset of these diagonals which
    define the soul $S^{-\ell}$.  But these special diagonals turn
    monotonically and by less than $\pi$ radians as
    we move from one to the next.
    Hence, 
    some of these diagonals nearly point in one direction
    along $L'$ and some point nearly in the opposite direction.
    But then $S^{-\ell}$ converges to a subset of $L'$.  This is a contradiction,
    \endproof

   The soul $S^{-\ell}$
   is a convex set, containing the origin, and is
    bounded by some of the $(k+1)$ diagonals.  If
    $S^{-\ell}$ does not converge to the whole plane,
    then some $(k+1)$-diagonal intersects a uniformly
    bounded region in $\R^2$ for each $\ell$.  But this
    produces a sequence of $(k+1)$-diagonals that
    does not converge to the line at infinity.
    This is a contradiction.
    Hence $S^{-\ell}$ converges to all of $\R^2$.
    But then so does $P^{-\ell}$.

\subsection{The Pre-Compact Case}
\label{compact}

The rest of the chapter is devoted to proving
Statrement 1 of Theorem \ref{aux}.
Let $P \in B_{n,k}$ and let $P^{\ell}=\Delta^{\ell}(P)$.
We take $\ell=0,1,2,3...$.

\begin{conjecture}
  \label{precompact}
  The sequence $\{P^{\ell}\}$ is pre-compact
  modulo affine transformations.   That is,
  this sequence has a convergent subsequence
  which converges to another element of $B_{n,k}$.
\end{conjecture}

In this section I will prove the $\widehat P_{\infty}$ is a single
point
under the assumption that $\{P^{\ell}\}$ is pre-compact.

We would like to see that the diameter of
$P^{\ell}$ steadily shrinks, but the notion of
diameter is not affinely natural.  We first develop
a notion of affinely natural diameter.
For each direction $v$ in the plane, we let
$\|S\|_v$ denote the maximum length of $L \cap S$ where
$L$ is a straight line parallel to $v$.  We then define
\begin{equation}
  \delta(S_1,S_2)=\sup_v \frac{\|S_1\|_v}{\|S_2\|_v} \in [0,1].
\end{equation}
The quantity $\delta(S_1,S_2)$ is affine invariant, and
(choosing a direction $\mu$ which realizes the diamater of $S_1$)
we have
\begin{equation}
  \frac{{\rm diam\/}(S_1)}{{\rm diam\/}(S_2)} \leq \frac{\|S_1\|_{\mu}}{\|S_2\|_{\mu}}\leq \delta(S_1,S_2).
\end{equation}

Let $S^{\ell}$ be the soul of $P^{\ell}$.  By Theorem \ref{soul1} we
have
$S^{\ell+n} \subset S^{\ell}$.    More precisely, the former set is
contained
in the interior of the latter set.
Under the pre-compactness assumption, there are infinitely many
indices $\ell_j$ and some $\epsilon>0$ such that
\begin{equation}
  \delta(S^{\ell_j+n},S^{\ell_j})<1-\epsilon.
\end{equation}
But then
\begin{equation}
  \frac{
    {\rm diam\/}(S^{\ell_j+n})
  }
  {
    {\rm diam\/}(S^{\ell_j})
  }
  <1-\epsilon
\end{equation}
infinitely often.  This forces ${\rm diam\/}(S^{\ell}) \to 0$.  But $\widehat
P_{\infty}$ is contained in this nested intersection and hence is a
point.

If we knew the truth of Conjecture \ref{precompact} then our
proof of Theorem \ref{aux} would be done.  Since we don't
know this, we have to work much harder to prove Statement 1
in general.

\subsection{Normalizing by Affine Transformations}

Henceforth we assume that the forward orbit
$\{P^{\ell}\}$ of $P$ under $\Delta_k$ is not pre-compact
modulo affine transformations.

\begin{lemma}
  \label{NNN}
  There is a sequence $\{T_{\ell}\}$ of affine transformations such that
  \begin{enumerate}
  \item $T_{\ell}(P^{\ell})$ has (the same) $3$ vertices which make a fixed
    equilateral triangle.
  \item $T_{\ell}$ expands distances on $P^{\ell}$ for all $\ell$.
  \item $T_{\ell}(P^{\ell})$ is contained in a uniformly bounded subset of $\R^2$.
    \end{enumerate}
\end{lemma}

\startproof
To $P^{\ell}$ we associate the triangle
$\tau_{\ell}$ made from $3$ vertices of
$P^{\ell}$ and having maximal area.
The diameter of $\tau_{\ell}$ is uniformly
small, so we can find a single equilateral
triangle $T$ and an expanding affine map
$T_{\ell}: \tau_{\ell} \to T$.  Let $d$ be the
side length of $T$. Every vertex
of $T_{\ell}(P^{\ell})$ is within $d$ of all
the sides of $T$, because otherwise we'd
have a triangle of larger area.
The sequence $\{T_{\ell}\}$ has the advertised properties.
\endproof

  Let $Q^{\ell}=T_{\ell}(P^{\ell})$.   By compactness  we can
  pass to a subsequence so that the limit polygon
  $Q$ exists, in the sense that the vertices and the
  edges converge.
  Let $Q_0, Q_1$, etc.
  be the vertices of $Q$.  Perhaps some of these
  coincide.
  Each distinguished diagonal of
$Q^{\ell}$ defines the unit vector which is parallel
to it.  Thus $Q^{\ell}$ defines a certain list of
$2n$ unit vectors.  We can pass to a subsequence
so that all these unit vectors converge.  Thus
$Q$ still has well defined distinguished diagonals
even when the relevant points coincide.

We now define the ``limiting soul''.
  Let $S^{\ell}=S(Q^\ell)$, the soul of $Q^{\ell}$.
  As in
  \S \ref{souldef}. let $S$ be the set of accumulation
  points of sequences $\{p^{\ell}\}$ with $p^{\ell} \in S^{\ell}$.
  Since $S^{\ell} \subset Q^{\ell}$ for all $\ell$ we have
  $S \subset Q$. Now we define a related object.
    We have a left half-plane associated to each diagonal of $Q$.
We define $\Sigma$ to be the intersection of all these
half-planes.   We will use the set $\Sigma$ at various places
below to get control over the set $S$.

\begin{lemma}
  \label{sigma}
  $S \subset \Sigma$.
\end{lemma}

\startproof
Fix $\epsilon>0$.  If this is not the case, then
by compactness
we can find a convergent sequence
$\{p^{\ell}\}$, with $p^{\ell\/} \in S^{\ell}$,
which does not converge to a point of $\Sigma$.
But $p^{\ell}$ lies in every left half plane associated
to $Q^{\ell}$.  But then, by continuity, the accumulation
point $p$ lies in every left half plane associated to $Q$.
Hence $p \in \Sigma$.
\endproof

\subsection{Structure of the Normalized Limits}
\label{limit}

We  work under the assumption that
$\widehat P_{\infty}$ is not a single point.
The goal of this section is to establish several
structural properties about the sets $S$ and $Q$.
Our first property guarantees that there is a
chord $S^*$ of $S$ connecting vertices of
$Q$.  Once we establish this, we show that
$Q$ is a union of two ``monotone'' arcs
joining the endpoints of $S^*$.
These structural properties will be used repeatedly
in subsequent sections of this chapter.

Let $H_Q$ denote the convex hull of $Q$.  Note that
$S \subset Q \subset H_Q$.

\begin{corollary}
  \label{SHRINK}
  Suppose that $\widehat P_{\infty}$ is not a single point.
  Then $\delta(S,H_Q)=1$.
\end{corollary}

\startproof
Suppose not. Note that $H_{Q^{\ell}} \subset S^{\ell-n}$ by Theorem
\ref{soul2} and convexity.  Then for $\ell$ large we have
$$\delta(Q^{\ell-n})=\delta(S^{\ell},S^{\ell-n}) \leq
\delta(S^{\ell},H_{Q^{\ell}}) <
\delta(S,H_Q)+\epsilon,$$
and we can make $\epsilon$ as small as we like.
This gives us a uniform $\delta<1$ such
that $\delta(Q^{\ell})<\delta$ once $\ell$ is large enough.
The argument in the compact case now shows that
$\widehat P_{\infty}$ is a single point.
\endproof

Corollary \ref{SHRINK} says that $S$ and $Q$ have
the same diameter.
Hence there is a chord $S^* \subset S$ which has the same diameter
as $Q$.  Since $Q$ is a polygon, this means that $Q$ must have
vertices at either endpoint of $S^*$.  We normalize so that
$S^*$ is the unit segment joining $(0,0)$ to $(1,0)$.

\begin{lemma}
  \label{structX}
  Let $Q' \subset Q$ be an arc of $Q$ that joins $(0,0)$ to $(1,0)$.
  \begin{enumerate}
\item The vertices of $Q'$ must have non-decreasing $x$-coordinates.
\item If consecutive vertices of $Q'$ have the
  same $x$-coordinate, they coincide.
\item Either $Q' \subset S^*$ or $Q'$ intersects $S^*$ only at $(0,0)$ and $(1,0)$.
  \end{enumerate}
\end{lemma}

\startproof
Suppose the Statement 1 is false. Then
we can find a vertical line $\Lambda$ which intersects
$S^*$ at a relative interior point and which intersects
$Q'$ transversely at $3$ points. But then once $\ell$ is sufficiently
large, $Q^{\ell}$ will intersect all vertical lines sufficiently
close to $\Lambda$ in at least $3$ points and moreover some of
these lines will contain points of $S^{\ell}$.  This contradicts
the fact that $Q^{\ell}$ is strictly star-shaped with respect
to all points of $Q^{\ell}$.

For Statement 2, we observe that
$Q'$ does not contain any point of the form
$(0,y)$ or $(1,y)$ for $y \not =0$. Otherwise
$Q$ has larger diameter than $1$.  This is to say
that once $Q'$ leaves $(0,0)$ it immediately moves
forward in the $X$-direction.  Likewise, once
$Q'$ (traced out the other way) leaves $(1,0)$
it immediately moves backward in the $X$-direction.
If Statement 2 is false, ten we can find a non-horizontal
line $\Lambda'$ which intersects $S^*$ in a relative interior point and which
intersects $Q'$ transversely at $3$ points.  The slope is
$\Lambda'$ depends on which of the two vertices of $Q'$ lies above
the other.  Once we have $\Lambda'$ we play the same game as
for the first statement, and get the same kind of contradiction.

Suppose Statement 3 is false.  We use the kind of argument
we had in \S \ref{ungood}.
By Statements 1 and 2 together,
$Q'$ must have an escape edge which touches $S^*$ in a relative
interior point.  Moreover, this one escape edge is paired with another
escape edge.  Thus we can find a point $x \in S^*$ which strictly
lies on the same side of both of these same-type escape edges.
The argument in \S \ref{ungood} now shows that
$Q^{\ell}$ is not strictly star-shaped with respect to points
of $S^{\ell}$ very near $x$.
\endproof

\begin{corollary}
  \label{pinch}
  Suppose $0 \leq a<b<n$ and $Q_a=Q_b$.  Then
  either we have $Q_a=Q_{a+1}=...=Q_b$ or else we have
  $Q_b=Q_{b+1}=...=Q_{a+n}$.
\end{corollary}

\startproof
In view of Lemma \ref{structX} it suffices to show that
our two monotone arcs comprising $Q$ are
disjoint except at their endpoints.

Let $U$ denote the open upper halfplane, bounded by
the $X$-axis. After reflecting
in the $X$-axis we can guarantee that one of our monotone
arcs $\alpha$ has a point in $U$.  But then, by Lemma \ref{structX},
all of $\alpha$ lies in $U$ except for its endpoints.
If the other monotone arc $\beta$ intersects $\alpha$
away from the endpoints, then $\beta$ has a point in $U$,
but then, by Lemma \ref{structX}, all of $\beta$ lies in
$U$ except for the endpoints.  But then $S$ lies in
$U$, except for the points $(0,0)$ and $(1,0)$.
This contradicts the fact that $S^* \subset S$.
\endproof

Our argument shows in particular that $Q$ is embedded,
up to adding repeated vertices.  However, we will not directly
use this property in our proof below.
  
\subsection{The Triangular Case}
\label{TLT}

We continue with the assumption that
$\widehat P$ is not a single point.
Here we pick off a special case:

\begin{itemize}
  \item There is a line $L$ such that
    $Q_0 \not \in L$.
   \item $Q_k,Q_{k+1},...,Q_{n-k-1},Q_{n-k} \in L$ and
   \item $Q_k \not = Q_{n-k}$.
   \end{itemize}
   Figure 8.1 shows the situation.  As always, the notation
   $Q_{-k}$ and $Q_{n-k}$ names the same point.
   All but $2k-1$ points are on $L$, and except for
   $Q_0$ we don't know where these other $2k-1$ points are.
   
\begin{center}
\resizebox{!}{2.6in}{\includegraphics{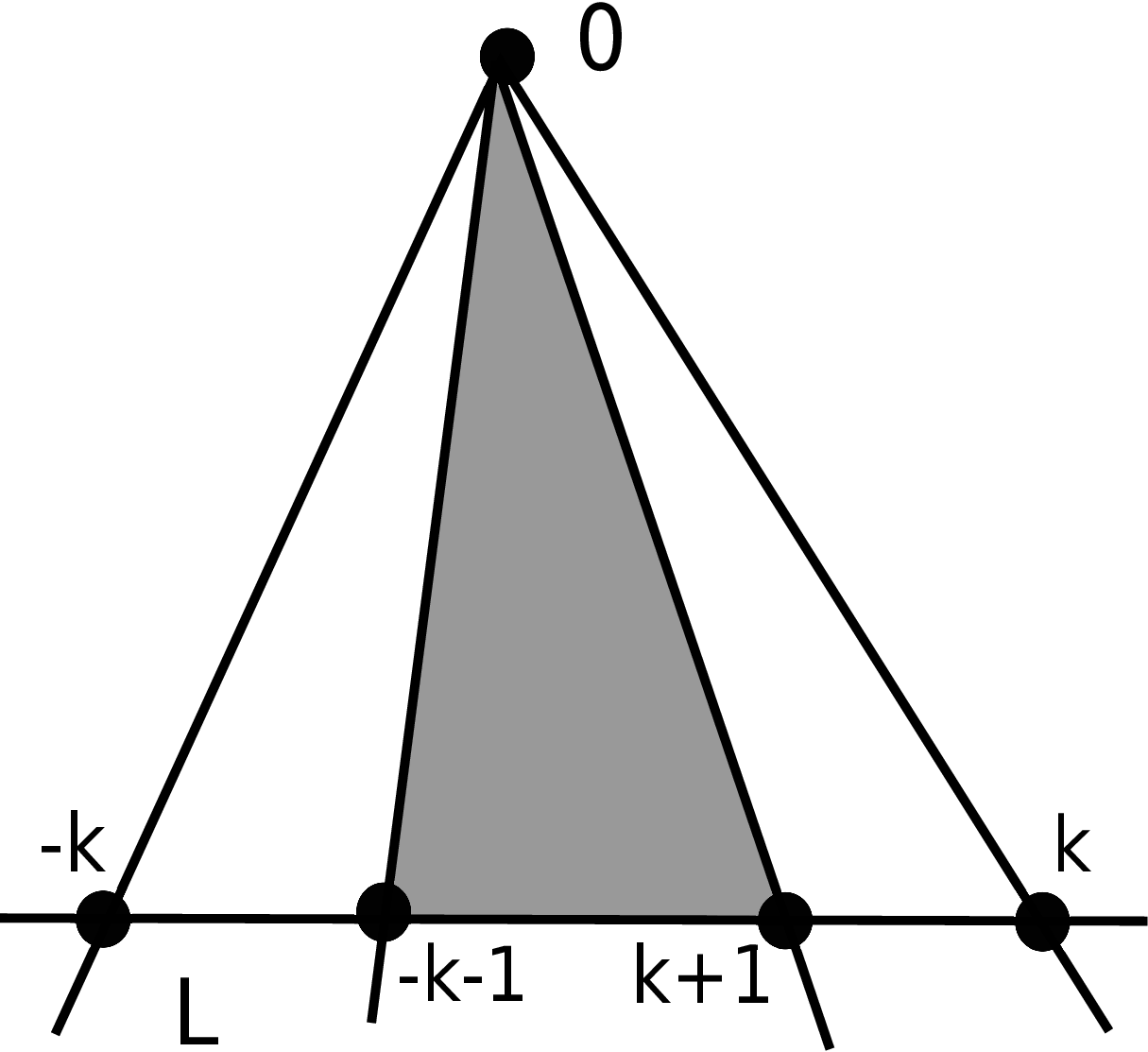}}
\newline
Figure 8.1: The triangular limit $Q$.
\end{center}

Given the constant energy of our orbit, the cross ratio of the lines
$$Q_{0,k},\  Q_{0,k+1},\  Q_{n-k-1,0},\ Q_{n-k,0}$$
is at least $\epsilon_0$. Also, these lines are cyclically
ordered about $0$ as indicated in Figure 8.1, thanks to the
$k$-niceness property and continuity. Also, the two
lines containing $Q_{0,k}$ and $Q_{-k,0}$ are not parallel because
$Q_0 \not \in L$.
Hence $S$ is contained in the shaded region
in Figure 8.1, namely the triangle with
vertices $Q_0$ and $Q_{\pm (k+1)}$.  But this shaded region has diameter strictly
smaller than the triangle $\tau$ with vertices
$Q_0$ and $Q_{\pm k}$.  Hence
${\rm diam\/}(S)<{\rm diam\/}(\tau) \leq {\rm diam\/}(Q)$.
This contradicts Corollary \ref{SHRINK} which says, in
particular, that $S$ and $Q$ have the same diameter.

    \subsection{The Case of No Folded Diagonals}
    \label{index}

    We work under the assumption that
    $\widehat P_{\infty}$ is not a single point.
The notions of collapsed diagonals, folded diagonals, and
aligned diagonals from \S \ref{DEGEN} make sense for $Q$ because
the concepts just involve the directions of the diagonals.
The proof of Lemma \ref{conseq} also works the same way.

\begin{lemma}
  \label{somecollapse}
  $Q$ must either have a trivial edge, a trivial
  distinguished diagonal, or collapsed diagonals, 
\end{lemma}

\startproof
As remarked in \S 5, the proof of the
Degeneration Lemma works for sequences as well as paths, and
only uses the fact that the limiting polygon has nontrivial
edges and nontrivial distinguished diagonals. So,
if $Q$ has no trivial edges and no trivial distinguished diagonals,
then all but one vertex of $Q$
lies in a single line.  But then $Q$ has collapsed diagonals.
\endproof

\noindent
{\bf Remark:\/}
Here is a second, more direct proof. If Lemma \ref{somecollapse} is
false then we have a picture as in the
    left side of Figure 7.1.
    The feathers defined in \S \ref{feather} would be all
    non-degenerate and the segments joining the tips of
    consecutive feathers would be nontrivial.  This would
    force $S$ to lie in the interior of $Q$.  But then
    ${\rm diam\/}(S)<{\rm diam\/}(Q)$, contradicting
    Corollary \ref{SHRINK}.
    \newline
        
   If $Q$ has a trivial
    distinguished diagonal, then by Lemma \ref{pinch}, we see
    that $Q$ also has a trivial edge.
    If $Q$ has a trivial edge, say $Q_{-1}=Q_0$, then the diagonals
    at $Q$ are collapsed at $Q_k$.   So, in all cases, $Q$ has
    collapsed diagonals.
    We assume in this section that $Q$ has no folded diagonals
    anywhere. This means that $Q$ has aligned diagonals, say
    at $Q_k$.   Thus $Q_{0,k}$ and $Q_{k,2k}$ are parallel.
      Since $Q$ does not lie in a line,
      Lemma \ref{conseq} tells us that
      the chain of $2k+1$ parallel distinguished diagonals:
    \begin{equation}
      \label{longchain}
    Q_{0,k}, Q_{0,k+1},
    Q_{1,k+1},Q_{1,k+2},...,Q_{k-1,2k},Q_{k,2k}
    \end{equation}
    Now we have a ``runaway situation''.  The two diagonals
    $Q_{2k,k}$ and $Q_{2k,k-1}$ (which are just the reversals of
    the last two in Equation \ref{longchain})  are parallel.
        Thus $Q$ has collapsed diagonals at $Q_{2k}$.   Since
    $Q$ has no folded diagonals, $Q$ has aligned diagonals at
    $Q_{2k}$.  But then, applying Lemma \ref{degen} again,
    we can extend that chain in Equation \ref{longchain} so that
    it contines as  $,...,Q_{2k-1,3k},Q_{2k,3k}$.
    But now $Q$ has collapsed diagonals at $Q_{3k}$.  And so on.
    Continuing this way, we end up with all points on $Q$.
    This is a contradiction.

    The only way out is that $Q$ must have folded diagonals somewhere

\subsection{The Case of Folded Diagonals}
\label{foldedcase}

We continue to work under the assumption that
$\widehat P_{\infty}$ is not a single point.  Now we consider
the case when $Q$ has folded diagonals at, say, $Q_0$.
What this means that the diagonals
$Q_{0,k+1}$,  $Q_{0,-k-1}$ are parallel.  (Again, these diagonals are
well defined even when their endpoints coincide; we are just
using a notational convention to name them here.)
But then the
corresponding half planes intersect along a single line $L$, forcing
$\Sigma \subset L$.  By Lemma \ref{sigma}, the soul $S$ is
contained in $\Sigma$.  Hence, $S \subset L$.
Letting $S^*$ be the chord from \S \ref{limit}, we also have
$S=S^*$.   This is because $S$ and $S^*$ are segments of
the same diagonal and in the same line.  We will use
$S$ and $S^*$ interchangeably below.

We normalize so that $S$ is the line segment connecting $(0,0)$ to
$(1,0)$.   As in \S \ref{limit}, both these points are vertices of
$Q$.   The folding condition forces $\Sigma$ (and hence $S$) to
lie to one side of these points.   Hence, we have either
$Q_0=(0,0)$ or $Q_0=(1,0)$.  Without loss of generality we consider
the case when $Q_0=(0,0)$.
Note that points of $Q-S$ do not belong to $L$, because
$Q$ and $S$ have the same diameter.
We break the analysis down into cases.
\newline
\newline
{\bf Case 1:\/}
Suppose that $Q_{k+1}$ is not an
endpoint of $S^*$ and $Q_{n-k-1} \not = (0,0)$.
Consider the arc $Q'$ given by
$Q_0 \to ... \to Q_{k+1} \to ... \to Q_{\beta}=(1,0)$.
Here $\beta$ is some index we do not know explicitly,
but we take $\beta$ as large as possible, in the sense
that $Q_{\beta+1} \not = (1,0)$.
The arc $Q'$ connects $(0,0)$ to $(1,0)$ and
intersets $S^*$ at
$Q_{k+1}$, a point which is neither
$(0,0)$ or $(1,0)$.
By Lemma \ref{structX}, we have
$Q' \subset S^*$.
We conclude that
$Q_0,...,Q_{\beta} \subset S^*$.

If $\beta$ does not lie in the index
interval $(k+1,n-k-1)$ then we have
just shown that
$Q_{k+1},...,Q_{n-k-1} \in S^*$.
If $\beta=n-k-1$ we have the same result.
Here is what we do if $\beta$ does lie in
$(k+1,n-k-1)$.     We apply our
same argument as in the previous paragraph
to the arc
$Q_{\beta} \to ... \to Q_{n-k-1}$, and
see that $Q_{\beta},...,Q_{n-k-1} \in S$.
So, in all cases, we see that
$Q_{k+1},...,Q_{n-k-1} \in S$.

In short, $Q_j \in L$ unless
$j \in \{-k,...,-1\}$.  All
but $k$ vertices belong to $L$.
In particular, we have an index
$h \in \{-k,...,-1\}$ such that
$Q_h \not \in L$ but
$Q_{h+k}, Q_{h+k+1},...,Q_{h+n-k-1},Q_{h+n-k} \in L$.
Now we are close to the Triangular case from \S \ref{TLT}
except that all the indices are shifted by $h$.
If it happens that $Q_{h+k} \not = Q_{h+n-k}$ then we
have the Triangular Case and we are done.

The other possibility is that
$Q_{h+k} = Q_{h+n-k}$.
In this case, Lemma \ref{pinch} gives
us $Q_{h+k}= Q_{h+k+1}=Q_{h+n-k-1}=Q_{h+n-k}$.
In particular, the diagonals
$Q_{h,h+k+1}$ and $Q_{h,h+n-k-1}$ are folded
at $Q_h$.  Since $Q_h \not \in L$ this means
that there is some other line $L'$ such that
$S \subset L'$.  This is a contradiction.
\newline
\newline
{\bf Case 2:\/}
Suppose $Q_{-k-1}=Q_{k+1}=(1,0)$.
Before analyzing this case, we remember a
lesson from the end of Case 1:  It is not
possible for $Q$ to have folded diagonals
at a point not on $S$.

Corollary \ref{pinch} says that
$Q_{k+1}=...=Q_{n-k-1}=(1,0).$
This is a run of $k+\beta$ points, where
$\beta=n - (3k+1) \geq 0$.
There is some index
$h \in \{\pm 1,...\pm k\}$ such
that $Q_h \not \in L$.
Without loss of generality we will
take $h \in \{1,...,k\}$.

Suppose first that $n>3k+1$.
Then there are at least $k+1$
consecutive vertices sitting at
$(1,0)$ and so both
diagonals $Q_{h,k+h}$ and
$Q_{h,k+h+1}$ point from
$Q_h$ to $(1,0) \not = Q_h$.
This means that $Q$ has
collapsed diagonals at $Q_h$.
Remembering our lesson, we know that
$Q$ does not have folded diagonals
at $Q_h$.    Hence $Q$ has aligned
diagonals at $Q_h$.

Now we have the same runaway situation we
had in \S \ref{index}.  The diagonals in the chain
$Q_{h-k,h}...Q_{h,h+k}$ point are all pointing along the
line connecting $(1,0)$ to $Q_h$, and they are pointing
away from $(1,0)$.
This gives us collapsed diagonals at $Q_{h+k}$.
Remembering our lesson, we see that $Q$ has
aligned diagonals at $Q_{h+k}$.   And so on.
All the points after $Q_h$ get stuck on $L'$ and
we have a contradiction.

If $n=3k+1$, then the same argument
works as long as $h \not = \pm k$.
So, we just have to worry about the case
when all points of $Q$ belong to $S$
except for $Q_k$ and $Q_{-k}$,
which do not belong to $S$.
Applying Lemma \ref{structX} to
the arc $Q_0 \to Q_1 \to ... \to Q_k \to (1,0)$
we conclude that $Q_0=...=Q_{k-1}=(0,0)$. 
Applying Lemma \ref{structX} to
the arc $Q_0 \to Q_{-1} \to ... \to Q_{-k} \to (1,0)$
we conclude that $Q_0=...=Q_{k-1}=(0,0)$.
But now we have a run of  $2k-1 \geq k+1$
points sitting at $(0,0)$ and we can run the
same argument as in the case $n>3k+1$,
with $(0,0)$ in place of $(1,0)$.
            \newline
            \newline
            {\bf Case 3:\/}  The only cases left to consider
            is when one or both of $Q_{\pm(k+1)}$ equals $(0,0)$.
            We suppose  without loss of generality that $Q_{-k-1}=(0,0)$.
            Since we also
            have $Q_0=(0,0)$, 
            Lemma \ref{pinch} gives
            $Q_{-k-1}=...=Q_0=(0,0)$.
            This is a run of $k+2$ consecutive points sitting
            at $(0,0)$.

            There is some smallest $h>0$ so that $Q_h \not \in S$.
            Applying Lemma \ref{structX} to
            the arc $Q_0 \to ... \to Q_k \to ... \to
            (1,0)$, we conclude that $Q_{h-1}=...=Q_1=(0,0)$.
            (Otherwise Lemma \ref{structX} would force
            $Q_h \in S$.)

            Now we know that $Q$ has collapsed diagonals
            at $Q_h \not \in L$.   We now get a contradiction
            from the same runaway situation as in Case 2.

\newpage

\section{Appendix}
\label{ANTON}

\subsection{The Energy Invariance Revisited}

In this section we sketch Anton Izosimov's proof
that $\chi_k \circ \Delta_k=\chi_k$.  This proof requires
the machinery from \cite{IK}.
(The perspective comes from \cite{IZOS1},
but the needed result for
$\Delta_k$ is in the follow-up paper \cite{IK}.)

Let $P$ be an $n$-gon.  We let
$V_1,...,V_n$ be points in $\R^3$ representing the consecutive
vertices of $P$.  Thus the vertex $P_j$ is the equivalence
class of $V_j$.   We can choose periodic sequences
$\{a_i\}$, $\{b_i\}$, $\{c_i\}$, $\{d_i\}$ such that
\begin{equation}
  \label{shift0}
  a_i V_i + b_i V_{i+k} + c_i V_{i+k+1} + d_i V_{i+2k+1}=0,
  \hskip 30 pt
  \forall i.
\end{equation}
Recall from \S \ref{duality} that
$\Delta_k=D_k \circ D_{k+1}$.

\begin{lemma}
  One of the cross ratio factors of
  $\chi_k \circ D_{k+1}$ is
  $(a_0d_{-k})/(c_0 b_{-k})$.
\end{lemma}

\startproof
One of the factors is the cross ratio of
$P_0, y, x, P_{k+1}$, where
$$x=P_{0,k+1} \cap P_{k,2k+1}, \hskip 30 pt
y=P_{-k,1} \cap P_{0,k+1}.$$
(Compare the right side of Figure 2.1, shifting
all the indices there by $k+1$.)

The points $x$ and $y$ respectively are represented by vectors
$$X=a_0V_0+c_0V_{k+1} = -b_0 V_k - d_0 V_{2k+1},$$
$$Y=-a_{-k}V_{-k} - c_{-k} V_1 = b_{-k} V_0 + d_{-k} V_{k+1}.$$
The point here is that the vector $X$ lies in the span of
$\{V_0,V_{k+1}\}$ and in the span of
$\{V_k,V_{2k+1}\}$ and projectively this is exactly what is
required.  A similar remark applies to $Y$.

Setting $\Omega=V_0 \times V_{k+1}$, we compute the relevant cross ratio as
\begin{equation}
  \frac{V_0 \times Y}{V_0 \times X} \cdot
  \frac{X \times V_{k+1}}{Y \times V_{k+1}}=
  \frac{d_{-k} \Omega}{c_0 \Omega} \times
  \frac{a_0 \Omega}{b_{-k} \Omega} =
  \frac{d_{-k}a_0}{b_{-k}c_0},
\end{equation}
which is just a rearrangement of the claimed term.
\endproof

The other cross ratio factors are obtained by
shifting the indices in an obvious way.
As an immediate corollary, we see that
\begin{equation}
  \label{key1}
  \chi_k(D_{k+1}(P))=
  \prod_{i=1}^n \frac{a_i d_i}{b_ic_i}.
\end{equation}
Let us call this quantity $\mu_k(P)$.

\begin{lemma}
  If  $\mu_k \circ \Delta_k = \mu_k$
  then $\chi_k \circ \Delta_k = \chi_k$.
\end{lemma}

\startproof
If $\mu_k \circ \Delta_k=\mu_k$ then
$\mu_k \circ \Delta_k^{-1}=\mu_k$.
Equation \ref{key1} says that
\begin{equation}
  \chi_k \circ D_{k+1} = \mu_k, \hskip 30 pt
  \mu_k \circ D_{k+1}=\chi_k.
\end{equation}
The first equation implies the second because $D_{k+1}$ is an involution.
Since $D_{k+1}$ conjugates $\Delta_k$ to $\Delta_k^{-1}$ we have
$$\chi_k \circ \Delta_k =
\chi_k \circ D_{k+1} \circ \Delta_k^{-1} \circ D_{k+1}=
\mu_k \circ \Delta_k^{-1} \circ D_{k+1}=\mu_k \circ D_{k+1}=\chi_k.$$
This completes the proof.
\endproof

Let
$\widetilde P=\Delta_k(P)$.  Let
$\{\widetilde a_i\}$, etc., be the sequences
associated to $\widetilde P$. We want to show that
\begin{equation}
  \label{key2}
  \prod_{i=1}^n \frac{a_i d_i}{b_ic_i}=  \prod_{i=1}^n \frac{\widetilde a_i \widetilde d_i}{\widetilde b_i\widetilde c_i}.
\end{equation}
This is just a restatement of the equation $\mu_k \circ \Delta_k=\mu_k$.

Now we use the formalism from \cite{IK} to establish
Equation \ref{key2}. 
We associate to our polygon $P$ operator $D$ on the space $\cal V$ of bi-infinite
sequences $\{V_i\}$ of vectors in $\R^3$.   The definition of $D$ is given
coordinate-wise as
\begin{equation}
  \label{shift1}
  D(V_i)=a_i V_i + b_i T^k(V_i) + c_i T^{k+1}(V_i) + d_i T^{2k+1}(V_i).
\end{equation}
Here $T$ is the shift operator, whose action is $T(V_i)=V_{i+1}$.
If we take $\{V_i\}$ to be a periodic bi-infinite sequence of vectors
corresponding to our polygon $P$, then
$D$ maps $\{V_i\}$ to the $0$-sequence.

Next, we write $D=D_++D_-$ where coordinate-wise
\begin{equation}
  D_+(V_i)=a_i V_i + c_i T^{k+1}(V_i), \hskip 30 pt
  D_-(V_i)=b_i T^k(V_i) + d_i T^{2k+1}(V_i).
\end{equation}
The pair $(D_+,D_-)$ is associated to the polygon $P$.

Let $\widetilde D$ and $(\widetilde D_+,\widetilde D_-)$ be the
corresponding operators associated to $\widetilde P$.
One of the main results of \cite{IK} is that the various choices can be made so that
\begin{equation}
  \label{shift3}
  \widetilde D_+ D_- = \widetilde D_- D_+.
\end{equation}
This is called {\it refactorization\/}.
Equating the lowest (respectively highest) terms of the
relation in Equation \ref{shift3} gives us the identity
$\widetilde a_i b_i = \widetilde b_i a_{i+k}$ (respectively
$\widetilde c_i d_{i+k+1} = \widetilde d_i c_{i+2k+1}$.)
These relations hold for all $i$ and together imply Equation \ref{key2}.

\subsection{Extensions of Glick's Formula}
\label{Glick}

Theorem 1.1 in
\cite{GLICK1} says that the coordinates for the collapse
point of the pentagram map $\Delta_1$  are algebraic functions
of the coordinates of the initial polygon.  In Equation
1.1 of \cite{GLICK1}, Glick goes further and gives a formula
for the collapse point.
I will explain
his formula.
Let $(x^*,y^*)$ denote the accumulation point of the
forward iterates of $P$ under $\Delta_1$.  Let
$\widehat P_{\infty}=(x^*,y^*,1)$ be the collapse point.
In somewhat different notation,
Glick introduces the operator
\begin{equation}
  T_P=nI_3 -G_P,  \hskip 30 pt
  G_{P}(v)=\sum_{i=1}^n \frac{|P_{i-1},v,P_{i+1}|}{|P_{i-1},P_i,P_{i+1}|} P_i.
\end{equation}
Here $|a,b,c|$ denotes the determinant of the matrix with rows $a,b,c$
and $I_3$ is the $3 \times 3$ identity matrix.
It turns out $T_P$ is a $\Delta_1$-invariant operator, in the sense that
$T_{\Delta_0(P)} = T_P$.  Moreover $P_{\infty}$ is an eigenvector of $T_P$.
This is Glick's formula for $\widehat P_{\infty}$.
Actually, one can say more simply that $G_P$ is a $\Delta_0$-invariant
operator and that $\widehat P_{\infty}$ is a fixed point of the
projective action of $G_p$.  This means that the vectors representing
these points in $\R^3$ are eigenvectors for the operator.
The reason Glick uses the more
complicated expression $nI_3-G_P$ is that geometrically
it is easier to work with.

Define $G_{P,a,b}$ by the formula
\begin{equation}
  G_{P,a,b}(v)=\sum_{i=1}^n \frac{|P_{i-a},v,P_{i+b}|}{|P_{i-a},P_i,P_{i+b}|} P_i.
\end{equation}
Let $\widehat P_{\infty,k}$ be the limit point of the
forward iterates of $P$ under $\Delta_k$.

  A lot of experimental evidence
    suggests the following conjecture.
\begin{conjecture}
  Let $k \geq 2$.  If $n=3k+1$ the point $\widehat P_{\infty}$
  is a fixed point for the projective action of
  $G_{P,k,k}$.
  If $n=3k+2$ the point $\widehat P_{\infty}$
is a fixed point for the projective action of
$G_{P,k+1,k+1}$.  In particular, in these cases
the coordinates
of $\widehat P_{\infty}$ are algebraic functions of
the vertices of $P$.
\end{conjecture}

Anton Izosimov kindly explained the following lemma,
which seems like a big step in proving the conjecture.
(I am still missing the geometric side of Glick's argument
in this new setting.)

\begin{lemma}
  \label{invariance}
  When $n=3k+1$ the operator $G_{P,k,k}$ is invariant under
  $\Delta_k$.  When $n=3k+2$ the operator $G_{P,k+1,k+1}$ is
  invariant under $\Delta_k$.
\end{lemma}

\startproof
These operators are Glick's operator in disguise.
    When $n=3k+1$ we can relabel our $n$-gons in
    a way that converts $\Delta_k$ to the pentagram map.
    The corresponding space of
    birds $B_{n,k}$ corresponds to some strange set of ``relabeled $k$-birds''.
    This relabeling converts $G_{P,k,k}$ respectively to
    Glick's original operator.  This proves the invariance of
    $G_{P,k,k}$ under $\Delta_k$ when $n=3k+1$. A similar
    thing works for $n=3k+2$, but this time the relabeling
    converts $\Delta_k$ to the inverse of the pentagram map.
    \endproof

I was not able to find any similar formulas when $n>3k+2$.

\begin{question}
  When $n>3k+2$ and $P$ is a $k$-bird, are
  the coordinates of the collapse
  point $\widehat P_{\infty}$ algebraic functions of the
  vertices of $P$?
\end{question}

Here is one more thing I
have wondered about.  Suppose that $n$ is very large and
$P$ is a convex $n$-gon.  Then $P$ can be considered as
a $k$-bird for all $k=1,2,...,\beta$, where $\beta$ is the largest
integer such that $n \geq 3\beta+1$.    When we apply the
map $\Delta_k$ for these various values of $k$ we get
potentially $\beta$ different collapse points.  All I can
say, based on experiments, is that these points are not generally collinear.

\begin{question}
  Does the collection of $\beta$ collapse points in this
  situation have any special meaning?
  \end{question}

    \subsection{Star Relabelings}
    \label{relabel000}

    Let us further take up the theme in the proof of
    Lemma \ref{invariance}.  Given an $n$-gon $P$ and
    and some integer $r$ relatively prime to $n$, we define
    a new $n$-gon $P^{*r}$ by the formula
    \begin{equation}
      P^{*r}_j = P_{rj}.
    \end{equation}
    Figure 1.5 shows the $P^{*(-3)}$ when $P$ is the regular $10$-gon.
    
    As we have already mentioned, the action of $\Delta_1$ on the
    $P^{*(-k)}$ is the same as the action of $\Delta_k$ on $P$
    when $n=3k+1$.  So, when $n=3k+1$, the pentagram map
    has another nice invariant set (apart from the set of
    convex $n$-gons), namely
    $$B_{k,n}^{*(-k)}=\{P^{*(-k)}|\ P \in B_{k,n}\}.$$
    The action of the pentagram map on this set is geometrically nice.
    If we suitably star-relabel, we get star-shaped (and hence embedded) polygons.
    A similar thing works when $n=3k+2$.

\end{document}